\newif\ifarxived
\arxivedtrue
\ifx\Arxived\Kaba\else\arxivedtrue\fi
\newif\ifextended
\ifx\Extended\Kaba\else\extendedtrue\fi
\newif\ifprivate
\ifx\Private\Kaba\else\privatetrue\extendedtrue\fi
\newif\ifJapanese
\newif\iftesting
\ifextended\testingtrue\fi
\ifarxived
\extendedtrue
\testingfalse
\privatefalse
\fi
\documentclass[12pt]{article}
\ifextended
\else
\usepackage{CJKutf8}
\usepackage{setspace}
\fi
\ifarxived
\usepackage{CJKutf8}
\usepackage{setspace}
\fi
\ifarxived
\usepackage[bookmarks=true,backref=page,bookmarksnumbered=true, 
  bookmarkstype=toc,
  pdftitle={Generic Absoluteness Revisited},%
  pdfsubject={},%
  pdfkeywords={Laver-generic large cardinal, Maximality Principles, Recurrence Axioms, 
    generic absoluteness}]{hyperref}
\else
\ifextended
\usepackage[dvipdfmx,backref=page,bookmarks=true,bookmarksnumbered=true,
  bookmarkstype=toc,
]{hyperref}
\else
\usepackage[bookmarks=true,backref=page,bookmarksnumbered=true, 
  bookmarkstype=toc,
  pdftitle={Generic Absoluteness Revisited},%
  pdfsubject={},%
  pdfkeywords={Laver-generic large cardinal, Maximality Principles, Recurrence Axioms, 
    generic absoluteness}]{hyperref}
\fi
\fi
\ifarxived
\else
\ifextended
\usepackage[dvipdfmx]{pxjahyper}% for uplatex bookmarks in Japanese fonts
\fi
\fi
\ifarxived
\hypersetup{
  colorlinks=true,
  linkcolor=blue,
  filecolor=blue,      
  citecolor=cyan,
  urlcolor=cyan,
  pdftitle={Generic Absoluteness Revisited},%
  pdfauthor={Sakaé Fuchino, Takehiko Gappo, and Francesco Parente},%
  pdfsubject={},%
  pdfkeywords={Laver-generic large cardinal, Maximality Principles, Recurrence Axioms, 
    generic absoluteness}
}
\else
\ifextended
\hypersetup{
  colorlinks=true,
  linkcolor=blue,
  filecolor=blue,      
  citecolor=cyan,
  urlcolor=cyan,
  pdftitle={Generic Absoluteness Revisited},%
  pdfauthor={Sakaé Fuchino (渕野 昌), Takehiko Gappo (合浦岳彦), and Francesco Parente},%
  pdfsubject={},%
  pdfkeywords={Laver-generic large cardinal, Maximality Principles, Recurrence Axioms, 
    generic absoluteness}
}
\else
\hypersetup{
  colorlinks=true,
  linkcolor=black,
  filecolor=black,      
  citecolor=black,
  urlcolor=black,
  pdftitle={Generic Absoluteness Revisited},%
  pdfauthor={Sakaé Fuchino, Takehiko Gappo, and Francesco Parente},%
  pdfsubject={},%
  pdfkeywords={Laver-generic large cardinal, Maximality Principles, Recurrence Axioms, 
    generic absoluteness}
}
\fi
\fi
\ifarxived
\usepackage{graphicx, color}
\else
\usepackage[dvipdfmx]{graphicx, color}
\fi
\usepackage{amsmath, amssymb}
\usepackage{bbm}% \mathbbm{} (for lowercase alhpabets) mathbb{} (for uppercase alphabets)
\usepackage{dsfont}
\usepackage{bbold}
\usepackage{accents} %% used in utilde
\usepackage{marginnote}
\usepackage[many]{tcolorbox}
\usepackage{mathtools} %% underbracket
\usepackage[normalem]{ulem}
\usepackage{setspace}
\definecolor{darkelectricblue}{rgb}{0.33, 0.39, 0.52}
\definecolor{darkgreen}{rgb}{0.31, 0.47, 0.26}
\ifextended
\newcommand{\extendedcolor}{\color{darkelectricblue}}

\newcommand{\privatecolor}{\color{darkgreen}}
\newcommand{\darkred}{\color[rgb]{0.8,0.1,0.1}}

\newcommand{\It}{\it\darkred{}}
\else

\newcommand{\darkred}{}

\newcommand{\It}{\it{}}
\fi
\usepackage{theorem}
\newcommand{\bbd}[1]{{\mathbb{#1}}}
\theorembodyfont{\ifJapanese\rm\else\it\fi}
\newcount\minute	% Current minute within the hour
\newcount\hour		% Current hour (24-hour type)
\newcount\hourMins  % Temporary for taking hour modulo 60
\ifJapanese
\def\today%
{% Displays today's date and time as ``yyyy-mm-dd hh:mm''
  \the\year 年\,\zeroPadTwo{\the\month}月\,\zeroPadTwo{\the\day}日%
}% --- \today ---
\fi
\def\now%
{% Displays today's time as ``hh:mm''
  \minute=\time    % Number of minutes since midnight
  \hour=\time \divide \hour by 60 % Get hours
  \hourMins=\hour \multiply\hourMins by 60
  \advance\minute by -\hourMins % Hours modulo 60
  \zeroPadTwo{\the\hour}:\zeroPadTwo{\the\minute}%
}% --- \now ---
\def\zeroPadTwo#1%
{% Left zero pad of the argument to 2 digits.  The argument
  \ifnum #1<10 0\fi    % Conditionally output a zero
  #1%    Then output the argument
}% --- \zeroPadTwo ---
\title{\scalebox{1.1}[1.1]{\bf Generic Absoluteness Revisited}}
\author{\ifextended
\protect\scalebox{1}[1.2]{\quad Sakaé Fuchino$^{\ast,\ddagger}$\quad\ \  Takehiko 
  Gappo$^{\dagger}$\quad\ \ Francesco Parente$^{\ast,\S}$}\\
\else\ifarxived\scalebox{1}[1.2]{\quad Sakaé Fuchino$^{\ast,\ddagger}$\quad\ \  Takehiko 
  Gappo$^{\dagger}$\quad\ \ Francesco Parente$^{\ast,\S}$}\\\else
\protect\scalebox{1}[1.2]{Sakaé Fuchino$^{\ast,\ddagger}$,\ \ Takehiko 
  Gappo$^{\dagger}$, \ and Francesco Parente$^{\ast,\S}$\quad}\\
\fi
\fi
  \ifarxived%%
  \protect\scalebox{1}[1.2]{\begin{CJK}{UTF8}{goth}\qquad\quad 渕野 昌\qquad\qquad\quad\ 合浦 岳彦\end{CJK}\quad\ \ 
    \begin{CJK}{UTF8}{goth}パレンテ・フランチェスコ{}
    \end{CJK}}\medskip\\
  \else\ifextended\protect\scalebox{1}[1.2]{\qquad\quad 渕野 昌\qquad\qquad\quad\ 合浦 岳彦\quad\ \ 
    パレンテ・フランチェスコ}\medskip\\ 
  \else
  \phantom{\protect\scalebox{1}[1.2]{\begin{CJK}{UTF8}{goth}\qquad\quad 渕野 昌\qquad\qquad\quad\ 合浦 岳彦\end{CJK}}\quad\ \ 
    \phantom{\begin{CJK}{UTF8}{goth}パレンテ・フランチェスコ{}
    \end{CJK}}}\medskip\\
  \fi
  \fi
  }
\date{}
\ifarxived
\else
\ifextended
\fi
\fi
\renewcommand{\baselinestretch}{1.2}
\renewcommand{\thefootnote}{(\arabic{footnote})\,}
\iftesting
\fi%\iftesting
\iftesting
\newcommand{\Label}[1]{\label{#1}\marginnote{{\color{cyan}%
      \renewcommand{\baselinestretch}{0.4}\tiny 
		  \rlap{#1}}}}

\else
\newcommand{\Label}[1]{\label{#1}}

\fi%\iftesting

\def\memo#1{\ifprivate\marginnote{{\privatecolor\normalsize%
      \renewcommand{\baselinestretch}{0.4}\tiny\mbox{}\vspace{-2.52ex}% 1.4*1.8=2.52
      \par\relax%
			#1\par\mbox{}}}\else\fi}%
\def\memox#1{}
\def\imemox#1{}

\newcounter{frml}[section]
\newcounter{frmla}[section]
\def\thefrml{{\arabic{section}.\arabic{frml}}}
\def\thefrmla{{$\aleph$\arabic{section}.\arabic{frmla}}}
\def\frmlabel#1{\refstepcounter{frml}{\def\baka{#1}\ifx\baka\empty\else\label{#1}\fi}%
{\rm({\thefrml})\hfill\hfill\hfill}}
\def\ifrmlabel#1{\refstepcounter{frml}{\def\baka{#1}\ifx\baka\empty\else\label{#1}\fi}%
{\iftesting\darkred\fi\rm({\thefrml})\,:\hspace{0.6em}}}
\def\frmlabela#1{\refstepcounter{frmla}{\def\baka{#1}\ifx\baka\empty\else\label{#1}\fi}%
{\rm({\thefrmla})\hfill\hfill\hfill}}
\def\ifrmlabela#1{\refstepcounter{frmla}{\def\baka{#1}\ifx\baka\empty\else\label{#1}\fi}%
{\iftesting\darkred\fi\rm({\thefrmla})\,:\hspace{0.6em}}}
\def\xitem[#1]{\item[\frmlabel{#1}]\mbox{}%
	\iftesting\marginnote{{\renewcommand{%
				\baselinestretch}{0.6}\color{cyan}\tiny#1}}\fi\ignorespaces}
\def\xitemq[#1]{\item[\frmlabel{#1}]\mbox{}%
				\baselinestretch}{0.6}\tiny#1}}\fi\ignorespaces}
\def\xitemd[#1]#2{\item[{\rm(\ref{#1})$#2$}\hfill\hfill\hfill]}
\def\xitemA[#1]{\item[\frmlabela{#1}]\mbox{}%
	\iftesting\marginnote{{\renewcommand{%
				\baselinestretch}{0.6}\tiny#1}}\fi\ignorespaces}
\def\xitemx[#1]{\item[]}
\def\smashx#1{{#1}}
\def\xitemsub[#1]#2{\item[\frmlabel{#1}$_{#2}$]\mbox{}%
	\iftesting\marginnote{{\renewcommand{%
				\baselinestretch}{0.6}\tiny#1}}\fi\ignorespaces}
\def\xxitem[#1][#2]{\item[{\rm(\ref{#1}{\makebox[1.6ex][r]{#2}})}]\mbox{}%
	\iftesting\marginnote{{\renewcommand{%
				\baselinestretch}{0.6}\tiny\{#1\}\{#2\}}}\fi\ignorespaces}
\def\xitemof#1{{\rm({\ref{#1}})}}
\def\Xitem[#1]{\item[{\makebox[7ex][l]{\rm(\ref{#1})}}]\iftesting\marginnote{{\renewcommand{%
				\baselinestretch}{0.6}\tiny#1}}\fi\ignorespaces}

\def\xitemsubof#1#2{{\rm({\ref{#1}})${}_{#2}$}}
\def\xxitemof#1#2{(\ref{#1}\,#2)}

\newenvironment{xitemize}{\begin{list}{}{\parsep=0.5\smallskipamount%
			\itemindent=-0.4ex%
			\itemsep=0.5\smallskipamount\leftmargin=4em\labelwidth=3em\labelsep=0.7em}}%
							 {\end{list}}
\iftesting
\def\ixitem[#1]{\ifrmlabel{#1}\marginnote{{\color{cyan}\renewcommand{%
				\baselinestretch}{0.6}\qquad\qquad\tiny\rlap{#1}\mbox{}}}\ignorespaces}
\else
\def\ixitem[#1]{\ifrmlabel{#1}\ignorespaces}
\fi
\iftesting
\def\ixitemq[#1]{\frmlabel{#1}\marginnote{{\color{cyan}\renewcommand{% ixitem without colon
				\baselinestretch}{0.6}\qquad\qquad\tiny\rlap{#1}\mbox{}}}\ignorespaces}
\else
\def\ixitemq[#1]{\frmlabel{#1}\ignorespaces}
\fi
\iftesting
\def\ixitemx[#1]#2{\ifrmlabel{#1}\marginnote{{\color{cyan}\renewcommand{%
				\baselinestretch}{0.6}{}#2{}\qquad\qquad\tiny\rlap{#1}\mbox{}}}\ignorespaces}
\else
\def\ixitemx[#1]#2{\ifrmlabel{#1}\ignorespaces}
\fi
\iftesting
\def\ixitema[#1]{\ifrmlabela{#1}\marginnote{{\color{cyan}\renewcommand{%
				\baselinestretch}{0.6}\tiny\mbox{}\hfill #1}}\ignorespaces}
\else
\def\ixitema[#1]{\ifrmlabela{#1}\ignorespaces}
\fi
\def\assert#1{\noindent\makebox[4.8ex][r]{\rm(\makebox[2.2ex][c]{#1})}\ \ \ignorespaces}
\def\wassert#1{\assert{#1}}
\def\wassertof#1{\makebox[4.8ex][r]{\rm(\makebox[2.2ex][c]{#1})\ }}%
\def\assertof#1{{({\rm #1})}}%

\def\daimaru#1{\makebox[1em][c]{\mbox{\leavevmode\lower.144ex\hbox{% 0.08*1.8=0.14400000000000002
        \rlap{\hbox to 
          0.76em{\hfil\mbox{}\hfill{}\raisebox{0.054ex}{\scalebox{1.2}{○}}\hfil}}% 0.03*1.8=0.054
        \raise0.342ex\hbox to 1em{\hfil{\hspace{0.16em}\footnotesize#1}\hfil}}}}\,}% 0.19*1.8=0.342
\newtheorem{Thm}{\ifJapanese{\bf 定理}\else {\bf Theorem}\fi}[section]
\ifextended\tcolorboxenvironment{Thm}{
  colback=blue!4!white,
  boxrule=0pt,
  boxsep=0pt,
  left=8pt,right=8pt,top=2pt,bottom=4pt,
  oversize=4pt,
  sharp corners,
  before skip=\topsep,
  after skip=\topsep,
  breakable
}\fi

\newtheorem{ThmA}{\ifJapanese{\bf 定理\,A\!}\else{\bf Theorem\,A\!}\fi}[section]
\ifextended\tcolorboxenvironment{ThmA}{
  colback=blue!4!white,
  boxrule=0pt,
  boxsep=0pt,
  left=8pt,right=8pt,top=2pt,bottom=4pt,
  oversize=4pt,
  sharp corners,
  before skip=\topsep,
  after skip=\topsep,
  breakable
}\fi

\ifextended\tcolorboxenvironment{Ex}{
  colback=blue!3!white,
  boxrule=0pt,
  boxsep=0pt,
  left=8pt,right=8pt,top=2pt,bottom=4pt,
  oversize=4pt,
  sharp corners,
  before skip=\topsep,
  after skip=\topsep,
  breakable
}\fi
\newtheorem{Prop}[Thm]{\ifJapanese{\bf 命題}\else{\bf Proposition}\fi}
\ifextended\tcolorboxenvironment{Prop}{
  colback=blue!4!white,
  boxrule=0pt,
  boxsep=0pt,
  left=8pt,right=8pt,top=2pt,bottom=4pt,
  oversize=4pt,
  sharp corners,
  before skip=\topsep,
  after skip=\topsep,
  breakable
}\fi
\newtheorem{Problem}[Thm]{\ifJapanese{\bf 未解決問題}\else{\bf Problem}\fi}

\newtheorem{Lemma}[Thm]{\ifJapanese{\bf 補題}\else{\bf Lemma}\fi}
\ifextended\tcolorboxenvironment{Lemma}{
  colback=blue!4!white,
  boxrule=0pt,
  boxsep=0pt,
  left=8pt,right=8pt,top=2pt,bottom=4pt,
  oversize=4pt,
  sharp corners,
  before skip=\topsep,
  after skip=\topsep,
  breakable
}\fi
\newtheorem{LemmaA}[ThmA]{\ifJapanese{\bf 補題\,A\!}\else{\bf Lemma\,A\!}\fi}
\ifextended\tcolorboxenvironment{LemmaA}{
  colback=blue!4!white,
  boxrule=0pt,
  boxsep=0pt,
  left=8pt,right=8pt,top=2pt,bottom=4pt,
  oversize=4pt,
  sharp corners,
  before skip=\topsep,
  after skip=\topsep,
  breakable
}\fi

\ifextended\tcolorboxenvironment{CorA}{
  colback=blue!4!white,
  boxrule=0pt,
  boxsep=0pt,
  left=8pt,right=8pt,top=2pt,bottom=4pt,
  oversize=4pt,
  sharp corners,
  before skip=\topsep,
  after skip=\topsep,
  breakable
}\fi
\newtheorem{Cor}[Thm]{\ifJapanese{\bf 系}\else{\bf Corollary}\fi}
\ifextended\tcolorboxenvironment{Cor}{
  colback=blue!4!white,
  boxrule=0pt,
  boxsep=0pt,
  left=8pt,right=8pt,top=2pt,bottom=4pt,
  oversize=4pt,
  sharp corners,
  before skip=\topsep,
  after skip=\topsep,
  breakable
}\fi

\ifextended\tcolorboxenvironment{Remark}{
  colback=blue!4!white,
  boxrule=0pt,
  boxsep=0pt,
  left=8pt,right=8pt,top=2pt,bottom=4pt,
  oversize=4pt,
  sharp corners,
  before skip=\topsep,
  after skip=\topsep,
  breakable
}\fi

\newtheorem{Claim}{{\bf Claim}}[Thm]
\ifextended\tcolorboxenvironment{Claim}{
  colback=pink!5!white,
  boxrule=0pt,
  boxsep=0pt,
  left=8pt,right=8pt,top=2pt,bottom=4pt,
  oversize=4pt,
  sharp corners,
  before skip=\topsep,
  after skip=\topsep,
  breakable
}\fi
\ifextended\tcolorboxenvironment{ClaimA}{
  colback=pink!5!white,
  boxrule=0pt,
  boxsep=0pt,
  left=8pt,right=8pt,top=2pt,bottom=4pt,
  oversize=4pt,
  sharp corners,
  before skip=\topsep,
  after skip=\topsep,
  breakable
}\fi
\newcommand{\prf}{\noindent\ifJapanese{\bf 証明．\ }\ignorespaces\else{\bf 
		Proof.\ \ }\ignorespaces\fi}

\newcommand{\prfofClaim}{\raisebox{-.4ex}{\Large $\vdash$\ \ }}

\newcommand{\Thmof}[1]{\ifJapanese{定理\,\ref{#1}}\else{Theorem~\ref{#1}}\fi}
\newcommand{\ThmAof}[1]{\ifJapanese{定理\,A\,\ref{#1}}\else{Theorem\,A\,\ref{#1}}\fi}

\newcommand{\Lemmaof}[1]{\ifJapanese{補題\,\ref{#1}}\else{Lemma~\ref{#1}}\fi}

\newcommand{\Propof}[1]{\ifJapanese{命題\,\ref{#1}}\else{Proposition~\ref{#1}}\fi}

\newcommand{\Corof}[1]{\ifJapanese{系\,\ref{#1}}\else{Corollary~\ref{#1}}\fi}

\newcommand{\sectionof}[1]{\ifJapanese{第\ref{#1}節}\else{Section~\ref{#1}}\fi}
\newcommand{\subsectionof}[1]{\ifJapanese{第\ref{#1}分節}\else{Subsection~\ref{#1}}\fi}
\newcommand{\footnoteof}[1]{Footnote \ref{#1}}
\newcommand{\Thmabove}{{\ifJapanese 定理\else Theorem\fi\ \number\theThm}}

\newcommand{\Propabove}{{\ifJapanese 命題\else Proposition\fi\ \number\theThm}}
\newcommand{\Lemmaabove}{{\ifJapanese 補題\else Lemma\fi\ \number\theThm}}

\newcommand{\Claimabove}{{Claim \number\theClaim}}

\newcommand{\ubecause}[3]{\underbrace{{}#1{}%
  \ifx\bakakaba#2\bakakaba\rule[-0.72ex]{0pt}{1pt}\else\rule[#2]{0pt}{1pt}\fi}_{\mbox{\footnotesize\clap{#3}}}}
\newcommand{\obecause}[3]{\overbrace{{}#1{}%
  \ifx\bakakaba#2\bakakaba\rule[1.62ex]{0pt}{1pt}\else\rule[#2]{0pt}{1pt}\fi}^{\mbox{\footnotesize\clap{#3}}}}
\newsavebox{\qedbox}\sbox{\qedbox}{% QED symbol (by U.Fuchs modified by S.F.)
{\unitlength=0.05mm \begin{picture}(40,60)
\put(0,0){\framebox(30,44)[cc]{}}
\put(30,-7){\rule{7\unitlength}{44\unitlength}}
\put(10,-7){\rule{27\unitlength}{7\unitlength}}
\end{picture}}}
\newcommand{\qed}{\mbox{}\hfill\usebox{\qedbox}}
\newcommand{\smallqed}%
{\mbox{}\smallskip\hfill\raisebox{-.4ex}{\Large $\dashv$}}
\newcommand{\qedof}[1]%
{\mbox{} \hspace*{\fill}{\usebox{\qedbox}{\tiny~(#1)}}}
\newcommand{\Qedof}[1]%
{\mbox{} \hspace*{\fill}{\usebox{\qedbox}%
{\tiny~(#1~\number\theThm)}}}
\newcommand{\QedAof}[1]%
{\mbox{} \hspace*{\fill}{\usebox{\qedbox}%
{\tiny~(#1~\number\theThmA)}}}
\newcommand{\qedofThm}{\Qedof{\ifJapanese 定理\else Theorem\fi}}

\newcommand{\qedofCor}{\Qedof{\ifJapanese 系\else Corollary\fi}}
\newcommand{\qedofProp}{\Qedof{\ifJapanese 命題\else Proposition\fi}}
\newcommand{\qedofLemma}{\Qedof{\ifJapanese 補題\else Lemma\fi}}

\newcommand{\qedskip}{\medskip}

\newcommand{\qedofClaim}%
{\mbox{}\hfill\raisebox{-.4ex}{\Large $\dashv$ }\nolinebreak%
\mbox{\tiny~(Claim~\number\theClaim)}}
\newcommand{\qedofClaimA}%
{\mbox{}\hfill\raisebox{-.4ex}{\Large $\dashv$ }\nolinebreak%
\mbox{\tiny~(Claim~A\,\number\theClaimA)}}
\newcommand{\qedofClaimAof}[1]%
{\mbox{}\hfill\raisebox{-.4ex}{\Large $\dashv$ }\nolinebreak%
\mbox{\tiny~(Claim~A\,\ref{#1})}}
\newcommand{\qedofSubclaim}%
{\mbox{}\hfill\raisebox{-.4ex}{\Large $\dashv$ }\nolinebreak%
\mbox{\tiny~(Subclaim~\number\theSubclaim)}}
\newcommand{\qedofSubsubclaim}%
{\mbox{}\hfill\raisebox{-.4ex}{\Large $\dashv$ }\nolinebreak%
\mbox{\tiny~(Subsubclaim~\number\theSubsubclaim)}}
\newcommand{\cardof}[1]{\mathopen{|\,}#1\mathclose{\,|}}

\newcommand{\Card}{{\it Card\/}}

\newcommand{\setof}[2]{\{#1\,:\,#2\}}
\newcommand{\ssetof}[1]{\{#1\}}

\newcommand{\subseteqand}[1]{\mathrel{\mathop{\subseteq}%
		\limits_{\scriptscriptstyle\hbox to 14pt{$\scriptscriptstyle #1$\hss}}}}

\newcommand{\mapping}[3]{#1:#2\rightarrow #3}
\newcommand{\isomrph}[3]{#1:#2\stackrel{\cong\hspace{0.8ex}}{\rightarrow}#3}

\newcommand{\Elembed}[4]{#1:#2\stackrel{\prec\hspace{0.8ex}}{\rightarrow}_{#4}#3}

\newcommand{\fnsp}[2]{\mbox{}^{{#1}\hspace{-0.02em}}#2}
\newcommand{\imageof}{{}^{\,{\prime}{\prime}}}

\newcommand{\seqof}[2]{\langle#1\,:\,#2\rangle}
\newcommand{\pairof}[1]{\langle#1\rangle}
\newcommand{\rpairof}[1]{(#1)}
\newcommand{\psof}[1]{{\mathcal P}\/(#1)}
\newcommand{\bvalof}[2]{{[}\!{[}\,#1\,{]}\!{]}_{#2}}
\newcommand{\forces}[2]{\,\|\hspace{-.35ex}\mbox{\sf--}_{\,#1\,}%
\mbox{\rm``}\,#2\,\mbox{\rm''}}

\newcommand{\checked}{^{\checkmark}}
\newcommand{\modelof}[1]{\models\!\mbox{\rm``\,}#1\mbox{\rm\,''}}

\newcommand{\bbone}{{\mathord{\mathbb{1}}}}
\newcommand{\bbzero}{{\mathord{\mathbb{0}}}}
\newcommand{\circleq}{\mathrel{{\leqslant}%
		\hspace{-0.86ex}{\lower-0.53ex\hbox{$\scriptscriptstyle\circ$}}}}
\newcommand{\symb}[1]{{\mathord{\hspace{0.08em}\underbracket[0.6pt][2pt]{#1}}\hspace{0.08em}}}
\newcommand{\vsymb}[1]{{\underline{#1}}}

\newcommand{\ol}[1]{\overline{#1}}
\newcommand{\ul}[1]{\underline{#1}}
\newcommand{\restr}{\restriction}
\newcommand{\id}{\mathop{id\,}}

\newcommand{\Col}{{\rm Col}}

\newcommand{\trcl}{\mathop{\mbox{\it trcl\/}}}

\newcommand{\natnums}{{\bbd{N}}}

\newcommand{\clpoP}{\mathfrak{P}}
\newcommand{\clgenG}{\mathfrak{G}}
\newcommand{\poC}{\bbd{C}}
\newcommand{\poD}{\bbd{D}}
\newcommand{\poP}{\bbd{P}}
\newcommand{\poQ}{\bbd{Q}}

\newcommand{\poS}{\bbd{S}}
\newcommand{\On}{{\rm On}}

\newcommand{\BaB}{\bbd{B}}

\newcommand{\genG}{\mathbb{G}}
\newcommand{\geng}{\mathbb{g}}

\newcommand{\utpoQ}{\utilde{\mathbb{Q}}}

\newcommand{\utgenG}{\utilde{\mathbb{G}}}

\newcommand{\genH}{\mathbb{H}}

\newcommand{\condp}{\mathbbm{p}}
\newcommand{\condq}{\mathbbm{q}}

\newcommand{\boundingno}{{\mathfrak b}}
\newcommand{\dominatingno}{{\mathfrak d}}

\newcommand{\LT}{{<}\,}
\newcommand{\LE}{{\leq}\,}
\newcommand{\GT}{{>}\,}
\newcommand{\GE}{{\geq}\,}

\newcommand{\ctentenc}{,{}\linebreak[0]\hspace{0.04ex}{{.}{.}{.}\hspace{0.1ex},\,}\linebreak[0]}

\newcommand{\xmbox}[1]{ $\relax{\rm #1}\relax$ }

\newcommand{\gmM}{\mathfrak{M}}

\newcommand{\continuum}{2^{\aleph_0}}
\newcommand{\calA}{{\mathcal A}}
\newcommand{\calB}{{\mathcal B}}

\newcommand{\calD}{{\mathcal D}}

\newcommand{\calH}{{\mathcal H}}

\newcommand{\calK}{{\mathcal K}}
\newcommand{\calL}{{\mathcal L}}
\newcommand{\calM}{{\mathcal M}}

\newcommand{\calP}{{\mathcal P}}

\newcommand{\calR}{{\mathcal R}}
\newcommand{\calS}{{\mathcal S}}

\newcommand{\calU}{{\mathcal U}}

\newcommand{\calX}{{\mathcal X}}

\newcommand{\utab}{\utilde{\ol{a}}} %% undertilde a bar
\newcommand{\utb}{\utilde{b}}

\newcommand{\utf}{\utilde{f}}
\newcommand{\utg}{\utilde{g}}

\newcommand{\utE}{\utilde{E}}
\newcommand{\utI}{\utilde{I}}

\newcommand{\utS}{\utilde{S}}
\newcommand{\utgmM}{\utilde{\gmM}}

\newcommand{\Lin}{{\calL}_{\in}}

\newcommand{\ZF}{{\sf ZF}}
\newcommand{\ZFC}{{\sf ZFC}}
\newcommand{\CH}{{\sf CH}}

\newcommand{\MA}{{\sf MA}}
\newcommand{\FA}{{\sf FA}}
\newcommand{\MM}{{\sf MM}}
\newcommand{\MMpp}{\mbox{\sf MM$^{++}$}}
\newcommand{\MP}{{\sf MP}}

\newcommand{\RA}{{\sf RA}}
\newcommand{\RcA}{{\sf RcA}}

\newcommand{\RcAp}{{\sf RcA}{$^+$}}

\newcommand{\PFA}{{\sf PFA}}
\newcommand{\BMM}{{\sf BMM}}
\newcommand{\BFA}{{\sf BFA}}
\newcommand{\NS}{{\sf NS}}
\newcommand{\INS}{I_\NS}

\newcommand{\CCA}{{\sf CCA}}
\newcommand{\GA}{{\sf GA}}

\newcommand{\refl}{{\mathfrak{r}\mathfrak{e}\mathfrak{f}\mathfrak{l}\,}}

\newcommand{\st}{such that}
\newcommand{\wrt}{with respect to}
\newcommand{\Wolog}{Without loss of generality}

\newcommand{\tfae}{the following are equivalent}

\newcommand{\po}{poset}
\newcommand{\pos}{posets}
\newcommand{\uniV}{\mathsf{V}}
\newcommand{\uniL}{\mathsf{L}}
\newcommand{\uniW}{\mathsf{W}}

\newcommand{\Pkl}[2]{\ifx\bakakaba#1\bakakaba\ifx\bakakaba#2\bakakaba{\mathcal 
    P}_\kappa(\lambda)\else{\mathcal P}_\kappa(#2)\fi\else{\mathcal P}_{#1}(#2)\fi}

\newcommand{\utildeT}[1]{%
  \underaccent{{\sim}}{#1}}
\newcommand{\utildeS}[1]{%
	\hbox to 0pt{\smash{$\mathop{\scriptstyle #1}\limits_{%
				\raisebox{0.6ex}[0pt]{$\scriptscriptstyle\sim$}}$}\hss}%
	\relax\phantom{\mathord{{#1}_{\rule[-0.6ex]{0pt}{1pt}}}}}
\newcommand{\utildeSS}[1]{%
	\hbox to 0pt{$\mathop{\scriptscriptstyle #1}%
		\limits_{\scriptscriptstyle\sim}$\hss}%
		\relax\phantom{\underline{#1}}}
\newcommand{\utilde}[1]{%
	\mathchoice{\utildeT{#1}}{\utildeT{#1}}{\utildeS{#1}}{\utildeSS{#1}}}

{\end{minipage}\end{trivlist}}

\newcommand{\baB}{\bbd{B}}

\begin{document}
\maketitle
\renewcommand{\thefootnote}{$\ast$\ }
  \footnotetext{Graduate School of System Informatics, Kobe University \\Rokko-dai 1-1, Nada, Kobe 657-8501 Japan
   \\
   \quad\scalebox{0.95}[1]{\tt $^\ddagger$ fuchino@diamond.kobe-u.ac.jp}\quad
  \quad\scalebox{0.95}[1]{\tt $^\S$ francesco.parente@people.kobe-u.ac.jp}}
\renewcommand{\thefootnote}{$\dagger$\ }
  \footnotetext{Technische Universität Wien
Institut für Diskrete Mathematik und Geometrie\\
Wiedner Hauptstraße 8-10/104
1040 Wien, Austria\\
    \quad\scalebox{0.95}[1]{\tt takehiko.gappo@tuwien.ac.at}}
\ifextended
\phantomsection
\addcontentsline{toc}{section}{* Generic Absoluteness Revisited} 
\ifarxived
\addcontentsline{toc}{section}{**** by S.Fuchino, T.Gappo, and F.Parente}
\else
\addcontentsline{toc}{section}{**** by S.Fuchino (渕野昌), T.Gappo (合浦岳彦), and 
  F.Parente (フランチェスコ・パレンテ)}
\fi
\fi

\ifextended
\addcontentsline{toc}{section}{Abstract}
\fi
\begin{abstract}
  The present paper is concerned with the relation between recurrence axioms and 
  Laver-generic large cardinal axioms in light of principles of generic absoluteness and the 
  Ground Axiom. 

  M.\,Viale %% \cite{viale-revisited}
  proved that Martin's Maximum$^{++}$ together with the 
  assumption that there are class many Woodin cardinals implies
  $\calH(\aleph_2)^\uniV\prec_{\Sigma_2}\calH(\aleph_2)^{\uniV[\genG]}$ for a 
  generic $\genG$ on any stationary 
  preserving $\poP$ which also preserves Bounded Martin's Maximum. We show that a
  similar but more general conclusion follows from each 
  of $(\calP,\calH(\kappa))_{\Sigma_2}$-\RcAp\ (which is a fragment of a  
  reformulation of the  
  Maximality Principle for $\calP$ and $\calH(\kappa)$), and the existence of the tightly
  $\calP$-Laver-generically huge cardinal.\renewcommand{\thefootnote}{$\star$\ }\footnotemark

  While under ``$\calP=$ all stationary preserving \pos'', our results are not very much more than 
  Viale's Theorem%% (even so they are proved without the stationary tower forcing 
  %% technique, and in case of $(\calP,\calH(\kappa))_{\Sigma_2}$-\RcAp, with much less 
  %% consistency strength)
  , 
  for other classes of \pos, ``$\calP=$ all proper \pos'' or ``$\calP=$ all ccc \pos'', for 
  example, our theorems are 
  not at all covered by his theorem. 

  The assumptions (and hence also the conclusion) of Viale's Theorem are compatible with the 
  Ground Axiom. In contrast, we show that the assumptions of our theorems (for most of the 
  common settings of $\calP$ 
  and with a modification of the large cardinal property involved) imply the negation of the 
  Ground Axiom. This fact is used to show that fragments of Recurrence Axiom
  $(\calP,\calH(\kappa))_\Gamma$-\RcAp\ can be different from the corresponding fragments of 
  Maximality Principle $\MP(\calP,\calH(\kappa))_\Gamma$ for $\Gamma=\Pi_2$. 
\end{abstract}
\renewcommand{\thefootnote}{$\star$\ }\footnotetext{We use here the definite article since 
  it is known that a  
  tightly $\calP$-Laver-generic large cardinal, if it exists, is the unique cardinal $\kappa_\refl$
  ($=\max(\ssetof{\aleph_2,\continuum})$) for (almost?) all reasonable non-trivial 
  instances of $\calP$ and notions of large cardinal, see \cite{sfetal-II}. }

\iftrue
{%%\extendedcolor 
\newpage
\phantomsection
\addcontentsline{toc}{section}{Contents}
\newcommand{\myscalebox}[1]{\scalebox{0.88}[1.06]{#1}}
%% contents 目次 index
\begin{quotation}
	\footnotesize
	\noindent
	\centerline{%%
      \normalsize\tt\quad\ Contents\hspace{6em}\mbox{}}\mbox{}\\
  %% {\mbox{}\hspace{-1.6em}\tt\makebox[3.4ex][l]{\ref{*}.}%
  %%  \hyperref[*]{\tt\myscalebox{*}}}\ \ \dotfill\ \ {\pageref{*}}\\  
  {\mbox{}\hspace{-1.6em}\tt\makebox[3.4ex][l]{\ref{theintro}.}%
   \hyperref[theintro]{\tt\myscalebox{Introduction and preliminaries}}}\ \ \dotfill\ \ {\pageref{theintro}}\\ 
  {\mbox{}\hspace{-1.6em}\tt\makebox[3.4ex][l]{\ref{rec-GA}.}%
   \hyperref[rec-GA]{\tt\myscalebox{Recurrence Axioms and the Ground Axiom}}}\ \ \dotfill\ \ {\pageref{rec-GA}}\\  
  \iffalse
  {\mbox{}\hspace{-0.6em}\tt\makebox[5ex][l]{\ref{*}}%
   \hyperref[*]{\tt\myscalebox{*}}}\ \ \dotfill\ \ {\pageref{*}}\\  
  %%
  \fi
  %%
  {\mbox{}\hspace{-0.6em}\tt\makebox[5ex][l]{\ref{hierarchies-rec-max}}%
   \hyperref[hierarchies-rec-max]{\tt\myscalebox{Hierarchies of Recurrence and Maximality}}}\ \ \dotfill\ \ {\pageref{hierarchies-rec-max}}\\  
  {\mbox{}\hspace{-0.6em}\tt\makebox[5ex][l]{\ref{ground-axiom}}%
   \hyperref[ground-axiom]{\tt\myscalebox{(In)compatibility of Recurrence and Maximality with Ground Axiom}}}\ \ \dotfill\ \ {\pageref{ground-axiom}}\\  
  {\mbox{}\hspace{-0.6em}\tt\makebox[5ex][l]{\ref{Laver-GA}}%
   \hyperref[Laver-GA]{\tt\myscalebox{Incompatibility of Laver genericity with Ground Axiom}}}\ \ \dotfill\ \ {\pageref{Laver-GA}}\\  
  {\mbox{}\hspace{-1.6em}\tt\makebox[3.4ex][l]{\ref{hierarchies}.}%
   \hyperref[hierarchies]{\tt\myscalebox{Hierarchies of restricted Recurrence Axioms and Maximality Principles}}}\ \ \dotfill\ \ {\pageref{hierarchies}}\\  
  {\mbox{}\hspace{-1.6em}\tt\makebox[3.4ex][l]{\ref{genabs-rec}.}%
   \hyperref[genabs-rec]{\tt\myscalebox{Generic absoluteness under Recurrence Axioms}}}\ \ \dotfill\ \ {\pageref{genabs-rec}}\\  
  {\mbox{}\hspace{-1.6em}\tt\makebox[3.4ex][l]{\ref{genabs-Laver}.}%
   \hyperref[genabs-Laver]{\tt\myscalebox{Generic absoluteness under Laver-genericity}}}\ \ \dotfill\ \ {\pageref{genabs-Laver}}\\  
  {\mbox{}\hspace{-1.6em}\tt\makebox[3.4ex][l]{\ref{misc}.}%
   \hyperref[misc]{\tt\myscalebox{Some more remarks and open questions}}}\ \ \dotfill\ \ {\pageref{misc}}\\  
  {\mbox{}\hspace{-0.6em}\tt\makebox[5ex][l]{\ref{RRAL}}%
   \hyperref[RRAL]{\tt\myscalebox{Restricted Recurrence Axioms under Laver-genericity}}}\ \ \dotfill\ \ {\pageref{RRAL}}\\  
  {\mbox{}\hspace{-0.6em}\tt\makebox[5ex][l]{\ref{Sep}}%
   \hyperref[Sep]{\tt\myscalebox{Separation of some other axioms and assertions}}}\ \ \dotfill\ \ {\pageref{Sep}}\\  
  {\ifextended
  {\mbox{}\hspace{-0.6em}\tt\makebox[5ex][l]{\ref{Yah}}%
   \hyperref[Yah]{\tt\myscalebox{\extendedcolor Yet another hierarchy of restricted Recurrence 
       Axioms}}}\ \ \dotfill\ \ {\pageref{Yah}}\\ \fi}
  {\ifprivate\privatecolor
  {\mbox{}\hspace{-1.6em}\tt\makebox[3.4ex][l]{\ref{class-many}.}%
   \hyperref[class-many]{\tt\myscalebox{\privatecolor Universe with class many large cardinals}}}\ \ \dotfill\ \ {\pageref{class-many}}\\  \fi}
  {\mbox{}\hspace{-1.6em}\hyperref[ref]{\tt References}}\ \ \dotfill\ \ 
  {\pageref{ref}}\medskip\\ 
\end{quotation}}
\fi%\ifextended
\renewcommand{\thefootnote}{}
\footnotetext{{\it Date:} June 2, 2024%% \ifextended
  \qquad {\it Last update:} 
  \today\ (\now\ \ifarxived UTC\else JST\fi)\vspace{-1\smallskipamount}%% \fi
}
\footnotetext{{\it MSC2020 Mathematical Subject Classification:} 03E45, 03E50, 03E55, 03E57, 03E65
  	\vspace{-1\smallskipamount}}
%%  03E50  	Continuum hypothesis and Martin's axiom
%% 	03E55  	Large cardinals
%%  03E57  	Generic absoluteness and forcing axioms
%%  03E65  	Other set-theoretic hypotheses and axioms

\footnotetext{{\it Keywords: generic large cardinal, Laver-generic large cardinal, 
    Maximality Principles, Continuum Hypothesis, Ground Axiom}
  }
%% 1.4*1.8=2.52
%%\footnotetext{\phantom{{\it Keywords: {}}}%% internal stationarity, 
  %% internal clubness
%% \footnotetext{A part of this work was done 
%%   during my stay in Bellaterra while I was a visiting researcher at the CRM 
%%   Barcelona joining the research programme Large Cardinals and Strong Logics. I 
%%   would like to thank Professor Joan Bagaria for the arrangement of the stay and 
%%   the CRM for its support and hospitality. }

\ifextended
\ifprivate
\footnotetext{\hspace{-1em}This is a private extended version of the paper in preparation.
  %% downloadable as:\smallskip\\  
%% \mbox{}\hfill{\tt http://fuchino.ddo.jp/papers/RIMS16-uncountable-reflection.pdf}\smallskip
  \par All additional 
  details not contained in the submitted version of the paper are either typeset in 
  {\extendedcolor ``dark electric blue''} in case the text is also included in the 
  (non-private) extended 
  version of the paper, or in {\privatecolor ``dark green''} if the comment are thought only for the 
  private version. 
  \iffalse
  \sout{The numbering of the assertions is kept identical with the submitted version.}
  Since the changes from the submitted version are now quite extensive, I am not trying to keep the 
  numbering of the theorems and assertions identical with the numbering in the submitted 
  version. \else 
  The numbering of the assertions is kept identical with the submitted version.
  \fi

  The most up-to-date pdf-file of this private extended version is downloadable as:\\
  \ifprivate\href{https://fuchino.ddo.jp/papers/generic-absoluteness-revisited-xx.pdf}{\tt 
  https://fuchino.ddo.jp/papers/generic-absoluteness-revisited-xx.pdf}\else
  \href{https://fuchino.ddo.jp/papers/generic-absoluteness-revisited-x.pdf}{\tt 
  https://fuchino.ddo.jp/papers/generic-absoluteness-revisited-x.pdf}\fi 
}
\else
\footnotetext{\extendedcolor This is an extended version of the paper with the 
  title ``Generic Absoluteness Revisited'' in preparation. 
  %% downloadable as:\smallskip\\  
%% \mbox{}\hfill{\tt http://fuchino.ddo.jp/papers/RIMS16-uncountable-reflection.pdf}\smallskip
  \par\extendedcolor  All additional 
  details not contained in the submitted version of the paper are either typeset in 
  dark electric blue (the color in which this paragraph is typeset) or put in a separate appendices. 
 %% \sout{The numbering of the assertions is kept identical with the submitted version.}
 %%  Since the changes from the submitted version are now quite extensive, I am not trying 
 %%  to keep the  
 %%  numbering of the theorems and assertions identical with the numbering in the submitted version. 

  The most up-to-date pdf-file of this extended version is downloadable as:\medskip\\
  \qquad\qquad\href{https://fuchino.ddo.jp/papers/generic-absoluteness-revisited-x.pdf}{\tt 
  https://fuchino.ddo.jp/papers/generic-absoluteness-revisited-x.pdf} \medskip\\
  %% The materials in this extended version may be reused in the forthcoming \cite{future}.
}
\fi%%\ifprivate ... \else
\else
\footnotetext{A pdf file of an updated and extended version of this paper 
  possibly with more details and 
  proofs is downloadable as:\\ \quad \href{https://fuchino.ddo.jp/papers/generic-absoluteness-revisited-x.pdf}{\tt https://fuchino.ddo.jp/papers/generic-absoluteness-revisited-x.pdf}}
\fi

\renewcommand{\thefootnote}{\arabic{footnote})\,}
\setcounter{footnote}{0}

\section{Introduction and preliminaries}
\Label{theintro}
In the following, we tried very hard to make the present paper as self-contained as 
possible. For notions and notation which remain unexplained, the reader may refer to 
\cite{millennium-book}, \cite{higher-inf}, or \cite{kunen-2011}. Set-theoretic forcing is 
treated here just as in \cite{kunen-2011} with the exception that $\poP$-names for a \po\ $\poP$ 
are represented with an under-tilde, e.g.\ as $\utpoQ$ or $\utilde{S}$. We adopt the (fake
but consistently interpretable) narration that generic filters ``exist'' though otherwise we
remain in the \ZFC\  
narrative so that all classes mentioned here are (meta-mathematically) definable classes. 

The main theorem of 
Viale \cite{viale-revisited} states:
\begin{Thm}\Label{p-theintro-a}{\rm(M.\ Viale, Theorem 1.4 in \cite{viale-revisited})}
  Assume that\/ \MMpp\ holds, and there are class many Woodin cardinals. Then, 
  for any 
  stationary preserving \po\ $\poP$ with $\forces{\poP}{\BMM}$, we have 
  \begin{xitemize}
  \item[]  
    $\calH(\aleph_2)^\uniV\prec_{\Sigma_2}\calH(\aleph_2)^{\uniV[\genG]}$\quad 
    for $(\uniV,\poP)$-generic $\genG$. \qed
  \end{xitemize}
\end{Thm}\memo{Example showing $\Sigma_2$ is optimal?}
Here \MMpp\ is the following strengthening of the Martin's Axiom (\MM):
\begin{xitemize}
\item[{\darkred(\MMpp): }] For any stationary preserving $\poP$, any family $\calD$ of dense 
  subsets  
  of $\poP$ with $\cardof{\calD}<\aleph_2$, and any family $\calS$ of $\poP$-names of 
  stationary subsets of $\omega_1$ with $\cardof{\calS}<\aleph_2$, there is 
  a $\calD$-generic filter $\genG$ on $\poP$ \st\ $\utilde{S}[\genG]$ is a stationary subset 
  of $\omega_1$ for all $\utilde{S}\in\calS$. 
\end{xitemize}

\BMM\ is the Bounded Martin's Maximum, a weakening of \MM\ which is an instance of 
Bounded Forcing Axioms: for a class $\calP$ of \pos\ closed under forcing equivalence and a 
cardinal $\kappa$, the 
{\It Bounded Forcing Axiom} for $\calP$ and $\LT\kappa$ is the axiom stating:
\begin{xitemize}
\item[{\darkred($\BFA_{\LT\kappa}(\calP)$): }] For any complete Boolean\footnotemark\ 
  $\poP\in\calP$, and a family $\calD$ of maximal antichains in $\poP$ \st\ 
$\cardof{\calD}<\kappa$ and $\cardof{I}<\kappa$ for all $I\in\calD$, there is 
  a $\calD$-generic filter $\genG$ on $\poP$. 
\end{xitemize}
\footnotetext{We say that a \po\ $\poP$ is {\It complete Boolean} if
  $\poP=\baB\setminus\ssetof{\bbzero_\baB}$ for a complete Boolean algebra. Note that the 
  definition of $\BFA_{\LT\kappa}(\calP)$ makes sense only when $\poP$ is complete Boolean 
  (since otherwise it can be the case that $\poP$ does not have any maximal antichains of 
  size $\LT\kappa$).}

The {\It Bounded Martin's Maximum} ({\darkred\BMM}) is
$\BFA_{\LT\aleph_2}(\mbox{stationary pres.\ \pos})$.  
Bounded Forcing Axioms were introduced by Goldstern and Shelah \cite{goldstern-shelah} 
answering a problem asked by the first author of the present paper in \cite{potential}.

In \Thmof{p-theintro-a}, the condition ``$\forces{\poP}{\BMM}$'' cannot be simply dropped. 
For example, the formula  
saying that  
there is a set which is the power set of $\omega$ is $\Sigma_2$. Since $\neg\CH$ holds 
in $\uniV$ under \MM, if $\poP$ forces \CH\ then
$\calH(\aleph_2)^\uniV\not\prec_{\Sigma_2}\calH(\aleph_2)^{\uniV[\genG]}$ 
for $(\uniV,\poP)$-generic $\genG$.

We show that a conclusion similar to and more general than that of \Thmof{p-theintro-a} 
follows from each of $(\calP,\calH(\kappa))_{\Sigma_2}$-\RcA$^+$ which  
is a fragment of Recurrence Axiom (a reformulation of Maximality Principle introduced in 
\cite{recurrence}) for $\calP$ and $\calH(\kappa)$, 
see \sectionof{rec-GA} below, and the existence of the tightly 
$\calP$-Laver-gen.\ large cardinal (\Thmof{p-genabs-rec-0}, \Thmof{p-genabs-Laver-1}).

The notion of Laver-generic large cardinal is introduced in Fuchino, Ottenbreit Maschio 
Rodrigues, and Sakai \cite{sfetal-II}. The definition we give here is the slightly modified 
version in later papers such as in Fuchino \cite{janos}:

For an iterable class $\calP$ of \pos\ (i.e.\ class $\calP$ of \pos\ satisfying 
\xitemof{x-theintro-0-0-a-0} and \xitemof{x-theintro-0-0-a-1} below) a cardinal $\kappa$ is 
said to be ({\It tightly, resp.}) {\It$\calP$-Laver-gen.\ supercompact}\/ if, for any
$\lambda>\kappa$ and $\poP\in\calP$, there is a $\poP$-name 
$\utpoQ$ with 
$\forces{\poP}{\utpoQ\in\calP}$, \st\  
for $(\uniV,\poP\ast\utpoQ)$-generic $\genH$, there  
are $j, M\subseteq\uniV[\genH]$ 
\st\ 
\smash{$\Elembed{j}{\uniV}{M}{\kappa}$},\footnotemark\ $j(\kappa)>\lambda$,
%% $j\imageof{j(\kappa)}$, 
$\poP,\genH, j\imageof\lambda\in M$ (and $\poP\ast\utpoQ$ is of size $\leq j(\kappa)$, 
resp.).%%\footnotemark  
\footnotetext{``$\Elembed{j}{\uniV}{M}{\kappa}$'' denotes the condition that $j$ is an elementary 
  embedding of $\uniV$ into a transitive $M$ with the critical point $\kappa$.}

This definition can be adopted to many other large cardinal notions other than supercompactness. 
The reader may refer to \cite{janos} for definitions of other variants of Laver-generic 
large cardinal.  
Defined as above, it is not obvious at first glance that the Laver-genericity is formalizable 
in the language $\Lin$ of \ZFC. That it is actually the case, is shown in Fuchino and Sakai 
\cite{fuchino-sakai}. 

A tightly $\calP$-Laver-generic large cardinal, if it exists, is unique and decided to be
$\kappa_\refl:=\max(\ssetof{\aleph_2,\continuum})$ for all known reasonable non-trivial 
instances of $\calP$ with a strong enough large cardinal notion (see \cite{sfetal-II}, or 
\cite{janos}, \cite{future}). This is the reason why we often simply talk about {\it the} tightly
$\calP$-Laver-generic large cardinal. 

While under ``$\calP=$ all semi-proper \pos'', our results are not much more than 
 slight variants of Viale's (but without relying on the stationary tower forcing 
technique),  
for other classes of \pos, for 
example ``$\calP=$ all proper \pos'' or ``$\calP=$ all ccc \pos'', they are not at all 
covered by Viale's result in \cite{viale-revisited} nor by its proof.

In the following we shall always 
assume that the classes $\calP$ of \pos\ we consider are {\It normal}, that is,
\begin{xitemize}
\xitem[x-theintro-0-0] 
  $\calP$ is 
  closed \wrt\ forcing equivalence, and $\ssetof{\bbone}\in\calP$. 
\end{xitemize}
In particular, we assume that for any $\poP_0\in\calP$ there is a complete 
Boolean\addtocounter{footnote}{-1}\footnotemark\ 
$\poP\in\calP$ which is forcing equivalent to $\poP_0$. In some cases like the 
case ``$\poP=$ all $\sigma$-closed \pos'' where the original class of \pos\ is not normal 
we just replace $\calP$ with its closure \wrt\ forcing equivalence without mention. 

In some cases (like in the definition of Laver-genericity above) it is natural to consider 
(normal) classes of \pos\ which are closed \wrt\ two-step 
iteration. A class $\calP$ of \pos\ is called {\It iterable} if 
\begin{xitemize}
\xitem[x-theintro-0-0-a-0] $\calP$ is closed \wrt\ restriction. That is, for $\poP\in\calP$ 
  and $\condp\in\poP$, we always have $\poP\restr\condp\in\calP$, and 
\xitem[x-theintro-0-0-a-1] For any $\poP\in\calP$, and any $\poP$-name $\utpoQ$ of a \po\ with
  $\forces{\poP}{\utpoQ\in\calP}$, we have $\poP\ast\utpoQ\in\calP$.
\end{xitemize}

Viale's Absoluteness \Thmof{p-theintro-a} is a result built upon the following \Thmof{p-theintro-0}.
We shall use the following notation for the formulation of the Theorem: For an ordinal
$\alpha$, let ${\darkred\alpha^{(+)}}:=\sup(\setof{\cardof{\beta}^+}{\beta<\alpha})$. Note that 
$\alpha^{(+)}=\alpha$ if $\alpha$ is a cardinal. Otherwise, we have
$\alpha^{(+)}=\cardof{\alpha}^+$. 

In Bagaria \cite{bagaria-bounded} the following theorem contains the extra assumption 
that, translated into the context of the following 
formulation, $\kappa$ is a successor of a cardinal of uncountable cofinality. However we can 
eliminate this assumption 
by slightly modifying the proof in \cite{bagaria-bounded}.

\begin{Thm}{\rm (Bagaria's Absoluteness Theorem, Theorem 5 in \cite{bagaria-bounded})}
  \Label{p-theintro-0}
  For an uncountable cardinal $\kappa$ and a class $\calP$ of \pos\ closed under forcing 
  equivalence, and restriction (in the sense of \xitemof{x-theintro-0-0-a-0}) \tfae: 
  \quad\wassertof{a} $\BFA_{\LT\kappa}(\calP)$.\smallskip

  \wassert{b} For any $\poP\in\calP$, $\Sigma_1$-formula $\varphi$ in $\Lin$ and
  $a\in\calH(\kappa)$, $\forces{\poP}{\varphi(a)}$ $\Leftrightarrow$ $\varphi(a)$.\smallskip

  \wassert{c} For any $\poP\in\calP$ and $(\uniV,\poP)$-generic $\genG$, we have 
  $\calH(\kappa)^\uniV\prec_{\Sigma_1}\calH((\kappa^{(+)})^{\uniV[\genG]})^{\uniV[\genG]}$. 
\end{Thm}
\prf
Note that \assertof{b} $\Leftrightarrow$ \assertof{c} is trivial since \ZFC\ proves that 
\begin{xitemize}
\xitem[x-theintro-0] 
  $\calH(\mu)\prec_{\Sigma_1}\uniV$ for any uncountable cardinal $\mu$
\end{xitemize}
(Lévy \cite{levy}%% see p.299 in Kanamori \cite{higher-inf}
).  
{\ifextended\extendedcolor\noindent [\![ 
      For $a\in\calH(\mu)$, let $\nu:=\cardof{\trcl^+(a)}<\mu$ 
      (here, $\trcl^+(a)$ denotes the variant of transitive closure which satisfies $a\in\trcl^+(a)$). 
      Then $\nu<\mu$. If $\calH(\mu)\models\varphi(a)$ for a $\Sigma_1$-formula $\varphi$ 
      then it follows $\uniV\models\varphi(a)$.

      Suppose that $\varphi(x)=\exists y\,\psi(x,y)$ where $\psi(x,y)$ 
      is a $\Sigma_0$-formula, and assume that $\uniV\models\varphi(a)$.  Let $b$ be \st\
      $\uniV\models\psi(a,b)$. Let $\delta$ be large enough \st\ $a$, $b\in V_\delta$. Let
      $M\prec V_\delta$ be \st\ $\trcl^+(a)\subseteq M$, $b\in M$, and $\cardof{M}=\nu$.
      Let $\isomrph{m}{M}{M_0}$ be the Mostowski-collapse. Note that
      $m\restr\trcl^+(a)=\id_{\trcl^+(a)}$. Thus $M_0\models\psi(a,m(b))$, and hence 
      $M_0\models\phi(a)$. Since $M_0\subseteq\calH(\mu)$, it follows that
      $\calH(\mu)\models\varphi(a)$. ]\!]\qedskip
  \fi}

Note also that if $\kappa=\continuum$, we also have the equivalence of \assertof{a}, \assertof{b},
\assertof{c} with
\begin{xitemize}
\item[\wassertof{c'}] {\it For any $\poP\in\calP$ and $(\uniV,\poP)$-generic $\genG$, we have
  $\calH(\continuum)^\uniV\prec_{\Sigma_1}\calH((\continuum)^{\uniV[\genG]})^{\uniV[\genG]}$.}
\end{xitemize}

{\ifextended{\extendedcolor \noindent By the remark above, it is 
  enough to prove \assertof{a} $\Leftrightarrow$ \assertof{b}.

  \assertof{a} $\Rightarrow$ \assertof{b}: Let $\poP\in\calP$. \Wolog, $\poP$ is 
  completely Boolean with $\poP=\BaB\setminus\ssetof{\bbzero_\BaB}$. Suppose
  $a\in\calH(\kappa)$ and $\varphi$ is a $\Sigma_1$-formula in 
  $\Lin$. If $\varphi(a)$ holds in $\uniV$, then clearly we also have
  $\forces{\poP}{\varphi(a)}$.

  Suppose now that $\varphi=\exists y\,\psi(x,y)$ for a bounded formula $\psi$ in $\Lin$, 
  and $\forces{\poP}{\varphi(a)}$. \Wolog, we may assume that $a\subseteq\mu$ for some 
  cardinal $\mu<\kappa$ (this is because $a$ can be reconstructed from $\trcl^+(a)$, and
  $\trcl^+(a)$ can be coded by a subset $a^*$ of $\cardof{\trcl^+(a)}$. The 
  formula $\varphi(a)$ can be replaced by the formula saying:
  \begin{xitemize}
  \item[] 
    $\exists x\,(\,x\mbox{ is the set ``}a
    \mbox{'' reconstructed from the transitive set coded 
      by }a^*\\
    \phantom{\exists x\,({}}\mbox{and }\varphi(x)\mbox{ holds}\,)$.
  \end{xitemize}
  Note that this formula is $\Sigma_1$ with 
  the parameter $a^*$ if $\varphi$ is $\Sigma_1$). We may also assume that $a$ is not an 
  ordinal (if necessary, we can replace $a$ with a subset of $\mu$ with some redundant 
  complexity to make $a\not\in\On$). 

  Let $\utb$ be a $\poP$-name \st\ $\forces{\poP}{\psi(a,\utb)}$. Let $\genG$ be a
  $(\uniV,\poP)$-generic filter and we work in $\uniV[\genG]$. Letting $b=\utb[\genG]$, we 
  have $\psi(a,b)$.

  Working further in $\uniV[\genG]$, let $\lambda$ be large enough \st\ $V_\lambda$ 
  satisfies large enough fragment of \ZFC, $a$, 
  $b\in V_\lambda$, and $V_\lambda\models\psi(a,b)$. Let $M\prec V_\lambda$ be \st\
  $\mu\subseteq M$, $a$, $b\in M$, and $\cardof{M}=\mu$. Let $\isomrph{m}{M}{M_0}$ be the 
  Mostowski collapse of $M$ and let $\nu=\On\cap M_0$. Note that we have
  $m\restr\mu\cup\ssetof{a}=\id_{\mu\cup\ssetof{a}}$. 

  Let 
  $\gmM:=\pairof{\nu+\mu,E,f}$ be the structure in the language
  $\calL:=\ssetof{\symb{E}, \symb{f}}$ \st\ there is an isomorphism 
  \begin{xitemize}
  \xitemA[x-theintro-0-0-a] $\isomrph{i}{\pairof{M_0, \in, rank}}{\pairof{\nu+\mu, E, f}}$\\
  \st\ $i\restr{\nu}=\id_{\nu}$, $i(a)=\nu$, and $i(m(b))=\nu+1$ 
  \end{xitemize}
  where $rank$ is the rank 
  function restricted to $M_0$.  Clearly we have 
  $\pairof{\nu+\mu,E}\models\psi^*(\nu,\nu+1))$, where $\psi^*$ is the formula 
  obtained from $\psi$ by replacing the symbol $\symb{\in}$ by $\symb{E}$.

  Let $\utgmM$, $\utE$, $\utf\in\uniV$ be $\poP$-names of $\gmM$, $E$ and $f$ respectively. By replacing 
  $\poP$ with $\poP\restr\condp$ for some $\condp\in\poP$ if necessary, we may assume that 
  \begin{xitemize}
  \xitemA[x-theintro-0-0-0] 
    all the properties of $\pairof{\nu+\mu, E, f}$ used below are forced  (as a statement 
    on $\pairof{\nu+\mu,\utE,\utf}$) by $\bbone_\poP$.
  \end{xitemize}

  In $\uniV$, let $\calD$ be the family of 
  dense subsets generated  in $\poP$ by each of the following sets of cardinality 
  $\leq\mu<\kappa$:
  \begin{xitemize}
  \xitemA[x-theintro-0-1]
    $\setof{\bvalof{\utilde{f}(\alpha)=\beta}{\BaB}}{\beta<\nu}\setminus\ssetof{\bbzero_\BaB}$,\qquad
    for all $\alpha\in\nu+\mu$.
  \xitemA[x-theintro-0-2] 
    $\ssetof{\bvalof{\utgmM\models\theta(a_0\ctentenc a_{k-1})}{\BaB},
    \bvalof{\utgmM\models\neg\theta(a_0\ctentenc a_{k-1})}{\BaB}}\setminus\ssetof{\bbzero_\BaB}$,\\
    \hfill for all $\Sigma_0$-formulas $\theta$ in $\calL$ 
    and $a_0\ctentenc a_{k-1}\in\nu+\mu$. 
  \xitemA[x-theintro-0-3]
    $\ssetof{\bvalof{\utgmM\models\eta\land\theta(a_0\ctentenc a_{k-1})}{\BaB},\,\\
    \ \bvalof{\utgmM\models\neg\eta(a_0\ctentenc a_{k-1})}{\BaB},\,
    \bvalof{\utgmM\models\neg\theta(a_0\ctentenc 
      a_{k-1})}{\BaB}}\setminus\ssetof{\bbzero_\BaB}$,\\[\jot]
    \hfill for all $\Sigma_0$-formulas $\eta$, 
    $\theta$ in $\calL$ and $a_0\ctentenc a_{k-1}\in\nu+\mu$.
  \xitemA[x-theintro-0-4]
    $\big(\ssetof{\bvalof{\neg(\exists x\,\symb{E}\, c)\,\eta(x,a_0\ctentenc a_{k-1})}{\BaB}}
    \ \cup\\
    \ \ \setof{\bvalof{d\,\symb{E}\,c
      \land\eta(d,a_0\ctentenc a_{k-1})}{\BaB}}{d\in \nu+\mu}\big)\setminus\ssetof{\bbzero_\BaB}$,
    \\[\jot]
    \hfill for all $\Sigma_1$-formulas $\eta=\eta(x,x_0\ctentenc x_{k-1})$  
    in $\calL$ and $c$, $a_0\ctentenc a_{k-1}\in\nu+\mu$.
  \end{xitemize}

  To see that each of the sets in \xitemof{x-theintro-0-1} is a maximal antichain in $\poP$ 
  of size $\leq\mu<\kappa$,
  suppose that $\alpha\in\nu+\mu$ and $\condp\in\poP$. Then there 
  is $\condq\leq_\poP\condp$ which decides $\utf(\alpha)$.  
  Since $\forces{\poP}{\utf(\alpha)\in\nu}$ by \xitemof{x-theintro-0-0-0}, if follows that 
  $\condq\forces{\poP}{\utf(\alpha)=\beta}$ for some $\beta\in\nu$. It is clear that elements 
  of each of the sets in \xitemof{x-theintro-0-1} are pairwise incompatible and these sets are of 
  size $\leq\mu<\kappa$.  

  It is also proved similarly that sets in \xitemof{x-theintro-0-4} are maximal antichains 
  in $\poP$ of size $\leq\mu<\kappa$. 

  Now, in $\uniV$, let $\genG$ be $\calD$-generic filter. $\genG$ exists by
  $\BFA_{\LT\kappa}(\calP)$, and since $\calD$ is  
  a family of maximal antichains of size $\LT\kappa$ with $\cardof{\calD}<\kappa$. 

  Let
  \begin{xitemize}
  \xitemx[] $\utgmM[\genG]:=\pairof{\nu+\mu, \utE[\genG], \utf[\genG]}$. 
  \end{xitemize}
  where
  \begin{xitemize}
  \xitemA[x-theintro-0-3-0]
    $\utE[\genG]:=
    \setof{\pairof{\xi,\eta}}{\xi,\eta\in\nu+\mu,\,\condp\forces{\poP}{
        \pairof{\xi,\eta}\in\utE}\mbox{ for some }\condp\in\genG}$,\quad and 
  \xitemA[x-theintro-0-3-1]
    $\utf[\genG]:=
    \setof{\pairof{\xi,\eta}}{\xi,\eta\in\nu+\mu,\,\condp\forces{\poP}{
        \pairof{\xi,\eta}\in\utf}\mbox{ for some }\condp\in\genG}$. 
  \end{xitemize}

  \begin{Claim}\extendedcolor\wassertof{1} $\utgmM[\genG]$ is an $\calL$-structure.\smallskip

    \wassert{2} For each $\Sigma_1$-formula $\theta=\theta(x_0\ctentenc x_{k-1})$ in $\calL$ and
    $a_0\ctentenc a_{k-1}\in\nu+\mu$,
    \begin{xitemize}
    \xitemA[x-theintro-0-4-a] 
      $\bvalof{\utgmM\models\theta(a_0\ctentenc a_{k-1})}{\BaB}\in\genG$\ \ if and only 
      if\/\ \ 
      $\utgmM[\genG]\models\theta(a_0\ctentenc a_{k-1})$. 
    \end{xitemize}

    \wassert{3} $\utE[\genG]$ is extensional and well-founded. $\utE[\genG]$ on $\nu+\mu$ coincides 
    with the canonical ordering on $\nu+\mu$. 
  \end{Claim}
  \prfofClaim \assertof{1}: Since the maximal antichains in \xitemof{x-theintro-0-1} 
  are in $\calD$, we have $\mapping{\utf[\genG]}{\nu+\mu}{\nu}$. \smallskip

  \assertof{2}: By induction on the construction of the formula $\theta$ using 
  \xitemof{x-theintro-0-2}, \xitemof{x-theintro-0-3}, and 
  \xitemof{x-theintro-0-4}.

  Suppose that $a_0$, $a_1\in\nu+\mu$. Then 
  \begin{xitemize}
  \xitemx[] $\bvalof{\utgmM\models a_0\,\symb{E}\,a_1}{\BaB}\in\genG$\ \ $\Leftrightarrow$\ \
    $\bvalof{a_0\,\utE\,a_1}{\BaB}\in\genG$\ \ $\Leftrightarrow$\ \
    $\exists\condp\in\genG\ \condp\forces{\poP}{a_0\,\utE\,a_1}$\ \
    {$\ubecause{\Leftrightarrow}{}{by \xitemof{x-theintro-0-3-0}}$}\ \ $a_0\,\utE[\genG]\,a_1$
    \ \ $\Leftrightarrow$\ \ $\utgmM[\genG]\models a_0\,\symb{E}\,a_1$.
  \end{xitemize}
  The proof for ``$a_0=\utf(a_1)$'' can be done similarly.

  The induction steps for $\neg$ and $\lor$ go through since $\genG$ is $\calD$-generic and 
  the antichains in  
  \xitemof{x-theintro-0-2} are in $\calD$. 

  Suppose now that the equivalence \xitemof{x-theintro-0-4-a} holds for a $\Sigma_1$-formula
  $\theta=\theta(x,x_0\ctentenc x_{k-1})$  and all other $\Sigma_1$-formulas with the 
  quantifier rank (\wrt\ bounded existential quantification) less than or equal to that of
  $\theta$.

  If $\bvalof{\utgmM\models(\exists x\,\symb{E}\,b)\,\theta(x,a_0\ctentenc a_{k-1})}{\BaB}\in\genG$,
  then, by  \xitemof{x-theintro-0-4}, there is
  $\bvalof{\utgmM\models d\in\nu_\mu\land\theta(d,a_0\ctentenc a_{k-1})}{\BaB}\in\genG$. 
  By the induction hypothesis, it follows that
  $\utgmM[\genG]\models d\,\symb{E}\,b\land\theta(d,a_0\ctentenc a_{k-1})$.
  Thus $\utgmM[\genG]\models (\exists x\,\symb{E}\,b)\,\theta(d,a_0\ctentenc a_{k-1})$.

  If $\utgmM[\genG]\models (\exists x\,\symb{E}\,b)\,\theta(d,a_0\ctentenc a_{k-1})$, then
  $\utgmM[\genG]\models d\,\symb{E}\,b \land \theta(d,a_0\ctentenc a_{k-1})$ for some
  $d\in\nu+\mu$. By induction hypothesis, it follows that
  $\bvalof{d\,\symb{E}\,b \land \theta(d,a_0\ctentenc a_{k-1})}{\BaB}\in\genG$. Since
  $\bvalof{d\,\symb{E}\,b \land \theta(d,a_0\ctentenc a_{k-1})}{\BaB}
  \leq_\BaB\bvalof{(\exists x\,\symb{E}\,b)\,\theta(x,a_0\ctentenc 
    a_{kp-1})}{\BaB}$, we have $\bvalof{(\exists x\,\symb{E}\,b)\,\theta(x,a_0\ctentenc 
    a_{kp-1})}{\BaB}\in\genG$ since $\genG$ is a filter.
\smallskip

  \assertof{3}: By \xitemof{x-theintro-0-0-0}, we have 
  $\forces{\poP}{\utgmM\models\mbox{Axiom of Extensionality}}$, and
  \begin{xitemize}
  \xitemA[x-theintro-0-5] 
    $\forces{\poP}{\utgmM\models\forall x\forall y\,(x\utE y\rightarrow \utf(x)<\utf(y))}$. 
  \end{xitemize}
  By \assertof{2}, it follows that $\utE[\genG]$ is extensional and the statement on the 
  structure $\utgmM[\genG]$ corresponding to \xitemof{x-theintro-0-5} holds. 
  A similar argument shows that the canonical ordering on $\nu$ coincides with
  $\utE[\genG]\restr \nu^2$. This and the property of $\utgmM[\genG]$ corresponding 
  to \xitemof{x-theintro-0-5} implies that $\utE[\genG]$ is well-founded. 
  \qedofClaim\qedskip
  
  By \Claimabove,\,\assertof{3}, we can take the Mostowski collapse of the structure
  $\utgmM[\genG]$ 
  $\isomrph{m^*}{\pairof{\nu+\mu,\utE[\genG]}}{\pairof{M_2,\in}}$. Since
  $\pairof{\nu+\mu,\utE[\genG]}\modelof{\psi^*(\nu,\nu+1)}$ by \xitemof{x-theintro-0-0-a}, 
  \xitemof{x-theintro-0-0-0} and \Claimabove,\,\assertof{2}, we have $m^*(\nu)=a$, and  
  $M_2\models\psi(a, m^*(\mu+1))$. Thus $M_2\models\varphi(a)$. 
  Since $\varphi$ is $\Sigma_1$, it follows 
  that $V\models\varphi(a)$. \smallskip

  \assertof{b} $\Rightarrow$ \assertof{a}: Suppose that $\poP\in\calP$ is complete
  Boolean and $\calD$ is a set of antichains each of size $<\kappa$ with $\cardof{\calD}<\kappa$.

  Let $X=\bigcup\calD$ then $\cardof{X}<\kappa$. Say, $\mu:=\cardof{X}$. Let $\lambda$ be 
  sufficiently large with $V_\lambda\prec_{\Sigma_n}\uniV$ for sufficiently large $n$. Let
  $M\prec V_\lambda$ be \st\ $\cardof{M}=\mu$, $\poP$, $\calD$, $X\in M$, and
  $\mu+1\subseteq M$. Note that $\calD\subseteq M$ and $I\subseteq M$ for each $I\in\calD$. 

  Let $\isomrph{m}{M}{M_0}$ be the Mostowski collapse 
  and $\pairof{\poP_0,\leq_{\poP_0}}:=m(\pairof{\poP,\leq_\poP})$.

  Since $(\uniV,\poP)$-generic filter $\genG$ generates an $(M_0,\poP_0)$-generic filter, 
  we have 
  \begin{xitemize}
  \xitemx[] 
    $\forces{\poP}{\calH(\kappa^{(+)})\models\mbox{\,there is a }(M_0,\poP_0)\mbox{-generic filter}}$. 
  \end{xitemize}

  By assumption it follows that
  $\calH(\kappa)\modelof{\mbox{\,there is a }(M_0,\poP_0)\mbox{-generic filter}}$ in $\uniV[\genG]$. Let 
  $\genG_0$ be such a filter. Then $m^{-1}\imageof{\genG_0}$ generates 
  a $\calD$-generic filter on $\poP$. }
  \else
  For ``\assertof{a} $\Leftrightarrow$ \assertof{b}'', see the proof 
  of ``\assertof{a} $\Leftrightarrow$ \assertof{b}'' of \Thmof{p-theintro-5} below. 
\fi}
\qedof{\Thmof{p-theintro-0}}\qedskip

The following is one of many nice applications of \Thmabove:  
\begin{Cor}
  If $\calP$ contains a \po\ adding a new real then $\BFA_{\LT\kappa}(\calP)$ for 
$\kappa>\aleph_1$ implies $\neg\CH$.
\end{Cor}
\prf Assume that $\BFA_{\LT\kappa}(\calP)$ holds for $\kappa>\aleph_1$, but also \CH\ holds in $\uniV$. Let
$a=\psof{\omega}$. We have $a\in\calH(\kappa)^\uniV$ by \CH\ and
$\calH(\kappa)\modelof{a\mbox{ is }\psof{\omega}}$. The statement can be formulated as a
$\Pi_1$-formula with the parameter $a$. But if $\poP\in\calP$ adds a real,
$\calH(\kappa)^{\uniV[\genG]}\modelof{a\mbox{ is not }\psof{\omega}}$. This is a 
contradiction to \Thmof{p-theintro-0},\,\assertof{c}.\qedofCor\qedskip

Suppose that $\calR$ is a definable class (proper or set). We shall say that a 
class $\calP$ of \pos\ is {\It provably correct for $\calR$} if the following is provable 
in \ZFC:
\begin{xitemize}
\xitem[p-theintro-2] for any $\poP\in\calP$ and $a$ $(\in\uniV)$, $a\in\calR$ $\Leftrightarrow$
  $\forces{\poP}{a\in\calR}$. 
\end{xitemize}
Thus if $\calP$ is provably correct for $\calR$ and $\genG$ is a $(\uniV,\poP)$-generic for 
a $\poP\in\calP$ then $(\uniV,\in,\calR^\uniV)$ is a (class) substructure of
$(\uniV[\genG],\in,\calR^{\uniV[\genG]})$. 

\memo{\mbox{$I_{NS}$ version of Bagaria's theorem}\\\mbox{Scan\_2024-05-14--12.42 gappo - annotated} pp.8-13]}
Let $\INS$ denote the non stationary ideal over $\omega_1$. Thus
\begin{xitemize}
\item[] 
  $\INS:=\setof{X\subseteq\omega_1}{X\mbox{ is non stationary}}$. 
\end{xitemize}\ifextended\qedskip\fi

\begin{Lemma}\Label{p-theintro-3} If (we can prove that) all $\poP\in\calP$ are stationary 
  preserving then $\calP$ is provably correct for $\INS$.\ifextended\else\qed\fi
\end{Lemma}
{\ifextended\extendedcolor\prf
  If $s\in\INS$ then there is a club $c\subseteq\omega_1$ \st\ $s\cap c=\emptyset$. 
  Since $\forces{\poP}{c\mbox{ is a club in }\omega_1\mbox{ and }s\cap c=\emptyset}$, we have
  $\forces{\poP}{s\in\INS}$.

  If $s\not\in\INS$ then $\forces{\poP}{s\not\in\INS}$ since $\poP$ is stationary 
  preserving. \qedofLemma\qedskip\fi}

Let $\calR$ be (the $\Lin$-definition of) a class. 
Let ${\calL}_{\in,\symb{\calR}}$ be the language which extends $\Lin$ with a new unary 
predicate symbol $\symb{\calR}$ where $M\models\symb{\calR}(a)$ is interpreted as $a\in\calR^M$ in 
an $\in$-structure $M$. In the following, we shall often identify the (definition) of the 
class  
$\calR$ with the symbol $\symb{\calR}$ of $\calR$,  and simply write $\calR$ and
$\calL_{\in,\calR}$ instead of $\symb{\calR}$ and $\calL_{\in,\symb{\calR}}$. This also 
applies when we are talking about $\INS$ and $\calL_{\in,\INS}$. 

The following \Lemmaof{p-theintro-4} can be proved in the same way as with the corresponding lemma 
for $\Lin$\ formulas: 
\ifextended\qedskip\fi

\begin{Lemma}\Label{p-theintro-4} For transitive (sets or classes) $M$, $N$ with
  $M\subseteq N$ and a class $\calR$ (i.e.\ an $\Lin$-formula with one single free variable)
  \st\ $\calR^M=\calR^N\cap M$, we have:\smallskip

  \wassert{1}
  $\pairof{M,\in,\calR^M}\models\varphi(\ol{a})$ $\Leftrightarrow$ 
  $\pairof{N,\in,\calR^N}\models\varphi(\ol{a})$ for 
  all $\Sigma_0$-formula $\varphi=\varphi(\ol{x})$ in 
  $\calL_{\in,\calR}$ and $\ol{a}\in M$. \smallskip

  \wassert{2}   $\pairof{M,\in,\calR^M}\models\varphi(\ol{a})$ $\Rightarrow$ 
  $\pairof{N,\in,\calR^N}\models\varphi(\ol{a})$ for 
  all $\Sigma_1$-formula $\varphi=\varphi(\ol{x})$ in 
  $\calL_{\in,\calR}$ and $\ol{a}\in M$. \qed
\end{Lemma}

\begin{Lemma}\Label{p-theintro-4-0}
  For a $\Sigma_1$-formula $\varphi$ in $\calL_{\in,\INS}$ we can find a 
  $\Sigma_2$-formula in $\Lin$ with the parameter $\omega_1$ equivalent to $\varphi$.
\end{Lemma}
\prf ``$x\in\INS$'' can be expressed by a $\Sigma_1$-formula in $\Lin$:
\begin{xitemize}
\item[] $\exists y\,(y\subseteq\omega_1\ \land\ x\subseteq\omega_1\ \land\ 
y\mbox{ is a club in }\omega_1\ \land\ x\cap y=\emptyset)$
\end{xitemize}
Thus ``$x\not\in\INS$'' can be expressed by a $\Pi_1$-formula in $\Lin$. \qedofLemma
\qedskip

\begin{Lemma}{\rm(A special case of Lemma 6.3 in Venturi and Viale 
    \cite{venturi-viale})}\Label{p-theintro-4-1} 
  For a cardinal $\lambda\geq2^{\aleph_1}$, we have $\pairof{\calH(\lambda),\,{\in},\,\INS}
\prec_{\Sigma_1}\rpairof{\uniV,\,{\in},\,\INS}$. \ifextended\else\qed\fi
\end{Lemma}
{\ifextended\extendedcolor
\prf If $\pairof{\calH(\lambda),\,{\in},\,I_\NS}\models\varphi(a)$ for a $\Sigma_1$-formula 
in $\calL_{\in,\INS}$, and $a\in\calH(\lambda)$, 
then $\rpairof{\uniV,\,{\in},\,\INS}\models\varphi(a)$ by 
\Lemmaof{p-theintro-4},\,\assertof{2}. 

Suppose now that $\rpairof{\uniV,\,{\in},\,\INS}\models\varphi(a)$ for $\varphi$ and $a$ as 
above. Suppose that $\varphi=\exists y\,\psi(x,y)$ for a $\Sigma_0$-formula $\psi$ in
$\calL_{\in,\INS}$ and let $b$ be \st\ $\rpairof{\uniV,\,{\in},\,\INS}\models\psi(a,b)$. 

Let $\alpha\in\On$ be sufficiently large \st\
$\pairof{V_\alpha,\,{\in},\,\INS}\prec_{\Sigma_n}\rpairof{\uniV,\,{\in},\,\INS}$ for 
sufficiently large $n\in\omega$. 

Let 
$\mu:=\sup\ssetof{\cardof{\trcl^+(a)},\,2^{\aleph_1}}$ and\footnote{We denote with
  $\trcl^+(a)$ the variant of transitive closure which is the minimal transitive set $T$ 
  with $a\cup\ssetof{a}\subseteq T$.} 
$\pairof{M,\,{\in}\,\INS}\prec\rpairof{\uniV,\,{\in},\,\INS}$ be \st\ $\cardof{M}=\mu$, 
$\trcl^+(a)$, $\mu+1$, $\INS\subseteq M$, and $b\in M$. Then we have
$\pairof{M,\,{\in},\,\INS}\models\psi(a,b)$.

Let $\isomrph{m}{M}{M_0}$ be the Mostowski collapse. Then we have 
$m\restr\trcl^+(a)=\id_{\trcl^+(a)}$ and $m\restr\psof{\omega_1}=id_{\psof{\omega_1}}$. It 
follows that $\pairof{M_0,\,{\in},\,\INS}\models\psi(a,m(b))$ and hence
$\pairof{M_0,\,{\in},\,\INS}\models\varphi(a)$. By \Lemmaof{p-theintro-4},\,\assertof{2}, 
it follows that $\pairof{\calH(\lambda),\,{\in},\,I_\NS}\models\varphi(a)$. 
\qedofLemma\qedskip\fi}

The following \Thmof{p-theintro-5} is an extension of Bagaria's Absoluteness 
\Thmof{p-theintro-0}. 
A special case of this theorem (the case where $\calP={}$the stationary preserving \pos) is 
also attributed to Bagaria in \cite{woodin-book}. Though \Thmof{p-theintro-5} in its 
generality must have been known, we included it here since we could not find any proof in 
the literature. 

We consider the following ``plus''-version of Bounded Forcing Axioms: For a (normal) class of 
\pos\ $\calP$,  
\begin{xitemize}
\item[{\darkred($\BFA^{+\LT\kappa}_{\LT\kappa}(\calP)$): }] For any complete Boolean 
  $\poP\in\calP$, a family $\calD$ of maximal antichains in $\poP$ \st\ 
$\cardof{\calD}<\kappa$ and $\cardof{I}<\kappa$ for all $I\in\calD$, and for a set $\calS$ 
  of $\calP$-names of cardinality $<\kappa$ \st\ each $\utS\in\calS$ is a $\poP$-name of a 
  stationary subset of $\omega_1$, 
  there is 
  a $\calD$-generic filter $\genG$ on $\poP$ \st\ $\utS[\genG]$ for all $\utS\in\calS$ are 
  stationary subsets of $\omega_1$. 
\end{xitemize}

%% In the formulation of \assertof{c} of the following theorem, we are using the notation we 
%% introduced right before \Thmof{p-theintro-0}. 

\begin{Thm}\Label{p-theintro-5} Suppose that $\calP$ is a class of \pos\ closed under 
  forcing 
  equivalence, and restriction (in the sense of \xitemof{x-theintro-0-0-a-0}) 
  \st\ all elements 
  of $\calP$ are stationary preserving and $\kappa=2^{\aleph_0}=2^{\aleph_1}$.\footnotemark\ Then \tfae:\smallskip\\
  \wassert{a}
  $\BFA^{+\LT\kappa}_{\LT\kappa}(\calP)$.\smallskip

  \wassert{b}
  For any $\Sigma_1$-formula $\varphi=\varphi(x)$ in $\calL_{\in,\INS}$, 
  $a\in\calH(\kappa)$, and $\poP\in\calP$, we 
  have 
  \begin{xitemize}
  \xitemx[] 
    $\forces{\poP}{\varphi(a)}$\  $\Leftrightarrow$\ $\varphi(a)$.
  \end{xitemize}

  \wassert{c} For any $\poP\in\calP$, and $(\uniV,\poP)$-generic $\genG$, we have
  \begin{xitemize}
  \item[] 
    $\pairof{\calH(\continuum)^\uniV,\in,\INS^\uniV}\prec_{\Sigma_1}
    \pairof{\calH(\left(2^{\aleph_1}\right)^{\uniV[\genG]})^{\uniV[\genG]},
      \in,\INS^{\uniV[\genG]}}$. 
  \end{xitemize}
\end{Thm}
\footnotetext{Note that this implies $\neg\CH$.}
\prf The equivalence of \assertof{b} and \assertof{c} follows from 
\Lemmaof{p-theintro-4-1}. 
%% \newpage
%% The equivalence mentioned in the last statement of the theorem 
%% holds since $\calL_{\in,\INS}$ version of \xitemof{x-theintro-0} 
%% holds.
%% 
%% Suppose that $\varphi$ is a $\Sigma_1$ formula in $\calL_{\in,\INS}$, $a\in\calH(\kappa)$ 
%% and $\poP\in\calP$.
%% If $\varphi(a)$, then $\forces{\poP}{\varphi(a)}$ by \Lemmaof{p-theintro-3} and 
%% \Lemmaof{p-theintro-4},\,\assertof{2}. 
%% 
%% If $\forces{\poP}{\varphi(a)}$, then by \Lemmaof{p-theintro-4-0} and
%% $(\calP,\calH(\kappa))_{\Sigma_2}$-\RcAp, there is a $\calP$-ground $\uniW$ of $\uniV$ \st\
%% $a\in\uniW$ and $\uniW\models\varphi(a)$.
%% By the assumption on $\calP$, \Lemmaof{p-theintro-3}, and 
%% \Lemmaof{p-theintro-4},\,\assertof{2}, it follows that $\uniV\models\varphi(a)$ as desired. 
%%\begin{spacing}{2}

\assertof{a} $\Rightarrow$ \assertof{b}: Assume that
$\BFA^{+\LT\kappa}_{\LT\kappa}(\calP)$ holds, and let $\poP\in\calP$. \Wolog, we may 
assume that $\poP$ is 
completely Boolean with $\poP=\BaB\setminus\ssetof{\bbzero_\BaB}$. Suppose
$a\in\calH(\kappa)$ and $\varphi$ is a $\Sigma_1$-formula in 
$\calL_{\in,\INS}$. If $\varphi(a)$ holds in $\uniV$, then we also have
$\forces{\poP}{\varphi(a)}$ by \Lemmaof{p-theintro-4}.

Suppose now that $\varphi=\exists y\,\psi(x,y)$ for a bounded formula $\psi$ in $\calL_{\in,\INS}$, 
and $\forces{\poP}{\varphi(a)}$. \Wolog, we may assume that $a\subseteq\mu$ for some 
cardinal $\mu<\kappa$ (this is because $a$ can be reconstructed from $\trcl^+(a)$, and
$\trcl^+(a)$ can be coded by a subset $a^*$ of $\cardof{\trcl^+(a)}$). The 
formula $\varphi(a)$ can be then replaced by another formula saying:
\begin{xitemize}
\item[] 
  $\exists x\,(\,x\mbox{ is the set ``}a
  \mbox{'' reconstructed from the transitive set coded 
    by }a^*\\
  \phantom{\exists x\,({}}\mbox{and }\varphi(x)\mbox{ holds}\,)$.
\end{xitemize}
Note that this formula is $\Sigma_1$ in $\calL_{\in,\INS}$ with 
the parameter $a^*$ if $\varphi$ is $\Sigma_1$ in $\calL_{\in,\INS}$.

We may also assume 
that $a$ is not an  
ordinal (if necessary, we can replace $a$ with a subset of $\mu$ with some redundant 
complexity to make $a\not\in\On$). 

Let $\utb$ be a $\poP$-name \st\ $\forces{\poP}{\psi(a,\utb)}$. Let $\genG$ be a
$(\uniV,\poP)$-generic filter and we work in $\uniV[\genG]$. Letting $b=\utb[\genG]$, we 
have $\psi(a,b)$.

Working further in $\uniV[\genG]$, let $\lambda$ be large enough \st\ $V_\lambda$ 
satisfies a large enough fragment of \ZFC, $a$, 
$b\in V_\lambda$, and $V_\lambda\models\psi(a,b)$. Let $M\prec V_\lambda$ be \st\
$\mu\subseteq M$, $a$, $b\in M$, and $\cardof{M}=\mu<\kappa$. Note that we have 
$\pairof{M,\in, \INS\cap M}\prec\pairof{V_\lambda,\in,\INS}$ since $\INS$ is definable in
$\pairof{V_\lambda,\in}$. 
Let $\isomrph{m}{M}{M_0}$ be the 
Mostowski collapse of $M$ and let $\nu=\On\cap M_0$. Note that we have
$m\restr\mu\cup\ssetof{a}\cup(\INS\cap M)=\id_{\mu\cup\ssetof{a}\cup(\INS\cap M)}$ and 
hence $\INS\cap M=\INS\cap M_0$. 

Let 
$\gmM:=\pairof{\nu+\mu,E,I,f,g}$ be the structure in the language
$\calL:=\ssetof{\symb{E}, \symb{I}, \symb{f}, \symb{g}}$ \st\ there is an isomorphism 
\begin{xitemize}
\xitem[a:x-theintro-0-0-a] $\isomrph{i}{\pairof{M_0, \in, \INS\cap M_0, rank, g_0}}{\pairof{\nu+\mu, E,I,f,g}}$\\
  \st\ $i\restr{\nu}=\id_{\nu}$, $i(a)=\nu$, and $i(m(b))=\nu+1$ 
\end{xitemize}
where $rank$ is the rank 
function restricted to $M_0$ and $\mapping{g_0}{M_0}{M_0}$ is a mapping \st\ $g_0\restr\mu$ 
is an enumeration of $(\psof{\omega_1}\cap M_0)\setminus\INS$
($=(\xmbox{the set of all stationary subsets of }\omega_1)^{M_0}$) and
$g_0\imageof{(M_0\setminus\mu)}=\ssetof{\emptyset}$. 

Clearly,
$\pairof{\nu+\mu,E, I}\models\psi^*(\nu,\nu+1))$ where $\psi^*$ is the formula obtained from
$\psi$ by replacing symbols $\symb{\in}$ and $\symb{\INS}$ in $\psi$ by 
$\symb{E}$ and $\symb{I}$.

Note that we have
\begin{xitemize}
\xitem[a:x-theintro-1-0] 
  $\pairof{\nu+\mu, E,I,f,g}\models
  \forall x\subseteq\omega_1\ (I(x)\lor \exists\alpha<\kappa\,(g(\alpha)=x))$. 
\end{xitemize}

Let $\utgmM$, $\utE$, $\utI$,  $\utf$, $\utg\in\uniV$ be $\poP$-names of $\gmM$, $E$, $I$,
$f$ and $g$ respectively. By replacing 
$\poP$ with $\poP\restr\condp$ for some $\condp\in\poP$ if necessary, we may assume that 
\begin{xitemize}
\xitem[a:x-theintro-0-0-0] 
  all the properties of $\pairof{\nu+\mu, E, I,f,g}$ used below are forced  (as a statement 
  on $\pairof{\nu+\mu,\utE, \utI, \utf,\utg}$) by $\bbone_\poP$.
\end{xitemize}

In $\uniV$, let $\calD$ be the family of 
dense sets in $\poP$ generated by the following predense sets in $\poP$, each of 
size $\leq\mu<\kappa$.
Note that since $\poP$ is assumed to be completely Boolean, each predense subset of $\calP$ 
can be replaced by a maximal antichain of at most the same size.
\begin{xitemize}
\xitem[a:x-theintro-0-1]
  $\setof{\bvalof{\utilde{f}(\alpha)=\beta}{\BaB}}{\beta<\nu}\setminus\ssetof{\bbzero_\BaB}$,
  and\\
  $\setof{\bvalof{\utilde{g}(\alpha)=\beta}{\BaB}}{\beta<\nu+\mu}\setminus\ssetof{\bbzero_\BaB}$, 
  \qquad
  for all $\alpha\in\nu+\mu$.
\xitem[a:x-theintro-0-2] 
  $\ssetof{\bvalof{\utgmM\models\theta(a_0\ctentenc a_{k-1})}{\BaB},
  \bvalof{\utgmM\models\neg\theta(a_0\ctentenc a_{k-1})}{\BaB}}\setminus\ssetof{\bbzero_\BaB}$,\\
  \hfill for all $\Sigma_0$-formulas $\theta$ in $\calL$ 
  and $a_0\ctentenc a_{k-1}\in\nu+\mu$. 
\xitem[a:x-theintro-0-3]
  $\ssetof{\bvalof{\utgmM\models\eta\land\theta(a_0\ctentenc a_{k-1})}{\BaB},\,\\
  \ \bvalof{\utgmM\models\neg\eta(a_0\ctentenc a_{k-1})}{\BaB},\,
  \bvalof{\utgmM\models\neg\theta(a_0\ctentenc 
    a_{k-1})}{\BaB}}\setminus\ssetof{\bbzero_\BaB}$,\\[\jot]
  \hfill for all $\Sigma_0$-formulas $\eta$, 
  $\theta$ in $\calL$ and $a_0\ctentenc a_{k-1}\in\nu+\mu$.
\xitem[a:x-theintro-0-4]
  $\big(\ssetof{\bvalof{\neg(\exists x\,\symb{E}\, c)\,\eta(x,a_0\ctentenc a_{k-1})}{\BaB}}
  \ \cup\\
  \ \ \setof{\bvalof{d\,\symb{E}\,c
      \land\eta(d,a_0\ctentenc a_{k-1})}{\BaB}}{d\in \nu+\mu}\big)\setminus\ssetof{\bbzero_\BaB}$,
  \\[\jot]
  \hfill for all $\Sigma_1$-formulas $\eta=\eta(x,x_0\ctentenc x_{k-1})$  
  in $\calL$ and $c$, $a_0\ctentenc a_{k-1}\in\nu+\mu$.
\end{xitemize}

To see that each of the sets in \xitemof{a:x-theintro-0-1} is a maximal antichain in $\poP$ 
of size $\leq\mu$,
suppose that $\alpha\in\nu+\mu$ and $\condp\in\poP$. Then there 
is $\condq\leq_\poP\condp$ which decides $\utf(\alpha)$.  
Since $\forces{\poP}{\utf(\alpha)\in\nu}$ by \xitemof{a:x-theintro-0-0-0}, if follows that 
$\condq\forces{\poP}{\utf(\alpha)=\beta}$ for some $\beta\in\nu$. 
For the sets corresponding to $\utg$ the argument is the same. 

It is clear that elements 
of each of the sets in \xitemof{a:x-theintro-0-1} are pairwise incompatible, and these sets are of 
size $\leq\mu<\kappa$.  

The sets in \xitemof{a:x-theintro-0-3} and \xitemof{a:x-theintro-0-4} are not necessarily 
maximal antichains but it can be proved similarly that each of them sums up to $\bbone_\BaB$ and 
is of size $\leq\mu<\kappa$. 

By \xitemof{a:x-theintro-0-0-0}, we have that 
\begin{xitemize}
\xitem[a:x-theintro-0-4-a]
  $\forces{\poP}{\setof{\xi\in\omega_1}{\xi\,\utE\,\utg(\alpha)}
  \mbox{ is a stationary subset of }\omega_1}$. 
\end{xitemize}

Now, in $\uniV$, let $\genG$ be $\calD$-generic filter 
\st\ 
\begin{xitemize}
\xitem[a:x-theintro-0-4-a-0] $\utg(\alpha)[\genG]:=\setof{\xi\in\omega}{
  \condp\forces{\poP}{\xi\,\utE\,\utg(\alpha)}\mbox{ for some }\condp\in\genG}$
  is a stationary subset of $\omega_1$ for all $\alpha<\mu$. 
\end{xitemize}
$\genG$ exists by
$\BFA^{+\LT\kappa}_{\LT\kappa}(\calP)$, by \xitemof{a:x-theintro-0-4-a}, and since $\calD$ is  
a family of dense sets generated by small sets with $\cardof{\calD}<\kappa$. 

Let
\begin{xitemize}
\xitemx[] $\utgmM[\genG]:=\pairof{\nu+\mu, \utE[\genG], \utI[\genG], \utf[\genG], \utg[\genG]}$. 
\end{xitemize}
where
\begin{xitemize}
\xitemx[]
  $\utE[\genG]:=
  \setof{\pairof{\xi,\eta}}{\xi,\eta\in\nu+\mu,\,\condp\forces{\poP}{
      \pairof{\xi,\eta}\in\utE}\mbox{ for some }\condp\in\genG}$,
\xitemx[]
  $\utI[\genG]:=
  \setof{\xi}{\xi\in\nu+\mu,\,\condp\forces{\poP}{
      \xi\in\utI}\mbox{ for some }\condp\in\genG}$,
\xitemx[]
  $\utf[\genG]:=
  \setof{\pairof{\xi,\eta}}{\xi,\eta\in\nu+\mu,\,\condp\forces{\poP}{
      \pairof{\xi,\eta}\in\utf}\mbox{ for some }\condp\in\genG}$,\quad and 
\xitemx[]
  $\utg[\genG]:=
  \setof{\pairof{\xi,\eta}}{\xi,\eta\in\nu+\mu,\,\condp\forces{\poP}{
      \pairof{\xi,\eta}\in\utg}\mbox{ for some }\condp\in\genG}$. 
\end{xitemize}

\begin{Claim}\wassertof{1} $\utgmM[\genG]$ is an $\calL$-structure.\smallskip

  \wassert{2} For each $\Sigma_1$-formula $\theta=\theta(x_0\ctentenc x_{k-1})$ in $\calL$ and
  $a_0\ctentenc a_{k-1}\in\nu+\mu$,
  \begin{xitemize}
  \xitem[a:x-theintro-0-4-0] 
    $\bvalof{\utgmM\models\theta(a_0\ctentenc a_{k-1})}{\BaB}\in\genG$\ \ if and only 
    if\/\ \ 
    $\utgmM[\genG]\models\theta(a_0\ctentenc a_{k-1})$. 
  \end{xitemize}

  \wassert{3} $\utE[\genG]$ is extensional and well-founded. $\utE[\genG]$ on $\nu+\mu$ coincides 
  with the canonical ordering on $\nu+\mu$. 
\end{Claim}
\prfofClaim \assertof{1}: Since the maximal antichains in \xitemof{a:x-theintro-0-1} 
are in $\calD$, we have $\mapping{\utf[\genG]}{\nu+\mu}{\nu}$ and
$\mapping{\utg[\genG]}{\nu+\mu}{\nu+\mu}$. \smallskip

\assertof{2}: By induction on the construction of the formula $\theta$ using 
\xitemof{a:x-theintro-0-2}, \xitemof{a:x-theintro-0-3}, and 
\xitemof{a:x-theintro-0-4}.\smallskip

\assertof{3}: By \xitemof{a:x-theintro-0-0-0}, we have 
$\forces{\poP}{\utgmM\models\mbox{Axiom of Extensionality}}$, and
\begin{xitemize}
\xitemA[a:x-theintro-0-5] 
  $\forces{\poP}{\utgmM\models\forall x\forall y\,(x\utE y\rightarrow \utf(x)<\utf(y))}$. 
\end{xitemize}
By \assertof{2}, it follows that $\utE[\genG]$ is extensional and the statement on the 
structure $\utgmM[\genG]$ corresponding to \xitemof{a:x-theintro-0-5} holds. 
A similar argument shows that the canonical ordering on $\nu$ coincides with
$\utE[\genG]\restr \nu^2$. This and the property of $\utgmM[\genG]$ corresponding 
to \xitemof{a:x-theintro-0-5} implies that $\utE[\genG]$ is well-founded. 
\qedofClaim\qedskip

By \Claimabove,\,\assertof{3}, we can take the Mostowski collapse of the structure
$\utgmM[\genG]$ 
$\isomrph{m^*}{\pairof{\nu+\mu,\utE[\genG],\utI[\genG]}}{\pairof{M_2,\in, I^*}}$. Since
$\pairof{\nu+\mu,\utE[\genG],\utI[\genG]}\modelof{\psi^*(\nu,\nu+1)}$ by \xitemof{a:x-theintro-0-0-a}, 
\xitemof{a:x-theintro-0-0-0} and \Claimabove,\,\assertof{2}, we have $m^*(\nu)=a$ and hence 
$\pairof{M_2,\in, I^*}\models\psi(a,m^*(\nu+1))$ where the predicate $\INS$ is interpreted 
as $I^*$. Thus $\pairof{M_2,\in, I^*}\models\varphi(a)$. By \xitemof{a:x-theintro-1-0}, 
\xitemof{a:x-theintro-0-0-0}, \xitemof{a:x-theintro-0-4-a-0} and 
\Claimabove,\,\assertof{2}, we have $I^*=\INS\cap M_2$. 

Since $\varphi$ is $\Sigma_1$, it follows 
that $V\models\varphi(a)$ by \Lemmaof{p-theintro-4},\,\assertof{2}. \smallskip

\assertof{b} $\Rightarrow$ \assertof{a}: Assume that \assertof{b} holds, and suppose that
$\poP\in\calP$ is complete 
Boolean, $\calD$ is a set of maximal antichains each of size $<\kappa$ 
with $\cardof{\calD}<\kappa$,  
and $\calS$ is a set of $\poP$-names of stationary subsets of $\omega_1$ with
$\cardof{\calS}<\kappa$. 

By replacing elements of $\calS$ by equivalent $\poP$-names which are sufficiently nice, 
we may assume that each element of $\calS$ is nice $\poP$-name of size $\aleph_1$ (this 
is possible since we assumed that $\poP$ is completely Boolean). 

Let $X=\bigcup\calD$ then $\cardof{X}<\kappa$. Let
$\mu:=\max\ssetof{\cardof{X},\cardof{\calS}}$. Let $\lambda$ be  
sufficiently large with $V_\lambda\prec_{\Sigma_n}\uniV$ for sufficiently large $n$. Let
$M\prec V_\lambda$ be \st\ $\cardof{M}=\mu$, 
\ixitem[x-theintro-0-6] $\poP$, $\calD$, $X$, $\calS\in M$, and 
  $\mu+1\subseteq M$. Note that \xitemof{x-theintro-0-6} implies $\calD,\calS\subseteq M$ 
  and $I,\utS\subseteq M$ for each $I\in\calD$ and $\utS\in\calS$. 

  Let $\isomrph{m}{M}{M_0}$ be the Mostowski collapse 
  and $\pairof{\poP_0,\leq_{\poP_0}}:=m(\pairof{\poP,\leq_\poP})$.
  Let $\calS_0:=\setof{m(\utS)}{\utS\in\calS}$. 

  Since $(\uniV,\poP)$-generic filter $\genG$ generates an $(M_0,\poP_0)$-generic filter, 
  we have 
  \begin{xitemize}
  \xitemx[] 
    $\forces{\poP}{%%
    \begin{array}[t]{@{}l}
      \mbox{\,there is a }(M_0,\poP_0)\mbox{-generic filter
        which realizes each element of }\calS_0\\
      \mbox{ to be a stationary subset of }\omega_1}.
    \end{array}
    $ 
  \end{xitemize}

  By the assumption \assertof{b}, it follows that
  \begin{xitemize}
  \xitemx[] 
    $\uniV\modelof{%%
    \begin{array}[t]{@{}l}
      \mbox{\,there is a }(M_0,\poP_0)\mbox{-generic filter
        which realizes each element of }\calS_0\\
      \mbox{ to be a stationary subset of }\omega_1}.
    \end{array}
    $ 
  \end{xitemize}

  Let 
  $\genG_0$ be such a filter. Then $m^{-1}\imageof{\genG_0}$ generates 
  a $\calD$-generic filter $\genG_1$ on $\poP$  which realizes each element of $\calS$ to be 
  a stationary subset of $\omega_1$. 
%%\end{spacing}
\qedofThm\qedskip

\memo{def. of Laver-genericity. Relation of \RcAp\ to Laver-genericity\\
Joel Paper, Ikegami and Trang Theorem 1.6
}

{\ifprivate\privatecolor
The recent progress in the research, in particular the following \ThmAof{p-theintro-a-0} in 
a 2021 paper \cite{aspero-schindler} by D.\ Asperó and R.\ Schindler in particular, 
clarifies further the relationship of (the assumptions in) Viale's 
\Thmof{p-theintro-a} to the $\Omega$-logic context. 

\begin{ThmA}\privatecolor{\rm (Asperó, and Schindler \cite{aspero-schindler})}\Label{p-theintro-a-0} 
  $\MM^{++}$ implies Woodin's $(*)$.\qed
\end{ThmA}

Theorem 4.69 in \cite{woodin-book} (see also Theorem 1.1 in 
\cite{larson-lumsdaine-yin})together with the \Thmabove\ above  
implies:\footnote{Hiroshi Sakai made the first author aware of \ThmAof{p-theintro-a-1}.}
\begin{ThmA}\Label{p-theintro-a-1}{\privatecolor Suppose that there are proper class many 
    Woodin cardinals,  
  and $\MM^{++}$ holds. Then, for any $\Pi_2$-sentence $\varphi$ (i.e.\ without parameters)
  in $\calL_{\in,\INS}$, and for any 
  \po\ $\poP$ with $(\uniV, \poP)$-generic $\genG$, \memox{$\Pi_1$??}
  \begin{xitemize}
  \item[] if $(\calH({\aleph_2}^{\uniV[\genG]}),\in,\INS)^{\uniV[\genG]}\models\varphi$ then
    $(\calH({\aleph_2}^{\uniV}),\in,\INS)^{\uniV}\models\varphi$. \qed
  \end{xitemize}}

\end{ThmA}

In the context of Laver-genericity we are going to discuss in \sectionof{genabs-Laver}, the 
axiom ``there is a tight $\calP$-Laver-generic hyperhuge cardinal for $\calP={}$the 
stationary preserving \pos''\footnote{We call this axiom ``the tight $\calP$-Laver 
  generic hyperhuge cardinal axiom. In general this type of axioms are called Laver-generic 
  large cardinal axioms.} implies all the assumptions  
of \Thmof{p-theintro-a} and \ThmAof{p-theintro-a-1} (see \cite{sfetal-II} and 
\cite{recurrence}), and hence also the conclusions of the theorems. \memo{Part of these 
  paragraphs should be put in in \sectionof{misc}.}

Note that it is not at all trivial that the existence of a tight $\calP$-Laver-generic hyperhuge 
cardinal for $\calP=$ stationary preserving \pos\ implies the existence of class many large 
cardinals: a tight $\calP$-Laver-generic  
hyperhuge cardinal $\kappa$ implies the existence of the bedrock $\ol{\uniW}$ \st\ $\kappa$ 
is hyperhuge in the bedrock. This implies that there are class many huge cardinals in
$\ol{\uniW}$. All of these huge cardinals except an initial segment of them survive when we 
set-generically extend $\ol{\uniW}$ to reach $\uniV$ (see \cite{recurrence} or \cite{janos}). 

Though logically, the assertion that the tight $\calP$-Laver generic hyperhuge cardinal axiom 
implies the conclusions of \Thmof{p-theintro-a} and \ThmAof{p-theintro-a-1}, is merely a 
weakening of these theorems, the 
fact can be seen as a part of the ``Laver-Generic Maximum'' discussed in \cite{recurrence} and 
\cite{janos}: all known important axioms 
and principles of set theory are either conclusions of one of the strongest form of Laver-generic 
large cardinal axiom, or theorems in (many) grounds of the universe satisfying the axiom. 

%% Note that the tight $\calP$-Laver gen.\ hyperhuge cardinal axiom for $\calP=$ stationary 
%% preserving \pos\ (i.e.\ the statement that $\aleph_2$ has this generic large cardinal property) 
%% implies the 
%% assumptions of \Thmof{p-theintro-a} and \Thmof{p-theintro-a-1} (see \cite{recurrence}).

Concerning \Thmof{p-theintro-a}, the conclusion of the theorem already follows from 
the weaker assumption of the tight 
Laver-gen. superhuge cardinal axiom (see \Thmof{p-genabs-Laver-1} and 
\footnoteof{fn-genabs-Laver-0} below).  \fi}
{\ifprivate\privatecolor
\begin{Problem} \wassertof{1} Does the tightly $\calP$-Laver-generic hyperhuge cardinal 
  axiom for $\calP$ other than stationary preserving \pos\ ( or equivalently semi-proper 
  \pos\ --- since 
  the Laver-genericity for both of them imply \MM) similar conclusion as that of 
  \ThmAof{p-theintro-a-1}?\smallskip

 \wassert{2} For $\calP={}$the stationary preserving \pos, 
  can the tight $\calP$-Laver-generic large cardinal axiom for a notion of large 
  cardinal with much less consistency strength than that of hyperhuge imply the conclusion 
  of \ThmAof{p-theintro-a-1}? \smallskip

  \wassert{3} Is the statement of \ThmAof{p-theintro-a-1} still consistent if $\INS$ is 
  replaced by $\setof{S\subseteq[\omega_1]^{\aleph_0}}{S\mbox{ is stationary}}$ (perhaps 
  no: this would imply an absoluteness which seems to be too strong)?
\end{Problem}
\fi}
\section{Recurrence Axioms and the Ground Axiom}
\Label{rec-GA}
\subsection{Hierarchies of Recurrence and Maximality}\Label{hierarchies-rec-max}
\memox{Definition of $(\calP,\calH(\kappa))_\Gamma$-\RcAp\ 
  $\BFA_{\LT\kappa}(\calP)$}

The term ``Recurrence Axiom'' was coined in Fuchino and Usuba \cite{recurrence} (see also 
Fuchino \cite{janos}). The Recurrence Axiom for a (normal) 
class $\calP$ of \pos, a set $A$ of parameters, and a set $\Gamma$ of $\Lin$-formulas
(\,{$(\calP,A)_\Gamma$-\RcA}, for short) is the statement \xitemof{x-theintro-1} below. 

A {\It ground of a }(transitive set or class) {\It model 
$\uniW$} (of some set theory) is an inner model $\uniW_0$ of $\uniW$ \st\ there is a \po\
$\poP\in\uniW_0$ \st\  
$\uniW$ is a $\poP$-generic extension of $\uniW_0$. For a class $\calP$ of \pos, a ground
$\uniW_0$ of $\uniW$ is a {\It$\calP$-ground of\/ $\uniW$}\/ if there is a \po\ $\poP\in\uniW_0$ \st\ 
$\uniW_0\modelof{\poP\in\calP\,}$ and $\uniW$ is a $\poP$-generic extension of $\uniW_0$.

The Recurrence Axiom {\It$(\calP,A)_\Gamma$-\RcA} is the following statement formulated in 
an axiom scheme in $\Lin$ (that this axiom is not formalizable in a single formula is 
discussed in \cite{future}): 

\begin{xitemize}
\xitem[x-theintro-1] 
  For any $\poP\in\calP$, $\varphi(\overline{x})\in\Gamma$, and $\overline{a}\in A$, if
  $\forces{\poP}{\varphi(\overline{a}\checked)}$, then there is a ground 
  $\uniW$ of $\uniV$ \st\ $\overline{a}\in \uniW$ and $\uniW\models\varphi(\overline{a})$. 
\end{xitemize}
The definition of a stronger variant {\It$(\calP,A)_\Gamma$-\RcAp} 
of {$(\calP,A)_\Gamma$-\RcA} is obtained when we replace ``ground'' in \xitemof{x-theintro-1} 
by ``$\calP$-ground''. If $\Gamma=\Lin$, we simply drop $\Gamma$ and talk about 
$(\calP,A)$-\RcA($^+$).

For definability of these axioms in the language of \ZFC, see the paragraphs 
around \Lemmaof{p-hierarchies-4-0}. 

As we have already followed these conventions without explanations, 
we often identify check names of sets with the sets themselves and drop the symbol
``$\checked$''. Also we shall often replace tuples $\overline{a}$ of parameters by a single 
parameter $a$ for simplicity (actually without loss of generality in most of the cases). 

Recurrence Axioms are almost identical with Maximality Principles introduced in 
\cite{hamkins} with the same parameters. 
For $\calP$, $A$ as above, the Maximality Principle for $\calP$ and $A$
($\MP(\calP,A)$ for short) is defined as below. 

For a class $\calP$ of \pos, an $\Lin$-formula $\varphi(\overline{a})$ with parameters
$\overline{a}$ ($\in\uniV$) is said to be a {\It$\calP$-button} if there is $\poP\in\calP$ \st\ 
for any $\poP$-name $\utpoQ$ of \po\ with $\forces{\poP}{\utpoQ\in\calP}$, we have
$\forces{\poP\ast\utpoQ}{\varphi(\overline{a}\checked)}$.

If $\varphi(\overline{a})$ is a $\calP$-button then we call $\poP$ as above a {\It push of 
  the  $\calP$-button $\varphi(\overline{a})$}.

For a class $\calP$ of \pos\ and a set $A$ (of parameters), the {\It Maximality Principle 
  for $\calP$ and $A$} ({\It $\MP(\calP,A)$}, for short)
is the following assertion which is  
formulated as an axiom scheme in $\Lin$: 
\begin{xitemize}
\item[\darkred$\MP(\calP,A)$: ] For any $\Lin$-formula $\varphi(\overline{x})$ and
  $\overline{a}\in A$, if $\varphi(\overline{a})$ is a $\calP$-button then
  $\varphi(\overline{a})$ holds.
\end{xitemize}

Similarly to the restricted versions of Recurrence Axiom, we define,  
for a set $\Gamma$ of $\Lin$-formulas:

\begin{xitemize}
\item[\darkred$\MP(\calP,A)_\Gamma$: ] For any $\varphi(\overline{x})\in\Gamma$ and
  $\overline{a}\in A$, if $\varphi(\overline{a})$ is a $\calP$-button then
  $\varphi(\overline{a})$ holds.
\end{xitemize}

\begin{Prop}{\rm (Barton et al. \cite{5a}, see also Proposition 2.2 in \cite{recurrence})}
  %%\Label{p-intro-0}
  \Label{p-intro-1} Suppose that $\calP$  is a class of \pos\ and $A$ a set 
  (of parameters).\smallskip\\ \wassert{1} $(\calP, A)$-\RcAp\ is equivalent to
  $\MP(\calP,A)$.\smallskip 

  \wassert{2} $(\calP, A)$-\RcA\ is 
  equivalent to the following assertion:  
\begin{xitemize}
\xitem[x-intro-5-0] For any $\Lin$-formula $\varphi(\overline{x})$ and
  $\overline{a}\in A$, if $\varphi(\overline{a})$ is a $\calP$-button then
  $\varphi(\overline{a})$ holds in a ground of\/ $\uniV$. \qed
\end{xitemize}
\end{Prop}
\noindent
See \Lemmaof{p-hierarchies-4} in \sectionof{hierarchies} below and its proof.

Recurrence Axiom ($\Leftrightarrow$ Maximality Principle) can be also characterized as the 
\ZFC\ version of Sy-David Friedman's Inner Model Hypothesis \cite{friedman-sy} (see Barton 
et al \cite{5a}, see also Fuchino and Usuba \cite{recurrence} or Fuchino \cite{janos}). 

In contrast to the proposition above, $(\calP, A)_\Gamma$-\RcAp\ is not necessarily 
equivalent to $\MP(\calP, A)_\Gamma$ for some set $\Gamma$ of formulas. In the next 
section,  
we prove that (under the consistency of certain large cardinal axioms) $\MP(\calP,A)_{\Sigma_2}$ 
does not imply $(\calP,A)_{\Sigma_2}$-\RcA\ and 
$\MP(\calP,A)_{\Pi_2}$ does not imply $(\calP,A)_{\Pi_2}$-\RcA\ (see \Corof{p-hierarchies-7}). 

Later we shall also consider a restricted form of \xitemof{x-intro-5-0} which we will call
$\MP^-(\calP,A)_\Gamma$:

\begin{xitemize}
\item[\darkred$\MP^-(\calP,A)_\Gamma$: ] For any $\varphi(\overline{x})\in\Gamma$ and
  $\overline{a}\in A$, if $\varphi(\overline{a})$ is a $\calP$-button then
  $\varphi(\overline{a})$ holds in a ground of $\uniV$.
\end{xitemize}

Writing $\MP^-(\calP,A)$ for $\MP^-(\calP,A)_{\Lin}$, the assertion of \Propof{p-intro-1},\,\assertof{2} 
is reformulated as  $(\calP,A)$-\RcA\ \ $\Leftrightarrow$\ \ $\MP^-(\calP,A)$. 

While Recurrence Axioms are assertions about the richness of the grounds of the universe
$\uniV$, their characterizations as Maximality Principles may be seen as a variation of 
generic absoluteness. {\ifextended\extendedcolor This is best seen in their further 
  characterization as the principle 
$\MP^*(\calP,A)$ defined around \xitemof{x-YAH-0}, see also \subsectionof{Yah}.  \fi}

The following is an immediate consequence of Bagaria's Absoluteness \Thmof{p-theintro-0}:

\begin{Thm}{\rm (Ikegami and Trang (reformulated for our hierarchy of restricted Recurrence 
    Axioms) \cite{ikegami-trang})}\Label{p-theintro-1} For a 
  (normal) class $\calP$ of \pos\ and a cardinal $\kappa$, \tfae:\smallskip 

  \wassert{a} $(\calP,\calH(\kappa))_{\Sigma_1}$-\RcAp.\smallskip

  \wassert{b} $(\calP,\calH(\kappa))_{\Sigma_1}$-\RcA.\smallskip

  \wassert{c} $\BFA_{\LT\kappa}(\calP)$. 
\end{Thm}
\prf \assertof{a} $\Rightarrow$ \assertof{b}: is trivial.\smallskip

\assertof{b} $\Rightarrow$ \assertof{c}: Assume 
that $(\calP,\calH(\kappa))_{\Sigma_1}$-\RcA\ holds, and suppose that $\poP\in\calP$,
$\varphi$ is a 
$\Sigma_1$-formula in $\Lin$ and $a\in\calH(\kappa)$. By Bagaria's Absoluteness 
\Thmof{p-theintro-0}, it is enough to show 
that $\forces{\poP}{\varphi(a)}$ $\Leftrightarrow$ $\varphi(a)$ holds.

$\forces{\poP}{\varphi(a)}$ $\Leftarrow$ $\varphi(a)$: is clear since $\varphi$ 
is $\Sigma_1$. 

$\forces{\poP}{\varphi(a)}$ $\Rightarrow$ $\varphi(a)$: Assume that
$\forces{\poP}{\varphi(a)}$. By $(\calP,\calH(\kappa))_{\Sigma_1}$-\RcA, there is a ground 
$\uniW$ of $\uniV$ \st\ $a\in\uniW$ and $\uniW\models\varphi(a)$. Since $\varphi$ is
$\Sigma_1$ it follows that $\uniV\models\varphi(a)$. 
\smallskip

\assertof{c} $\Rightarrow$ \assertof{a}: Assume that $\BFA_{\LT\kappa}(\calP)$ holds.
Suppose that $\forces{\poP}{\varphi(a)}$  for
$\poP$, $\varphi$, $a$ as above.  Then, by Bagaria's Absoluteness 
\Thmof{p-theintro-0}, we have $\varphi(a)$. In particular, since $\ssetof{\bbone}\in\calP$ 
(remember the convention set at \xitemof{x-theintro-0-0}),
$\varphi(a)$ holds in a $\calP$-ground of $\uniV$ (namely $\uniV$ itself). 
\qedofThm\qedskip

According to Joel Hamkins \cite{hamkins} it is an observation of his former PhD student 
George Leibman that \MA\ follows from Maximality Principle for $\calP=$ ccc \pos, and the 
set of parameters $\calH(\continuum)$. Actually the corresponding statement had been proved 
much earlier by Stavi and Väänänen \cite{stavi-vaananen} (see also the introduction of 
\cite{stavi-vaananen}). 
This observation is now a part of \Thmof{p-theintro-1}, since \MA\ is 
equivalent to $\BFA_{\LT\continuum}(\calP)$ for this $\calP$. 

Strictly speaking, \Thmof{p-theintro-1} is different from the original theorem in 
Ikegami-Trang \cite{ikegami-trang} (Theorem 1.13 there) in that Ikegami and Trang are 
talking about the 
$\Sigma_n$, $\Pi_n$-hierarchy %% based on the formulation of Maximality Principle (denoted 
%% here as 
$\MP(\cdots)_\Gamma$ for $\Gamma=\Sigma_n$, $\Pi_n$ etc., 
which is shown to be  
different from 
$(\cdots)_{\Gamma}$-\RcAp\ hierarchy (see \Corof{p-hierarchies-7} and the remark 
before the corollary). The proof 
above together with the proof of Theorem 1.13 in \cite{ikegami-trang} actually shows that for a 
normal class of \pos, $(\calP,\calH(\aleph_2))_{\Sigma_1}$-\RcAp\ coincides with
$\MP(\calP,\calH(\aleph_2))_{\Sigma_1}$ (see \Thmof{p-hierarchies-1} and \Corof{p-hierarchies-2}). 

$(\calP,\calH(\aleph_2))_{\Gamma}$-\RcAp\ and 
$\MP(\calP,\calH(\aleph_2))_{\Gamma}$ in general can be different principles. 
We will address to this subtle difference in the next \sectionof{hierarchies}, 
and show that these two hierarchies can split up drastically on the $\Pi_2$ and $\Sigma_2$
levels (see \Corof{p-hierarchies-7}). 

\subsection{(In)compatibility of Recurrence and Maximality with Ground Axiom}
\Label{ground-axiom}
The {\It Ground Axiom} (abbreviation: {\It \GA}) is the axiom asserting that there is no 
proper ground of the  
universe $\uniV$. The axiom is introduced by Joel Hamkins and Jonas Reitz. Its basic 
properties including the formalizability of the axiom in $\Lin$ are proved in Reitz 
\cite{reitz}. 

The relative consistency of \GA\ with \PFA\ is proved in 
\cite{reitz} (see also the proof of \Thmof{p-hierarchies-6} below; actually \GA\ is even 
consistent with $\MM^{++}$, see \Thmof{p-Lg-RcA-1-0}). In particular, this and 
Ikegami-Trang \Thmof{p-theintro-1} imply: 

\begin{Thm}\Label{p-rec-GA-0}
  \GA\ is relatively consistent with $(\calP,\calH(\aleph_2))_{\Sigma_1}$-\RcAp\ for 
  a class of \pos\ $\calP$ whose elements are proper.\qed
\end{Thm}

Since the Recurrence Axiom implies that there are ``many'' different grounds, it is clear 
that \Propabove\ cannot be generalized for $(\cdots)_\Gamma$-\RcA\ for arbitrary $\Gamma$. 
In particular, since Ground Axiom  
itself is formalizable in a $\Pi_3$-sentence in $\Lin$ (see the remark after
\Lemmaof{p-hierarchies-4}), and it is not true in any  
non-trivial generic extension of the ground model, we obtain: 
\begin{Thm}\Label{p-rec-GA-0-a}
  Suppose $(\calP,\emptyset)_{\Sigma_3}$-\RcA\ holds for a non-trivial class $\calP$ of \pos.
  Then \GA\ does not hold.

  $\MP^-(\calP,\emptyset)_{\Sigma_3}$ for a non-trivial $\calP$ also implies $\neg\GA$. 
\end{Thm}
\prf Assume toward a contradiction that $(\calP,\emptyset)_{\Sigma_3}$-\RcA\ holds for a 
non-trivial class $\calP$ of \pos, and \GA\ also holds.
Let $\psi$ be the $\Pi_3$-sentence expressing that the universe does not have a 
non-trivial ground. Let $\poP\in\calP$ be non-trivial forcing. Then
$\forces{\poP}{\neg\psi}$. Since $\neg\psi$ is a $\Sigma_3$-sentence,
$(\calP,\emptyset)_{\Sigma_3}$-\RcA\ implies that there 
is a ground $\uniW$ of $\uniV$ \st\ $\uniW\models\neg\psi$. Since we also assumed 
\GA, $\uniW$ must be identical with $\uniV$. Thus $\uniV\models\neg\psi$. This is a 
contradiction.

Since $\psi$ above is a $\calP$-button, the same proof leads to a contradiction under
$\MP^-(\calP,\emptyset)_{\Sigma_3}$. 
\qedofProp\qedskip

Actually we also have the following delimitation, which shows that \Thmof{p-rec-GA-0} is 
optimal in many instances of $\calP$. 

In the following, we use a variant of the cardinal invariant $\boundingno^*$ introduced in 
Eda-Kada-Yuasa \cite{eda-kada-yuasa}:

\begin{xitemize}
\item[] $\boundingno^{**}:=\min\setof{\kappa\in\Card}{{}
  \begin{array}[t]{@{}l}\mbox{for any }B\subseteq\fnsp{\omega}{\omega}
    \mbox{, if }B\mbox{ is unbounded in }\fnsp{\omega}{\omega}\mbox{ with }\\
    \mbox{respect to }\leq^*\mbox{, then there is }B'\subseteq B\mbox{ with }\\
    \cardof{B'}\leq\kappa
    \mbox{ \st\ }B'\mbox{ is unbounded in }B}.
  \end{array}$
\end{xitemize}

\begin{Lemma}\Label{p-rec-GA-0-0}
  \wassertof{1} $\boundingno=\aleph_1$ can be formulated as a $\Sigma_2$-sentence $\varphi$ 
  in $\Lin$. 
  \smallskip

  \wassertof{2} $\boundingno^{**}=\aleph_1$ can be formulated as a $\Pi_2$-sentence $\psi$ 
  in $\Lin$. \smallskip

  \wassertof{3} $\boundingno<\dominatingno$ can be formulated as a $\Pi_2$-sentence $\eta$ in
  $\Lin$. 
\end{Lemma}
\prf \assertof{1}: The following formula $\varphi$ will do: 
\begin{xitemize}
\xitemx[] $\exists B\ \exists R\ \exists F\ (\,
  \begin{array}[t]{@{}l}
    B\subseteq\fnsp{\omega}{\omega} \land\ 
    \mbox{``}R\mbox{ is an }\omega_1\mbox{-like linear ordering on }B\mbox{ which is}\\ 
    \mbox{witnessed by }F\,\mbox{''}\ \land\ 
    \forall f\,(f\in\fnsp{\omega}{\omega}\ \rightarrow\ (\exists g\in B)\,(g\not<^* f)))\,.
  \end{array}$
\end{xitemize}

\assertof{2}: The following formula $\psi$ will do:

\begin{xitemize}
\xitemx[] $\forall B\ (\,
  \begin{array}[t]{@{}l}
    \obecause{B\subseteq\fnsp{\omega}{\omega}}{}{$\Sigma_0$}\ \land\ 
    \obecause{B\mbox{ is unbounded in }\fnsp{\omega}{\omega}}{}{
      $\Pi_1$ and hence its negation is $\Sigma_1$}\ \ \rightarrow\\
    \exists B'\ \exists R\ \exists F\ (\,
    \begin{array}[t]{@{}l}
      \obecause{B'\subseteq B}{}{$\Sigma_0$}\ \land\
      \obecause{B'\mbox{ is unbounded in }B}{}{$\Sigma_0$}\\
      \land\ {\ubecause{R\mbox{ is }\omega_1\mbox{-like order on }B'
        \mbox{ which is witnessed by }F}{}{$\Sigma_0$}}))\,.
    \end{array}
  \end{array}$
\end{xitemize}

\assertof{3}: ``$\boundingno=\dominatingno$'' is characterized by the existence of a bounding 
family $\subseteq\fnsp{\omega}{\omega}$ which is well ordered \wrt\ $\leq^*$. Similarly to 
above, this can be 
formulated by a $\Sigma_2$-sentence. Hence, the negation of the equality
($\Leftrightarrow$ $\boundingno<\dominatingno$) is $\Pi_2$ in $\Lin$. 
\qedofLemma

\begin{Lemma}\Label{p-rec-GA-0-1}
  $\boundingno\leq\boundingno^{**}\leq\dominatingno$. 
\end{Lemma}
\prf Let $\seqof{f_\alpha}{\alpha<\boundingno}$ be \st\ $f_\alpha\leq^*f_{\alpha'}$ for all
$\alpha<\alpha'<\boundingno$, and $\setof{f_\alpha}{\alpha<\boundingno}$ is unbounded in
$\fnsp{\omega}{\omega}$ (this can be done by letting $\setof{g_\alpha}{\alpha<\boundingno}$ 
be an unbounded subset of $\fnsp{\omega}{\omega}$, and defining $f_\alpha$,
$\alpha<\boundingno$ inductively \st\ we have $g_\alpha\leq^* f_\alpha$ for all
$\alpha<\boundingno$).  Let $B=\setof{f_\alpha}{\alpha<\boundingno}$. Then no 
$B'\subseteq B$ with $\cardof{B'}<\cardof{B}=\boundingno$ is unbounded in $\calB$.

Note 
that the sequence $\seqof{f_\alpha}{\alpha<\boundingno}$ as above also shows that
$\boundingno$ is a regular cardinal.

This proves that $\boundingno\leq\boundingno^{**}$.

To show that $\boundingno^{**}\leq\dominatingno$, suppose that 
$D\subseteq\fnsp{\omega}{\omega}$ is dominating in $\fnsp{\omega}{\omega}$ (\wrt\ $\leq^*$), 
and $\cardof{D}=\dominatingno$.

For any unbounded $B\subseteq\fnsp{\omega}{\omega}$, and for each $d\in D$ let $b_d\in B$ 
be \st\ $b_d\not\leq^* d$. Then $B'=\setof{b_d}{d\in D}\subseteq B$ is unbounded in 
$\fnsp{\omega}{\omega}$ and hence also unbounded in $B$ and is of cardinality
$\leq\dominatingno$. This shows that $\boundingno^{**}\leq\dominatingno$. 
\qedofLemma\qedskip

In the following, we denote with $\poC_\kappa$ the finite support 
$\kappa$-product of Cohen forcing, and with $\poD$ the finite support iteration of Hechler 
forcing of length $\omega_1$. It is easy to see that (over an arbitrary ground model
$\uniV$) we have $\forces{\poC_\kappa}{\boundingno=\aleph_1, \dominatingno\geq\kappa}$ for 
any regular $\kappa\geq\aleph_1$, and 
$\forces{\poD}{\dominatingno=\aleph_1}$. More generally, letting $\poD_\kappa$ be the 
FS-iteration of Hechler forcing of length $\kappa$ for regular $\kappa$, we have
$\forces{\poD_\kappa}{\boundingno=\dominatingno=\kappa}$. 

\begin{Prop}\Label{p-rec-GA-1}
  Suppose $\calP$ is a class of \pos\ with $\poD\in\calP$ and
\ixitem[x-rec-GA-a-0] $\boundingno>\aleph_1$  
  %% $(\calP,\calH(\continuum))_{\Sigma_1}$-\RcA\ 
  holds.\smallskip 

  \wassert{1} If $(\calP,\emptyset)_{\Sigma_2}$-\RcA\ holds, then \GA\ does not 
  hold.\smallskip

  \wassert{2} If $(\calP,\emptyset)_{\Pi_2}$-\RcA\ holds, then \GA\ does not 
  hold.
\end{Prop}
%% \footnotetext{Note that, by Ikegami-Trang \cite{ikegami-trang} (see \Thmof{p-theintro-1} 
%%   below), this condition is  
%%   equivalent with $\BFA_{\LT\continuum}(\calP)$.}
\prf 
Suppose that $\calP$ is as above, and \xitemof{x-rec-GA-a-0} 
%% $(\calP,\calH(\continuum))_{\Sigma_1}$-\RcA\ 
holds. \smallskip

%% \begin{Claim}
%%   $\boundingno>\aleph_1$. 
%% \end{Claim}
%% \prfofClaim Suppose otherwise then there is $B\in\calH(\continuum)$ \st\ 
%% $B\subseteq\fnsp{\omega}{\omega}$ and $B$ is unbounded in $\fnsp{\omega}{\omega}$ \wrt\
%% $\leq^*$. We have $\forces{\poD}{B\mbox{ is bounded}}$. Since the statement 
%% ``$B$ is bounded'' is $\Sigma_1$, there is a ground $\uniW_0$ 
%% of $\uniV$ \st\ $\uniW_0\modelof{B\mbox{ is bounded}}$ by
%% $(\calP,\calH(\continuum))_{\Sigma_1}$-\RcA. This is a  
%% contradiction.\qedofClaim\qedskip

\assertof{1}: Suppose that $(\calP,\emptyset)_{\Sigma_2}$-\RcA\ and \GA\ hold. 

Since we have $\forces{\poD}{\boundingno=\aleph_1}$, and since ``$\boundingno=\aleph_1$'' 
is expressible in  a $\Sigma_2$-sentence in $\Lin$ by \Lemmaof{p-rec-GA-0-0},\,\assertof{1}, it 
follows by $(\calP,\emptyset)_{\Sigma_2}$-\RcA, that there is a ground $\uniW_0$ 
with $\uniW_0\modelof{\boundingno=\aleph_1}$. Since $\uniV=\uniW_0$ by \GA, this is a 
contradiction to \xitemof{x-rec-GA-a-0}.\smallskip

\assertof{2}: Suppose that $(\calP,\emptyset)_{\Pi_2}$-\RcA\ and \GA\ hold. 

Note that by \xitemof{x-rec-GA-a-0} and \Lemmaof{p-rec-GA-0-1}, we have
\ixitem[x-rec-GA-0] $\uniV\models\boundingno^{**}>\aleph_1$. 

Since we have $\forces{\poD}{\dominatingno=\aleph_1}$, we have $\forces{\poD}{\boundingno^{**}=\aleph_1}$
by \Lemmaof{p-rec-GA-0-1}. Since ``$\boundingno^{**}=\aleph_1$'' 
is expressible in  a $\Pi_2$-sentence in $\Lin$ by \Lemmaof{p-rec-GA-0-0},\,\assertof{2}, it 
follows by $(\calP,\emptyset)_{\Pi_2}$-\RcA, that there is a ground $\uniW_0$ 
with $\uniW_0\modelof{\boundingno^{**}=\aleph_1}$. Since $\uniV=\uniW_0$ by \GA, 
it follows that $\uniV\modelof{\boundingno^{**}=\aleph_1}$. 
This is a 
contradiction to \xitemof{x-rec-GA-0}.\smallskip
\qedofProp\qedskip

The following can be proved similarly to \Propof{p-rec-GA-1},\,\assertof{1}. 
\begin{Prop}\Label{p-rec-GA-1-0}
  Suppose that $\calP$ is a class of \pos\ with $\poC_{\aleph_1}\in\calP$ and 
  $\boundingno\geq\aleph_2$ holds. Then $(\calP,\emptyset)_{\Sigma_2}$-\RcA\ implies 
  that $\GA$ does not hold.\qed
\end{Prop}

Note: 
\begin{Lemma}\Label{p-rec-GA-1-1}
  \CH\ can be formulated both as $\Sigma_2$-sentence and $\Pi_2$-sentence in $\Lin$ 
  without parameters.\ifextended\else\qed\fi
\end{Lemma}
{\ifextended\extendedcolor
\prf Consider 
\begin{xitemize}
\xitemx[] $\exists F\ \exists R\ (\,
  \begin{array}[t]{@{}l}
    \smashx{\obecause{R\mbox{ is }\omega_1\mbox{-like ordering on }\fnsp{\omega}{\omega}
      \mbox{ and }R\mbox{ witnesses the }\omega_1\mbox{-likeness}
      }{}{$\Pi_1$: note that we need here the quantification of the sort
        $\forall f\ (f\in\fnsp{\omega}{\omega}\ \rightarrow\ ...)$}}\ 
    )
  \end{array}$
\end{xitemize}
and 
\begin{xitemize}
\xitemx[] $\forall F\ \forall S\ \forall\ A\ (\,
  \begin{array}[t]{@{}l}
    (\,\smash{\obecause{S=\mbox{``\,}\fnsp{\omega}{\omega}\mbox{\,''}\ \land\ 
         A\mbox{ is an ordinal}\ \land\ \mapping{F}{S}{A}\mbox{ is a surjection}
    }{}{$\Pi_1$, so the negation is $\Sigma_1$}}\,)\\
    \rightarrow\ \ubecause{A<\omega_2}{}{$\Sigma_1$}\ ).\\[-1em]
    \hfill{\rm\qedofLemma}
  \end{array}$
\end{xitemize}\fi}

\begin{Prop}\Label{p-rec-GA-1-1-0}
\wassertof{1}  Suppose that $\neg\CH$ holds and $\calP$ contains a \po\ collapsing
$\continuum$ to $\aleph_1$ without adding reals. Then each of 
$(\calP,\emptyset)_{\Sigma_2}$-\RcA\ and $(\calP,\emptyset)_{\Pi_2}$-\RcA\ implies $\neg\GA$. 
\smallskip

\wassert{2} Suppose that $\CH$ holds and $\calP$ contains a \po\ $\poQ$ 
adding $\GE\aleph_2$ reals without collapsing cardinals $\LE\aleph_2$. Then 
each of $(\calP,\emptyset)_{\Sigma_2}$-\RcA\ and $(\calP,\emptyset)_{\Pi_2}$-\RcA\ implies $\neg\GA$.\smallskip

\wassert{3} Suppose that $\calP$ contains sufficiently many ccc \pos\ (containing enough 
CS-iterations of  
Cohen and Hechler \pos\ would suffice), then each of $(\calP,\emptyset)_{\Sigma_2}$-\RcA\ and
$(\calP,\emptyset)_{\Pi_2}$-\RcA\ implies $\neg\GA$. 
\end{Prop}
\prf Similarly to the proof of \Propof{p-rec-GA-1} using \Lemmaof{p-rec-GA-1-1}. 
For \assertof{3}, we consider cases where (a) $\aleph_1=\boundingno=\dominatingno$, (b)
$\aleph_1<\boundingno=\dominatingno$, or (c) $\aleph_1<\boundingno<\dominatingno$, and 
apply  \Lemmaof{p-rec-GA-0-0} in all of the cases. 
\qedofProp\qedskip

\subsection{Incompatibility of Laver genericity with Ground Axiom}\Label{Laver-GA}
In the following we want to discuss the impact of the results we obtained above on axioms 
stating that there is a Laver-generic large cardinal (Laver-gen.\ large cardinal axioms). 

The strongest variant of Laver-generic large cardinal axiom which has been considered so far, is 
the tightly super-$C^{(\infty)}$ $\calP$-Laver-generically hyperhuge cardinal (see Fuchino 
and Usuba \cite{recurrence}).

Here, a 
cardinal $\kappa$ is said to be ({\It tightly}, resp.) {\It 
  super-$C^{(\infty)}$ $\calP$-Laver-generically hyperhuge} if for all $n\in\natnums$ and 
for any $\lambda_0>\kappa$ and $\poP\in\calP$, there are $\lambda\geq\lambda_0$ with  
$V_\lambda\prec_{\Sigma_n}\uniV$, and $\poP$-name $\utpoQ$ with
$\forces{\poP}{\utpoQ\in\calP}$ \st\ for $(\uniV,\poP\ast\utpoQ)$-generic $\genH$, there 
are $j$, $M\subseteq\uniV[\genH]$ with 
$\Elembed{j}{\uniV}{M}{\kappa}$, $j(\kappa)>\lambda$, $j\imageof{j(\lambda)}$, $\poP, \poP\ast\utpoQ,\genH\in M$, 
${V_{j(\lambda)}}^{\uniV[\genH]}\prec_{\Sigma_n}\uniV[\genH]$ (and $\poP\ast\utpoQ$ is of 
size $\leq j(\kappa)$,  
resp.).%% \addtocounter{footnote}{-1}
\footnote{When we say ``a \po\ $\poP$ is of 
  cardinality $\leq\mu$'' we actually mean that there is a \po\ $\poQ$ forcing equivalent 
  to $\poP$ \st\ $\cardof{\poQ}\leq\mu$. } 

We can also define (tightly) super-$C^{(\infty)}$ $\calP$-Laver gen.\ large cardinal 
analogously  for notions of large cardinal other than hyperhugeness (see \cite{recurrence} 
or \cite{janos}). For an iterable class $\calP$ of \pos\ which also permits transfinite 
iteration with some suitable support, we can prove that the existence of the tightly super
$C^{(\infty)}$ $\calP$-Laver gen.\ hyperhuge cardinal is consistent under a 2-huge 
cardinal (\cite{recurrence}). 

Note that we cannot formulate the (genuine) large cardinal property corresponding to 
(tightly) super $C^{(\infty)}$ $\calP$-Laver-generically large cardinal in $\Lin$. However, for a 
natural class $\calP$ of \pos\ like proper \pos, semiproper \pos, ccc \pos, etc.\ we can 
formulate the notion of (tightly) super $C^{(\infty)}$ $\calP$-Laver-generically 
hyperhugeness in an axiom scheme in $\Lin$. 
This is because $\calP$-Laver-generically hyperhugeness of a cardinal $\kappa$ implies
$\kappa=\kappa_\refl$ ($:=\max\ssetof{\aleph_2,\continuum}$) for these classes $\calP$ of 
\pos, and hence we can formulate the 
(tightly)  
super $C^{(\infty)}$ $\calP$-Laver-generically hyperhugeness of such $\kappa$ in infinitely 
many formulas without introducing a new constant symbol for the cardinal. 

In Fuchino and Usuba \cite{recurrence}, it is proved that if $\kappa$ is tightly 
super $C^{(\infty)}$ $\calP$-Laver-generically ultrahuge, 
then $(\calP,\calH(\kappa))$-\RcAp\ holds. Here the tightly 
super $C^{(\infty)}$ $\calP$-Laver-generically ultrahugeness is 
apparently much weaker than tightly 
super $C^{(\infty)}$ $\calP$-Laver-generically hyperhugeness. 

Note that 
$(\calP,\calH(\kappa))$-\RcAp\ is also 
an assertion formalizable only in infinitely many formulas. In contrast, it is proved in 
\cite{future}, that, in a sense, Laver-genericity without ``super $C^{(\infty)}$'' details never 
implies the full $(\calP,\calH(\kappa))$-\RcAp. 

By the result mentioned above and by \Thmof{p-rec-GA-0-a}, it follows immediately that:

\begin{Prop}\Label{p-rec-GA-1-1-1}
  For any iterable class $\calP$ of \pos, if $\kappa$ is tightly 
  super $C^{(\infty)}$ $\calP$-Laver-generically ultrahuge, then $\GA$ does not hold.\qed
\end{Prop}

In \cite{recurrence}, it is proved that if $\kappa$ is tightly $\calP$-generically 
hyperhuge (not necessarily Laver-generic) then there is the bedrock (i.e.\ the ground 
satisfying \GA) and $\kappa$ is hyperhuge in the bedrock.

On the other hand, it is shown in 
\cite{janos} that a tightly $\calP$-Laver-generically ultrahuge cardinal for nice iterable
$\calP$ is $\leq\kappa_\refl$.
Here an iterable class $\calP$ of \pos\ is said to be {\It nice} if either $\calP$ preserves 
$\omega_1$ and $\Col(\omega_1,\ssetof{\omega_2})\in\calP$, or $\calP$ contains a \po\ which 
adds a new real. Actually the following lemma is one of the main rationales of the definition of 
the cardinal $\kappa_\refl$. 
\begin{Lemma}\Label{p-rec-GA-1-1-2}
  Suppose that $\calP$ is a nice iterable class of \pos. If $\kappa$ is $\calP$-Laver-gen.\ 
  supercompact, then $\kappa\leq\kappa_\refl$. 
\end{Lemma}
\prf By Lemma 6.,\assertof{2} and \assertof{3} in \cite{janos}. \qedofLemma\qedskip

Thus we obtain:
\begin{Prop}\Label{p-rec-GA-1-2}
  For a nice iterable class $\calP$ of \pos, suppose that there is a tightly $\calP$-Laver-generically 
  hyperhuge cardinal. Then the bedrock exists and it is different from $\uniV$. 
  In particular, \GA\ does not hold.
\end{Prop}
\prf By \Lemmaof{p-rec-GA-1-1-2}, the tightly $\calP$-Laver-generically 
  hyperhuge cardinal is $\leq\kappa_\refl$ (in $\uniV$) while $\kappa$ is hyperhuge in the 
  bedrock $\ol{\uniW}$. This implies that $\uniV\not=\ol{\uniW}$. In particular, \GA\ does not 
  hold. \qedofProp\qedskip

At the moment we do not know if the existence of a tightly $\calP$-generically 
hyperhuge in theorem in 
\cite{recurrence} mentioned above can be weakened to the existence of some tight generic 
large cardinal of lower consistency strength. 
However, in \cite{janos}, it is proved that for an iterable class $\calP$ of \pos, if $\kappa$ is 
tightly $\calP$-Laver-gen.\ ultrahuge then $(\calP,\calH(\kappa))_{\Sigma_2}$-\RcAp\ holds 
(Theorem 21 in \cite{janos}). 
Note that ultrahuge cardinal is apparently much weaker than hyperhuge cardinal. 

\begin{Thm}\Label{p-rec-GA-2}
  Suppose that $\calP$ is an iterable class of \pos\ 
  satisfying one of the conditions in \Propof{p-rec-GA-1-1-0}.

  If $\kappa_\refl$ is tightly $\calP$-Laver-gen.\ ultrahuge then \GA\ does not hold. 
\end{Thm}
\prf If $\kappa$ is tightly $\calP$-Laver-gen.\ ultrahuge then
$(\calP,\calH(\kappa))_{\Sigma_2}$-\RcAp\ holds by Theorem 21 in \cite{janos}. Thus 
if $\poD\in\calP$ then (by \Propof{p-genabs-Laver-a} below and) by 
\Propof{p-rec-GA-1},\assertof{1}, it follows that \GA\ does not  
hold.

Other cases can treated similarly by applying other assertions of \Propof{p-rec-GA-1-1-0}. \\
\qedofThm

\section{Hierarchies of restricted Recurrence Axioms and Maximality Principles}
\Label{hierarchies}

{\ifextended\extendedcolor
Ikegami and Trang \cite{ikegami-trang} formulated 
Maximality Principle slightly different from our 
notation. Their Maximality Principle
 %% (which we denote with ``$\MPIT$''\footnote{This is a 
 %%  local defintion only introduced to formulate the following \Lemmaof{p-hierarchies-0}.}) 
in restricted form 
is defined for a class $\calP$ of \pos, cardinal $\kappa$ and set $\Gamma$ of formulas (in
$\Lin$) as:
\begin{xitemize}
%% \item[{\darkred$\MPIT(\calP,\kappa)_\Gamma$: }]
\xitem[x-hierarchies-a] 
  For any formula $\varphi\in\Gamma$ and
  $A\subseteq\kappa$, if $\varphi(A)$ is a $\calP$-button then $\varphi(A)$ holds in $V$.
\end{xitemize}

Since (tuples of) elements of $\calH(\kappa^+)$ can be coded by subsets of $\kappa$ we have:
\begin{LemmaA}\Label{p-hierarchies-0}\extendedcolor
  For any class of \pos\ $\calP$, cardinal $\kappa$, and a set $\Gamma$ of $\Lin$-formulas, 
  we have:\\ The Maximality Principle of Ikegami and Trang \xitemof{x-hierarchies-a}\ \
  $\Leftrightarrow$\ \ $\MP(\calP,\calH(\kappa^+))_\Gamma$.
  \qed
\end{LemmaA}

Thus, Ikegami and Trang's Theorem (Theorem 1.13 in \cite{ikegami-trang}) is reformulated as:
\begin{xitemize}
\xitem[x-hierarchies-0] $\MP(\calP,\calH(\aleph_2))_{\Sigma_1}$\ \ $\Leftrightarrow$\ \
  $\BFA_{\LT\aleph_2}(\calP)$\quad for any (normal) class $\calP$ of \pos. 
\end{xitemize}
\fi}

\ifextended
{\extendedcolor Actually, almost the same argument as the proof of \Thmof{p-theintro-1} 
  given above, we can  
  show the following more general theorem:}\else
The proof of \Thmof{p-theintro-1} can be reused almost without any change to prove the following 
generalization of Theorem 1.13 in Ikegami-Trang \cite{ikegami-trang}:
\fi

\begin{Thm}{\rm (Generalization of the original Ikegami-Trang Theorem)}\Label{p-hierarchies-1}
  For a (normal) class $\calP$ of \pos, and a cardinal $\kappa$ \tfae:\smallskip

  \wassert{a} $\MP(\calP,\calH(\kappa))_{\Sigma_1}$.\smallskip

  \wassert{b} $\MP^-(\calP,\calH(\kappa))_{\Sigma_1}$.\smallskip

  \wassert{c} $\BFA_{\LT\kappa}(\calP)$. \ifextended\else\qed\fi
\end{Thm}
{\ifextended\extendedcolor\prf \assertof{a} $\Rightarrow$ \assertof{b}: is 
  trivial,\smallskip

  \assertof{b} $\Rightarrow$ \assertof{c}: Assume $\MP^-(\calP,\calH(\kappa))_{\Sigma_1}$. 
  By Bagaria's Absoluteness \Thmof{p-theintro-0}, it is enough to show that
  $\varphi(\overline{a})$ $\Leftrightarrow$ $\forces{\poP}{\varphi(\overline{a})}$ for all 
  $\poP\in\calP$, $\Sigma_1$-formula $\varphi$ in $\Lin$ 
  and $\overline{a}\in\calH(\kappa)$.

  If $\varphi(\overline{a})$ holds then, since $\varphi$ is $\Sigma_1$,
  $\forces{\poP}{\varphi(\overline{a})}$ also holds.

  Suppose that $\forces{\poP}{\varphi(\overline{a})}$. Then, for any $\poP$-name $\utpoQ$ 
  of \po\ $\forces{\poP}{\forces{\utpoQ}{\varphi(\overline{a})}}$ since $\varphi$ is
  $\Sigma_1$. In particular $\varphi(\overline{a})$ is a $\calP$-button with the push $\poP$ 
  of the button. By $\MP^-(\calP,\calH(\kappa))_{\Sigma_1}$, there is a ground $\uniW_0$ of 
  $\uniV$ \st\ $\overline{a}\in\uniW_0$ and $\uniW_0\models\varphi(\overline{a})$. Since 
  $\varphi$ is $\Sigma_1$, it follows that $\uniV\models\varphi(\overline{a})$.

  \assertof{b} $\Rightarrow$ \assertof{c}: Assume $\BFA_{\LT\kappa}(\calP)$. 
  By Bagaria's Absoluteness \Thmof{p-theintro-0}, this means that
  $\varphi(\overline{a})$ $\Leftrightarrow$ $\forces{\poP}{\varphi(\overline{a})}$ for all 
  $\poP\in\calP$, $\Sigma_1$-formula $\varphi$ in $\Lin$ 
  and $\overline{a}\in\calH(\kappa)$.

  Suppose that $\forces{\poP}{\forces{Q}{\varphi(\overline{a})}\mbox{ for all }Q\in\calP}$ 
  for $\poP$, $\varphi$, $\overline{a}$ as above. Since $\ssetof{\bbone}\in\calP$ it 
  follows that $\forces{\poP}{\varphi(\overline{a})}$. By assumption, it follow that 
  $\varphi(\overline{a})$ holds. 
\qedofThm
\fi}

\begin{Cor}\Label{p-hierarchies-2}
  For a class $\calP$ of \pos\ and for an infinite cardinal $\kappa$, we have\vspace{-1\smallskipamount}
  \begin{xitemize}
  \item[] 
    $\MP(\calP,\calH(\kappa))_{\Sigma_1}$,\ \ $\Leftrightarrow$\ \
    $\MP^-(\calP,\calH(\kappa))_{\Sigma_1}$,\ \ $\Leftrightarrow$\ \  
    $\BFA_{\LT\kappa}(\calP)$, \\ $\Leftrightarrow$\ \ 
    $(\calP,\calH(\kappa))_{\Sigma_1}$-\RcA,\ \ $\Leftrightarrow$\ \
    $(\calP,\calH(\kappa))_{\Sigma_1}$-\RcAp. 
  \end{xitemize}

\end{Cor}
\prf By \Thmof{p-theintro-1} and \Thmof{p-hierarchies-1}. \qedofCor
\qedskip

The following lemma holds since $\Pi_1$-formulas are downward absolute. 

\begin{Lemma}\Label{p-hierarchies-3}
  For any class $\calP$ of \pos, and any set $A$, $(\calP,A)_{\Pi_1}$-\RcAp\ and
  $\MP(\calP,A)_{\Pi_1}$ hold (in \ZFC). In particular, we 
  have \vspace{-1\smallskipamount}
  \begin{xitemize}
  \item[] 
    $(\calP,A)_{\Pi_1}$-\RcA\ \ $\Leftrightarrow$\ \ 
    $(\calP,A)_{\Pi_1}$-\RcAp\ \ $\Leftrightarrow$\ \ $\MP^-(\calP,A)_{\Pi_1}$
    \ \ $\Leftrightarrow$\ \ $\MP(\calP,A)_{\Pi_1}$\,.\\
    \qed
  \end{xitemize}
\end{Lemma}

In the following, we show that the equivalence in \Lemmaabove\ does not hold for $\Pi_2$. 

Nevertheless, we have the following implications. %% This is obtained 
%% by examining the proof of \Propof{p-intro-1},\,\assertof{1} (e.g.\ the one given in 
%% \cite{recurrence}):   

\begin{Lemma}\Label{p-hierarchies-4}
  Suppose that $\calP$ is a (normal) class of \pos\ defined by a $\Sigma_m$-formula 
  without parameters  for some number $m$, and $A$ a set.\footnotemark \smallskip

  \wassertof{1} $(\calP,A)_{\Pi_n}$-\RcAp\ \ $\Rightarrow$\ \ 
  $\MP(\calP,A)_{\Pi_n}$, for all $n\geq\max\ssetof{m,1}$.\smallskip

  \wassertof{2} $\MP(\calP,A)_{\Sigma_n}$\ \ $\Rightarrow$\ \ $(\calP,A)_{\Sigma_n}$-\RcAp, 
  for all $n\geq\max\ssetof{m,3}$. 
\end{Lemma}
\footnotetext{Note that ``$x$ is ccc\ \po'', ``$x$ is proper \pos'', ``$x$ is 
  semi-proper \po" are all $\Sigma_2$-statements. In case of ``$x$ is (semi-)proper \po'', 
  this can be seen in the formulation: 
  %% \begin{xitemize}
  %% \item[] 
  %%   ``$\exists K\, (K\mbox{ is a 
  %%     cardinal }\geq(\beth_\omega)^+(\cardof{x})\ \land\ \cdots)$''
  %% \end{xitemize}
  %% and in that this can be further 
  %% reformulated as: 
  \begin{xitemize}
    \item[] 
%%  \xitem[x-hierarchies-1] 
    ``$\exists \ul{\kappa}\, \exists F\, (\ul{\kappa}\mbox{ is a 
      cardinal }\land\,F\mbox{ ``codes'' the fact ``}
    \ul{\kappa}\geq(\beth_\omega)^+(\cardof{x})\mbox{''}\ \land\ \cdots)$.''
  \end{xitemize}
  Here, the underline to $\kappa$ is added to suggest that the symbol does not denote a 
  constant symbol but rather  a variable 
  in the language $\Lin$. We shall keep this convention in the following. }
\prf The following proofs are  just re-examinations of the easy proof of 
\Propof{p-intro-1},\,\assertof{1} (e.g.\ the one given in Fuchino and Usuba
\cite{recurrence}).\smallskip 

\assertof{1}: Note that, for $n=1$, the claim also follows from \Lemmaof{p-hierarchies-3}.
%% Hence we may
%% assume that we have
%% \begin{xitemize}
%% \xitem[x-hierarchies-1-a] 
%%   $n\geq\max\ssetof{m,2}$.
%% \end{xitemize}

Assume that $(\calP,A)_{\Pi_n}$-\RcAp\ holds for
$n\geq\max\ssetof{m,2}$. To show that $\MP(\calP,A)_{\Pi_n}$ holds, suppose that
$\varphi=\varphi(\ol{x})$ is 
a $\Pi_n$-formula, $\ol{a}\in A$, and $\poP\in\calP$ is \st\
$\forces{\poP}{\forall P\in\calP\,(\forces{P}{\varphi(\ol{a})})}$ holds in $\uniV$.

``\,$\forall P\in\calP\,(\forces{P}{\varphi(\ol{x})})$'' is 
$\Pi_n$ by the choice of $n$.  Let us denote this formula by $\varphi^*$. Thus, we have 
$\forces{\poP}{\varphi^*(\ol{a})}$. 

By $(\calP,A)_{\Pi_n}$-\RcAp, it follows 
that there is a $\calP$-ground $\uniW$ of $\uniV$ \st\ $\ol{a}\in\uniW$ and
$\uniW\models\varphi^*(\ol{a})$. By the definition of $\varphi^*$, and since $\uniW$ is a
$\calP$-ground, it follows that $\uniV\models\varphi(\ol{a})$.\smallskip

\assertof{2}: Assume that $\MP(\calP,A)_{\Sigma_n}$ holds. Suppose that $\varphi$ is
$\Sigma_n$-formula, $\ol{a}\in A$, and $\poP\in\calP$ is \st\
\begin{xitemize}
\xitem[x-hierarchies-1-0] 
  $\forces{\poP}{\varphi(\ol{a})}$. 
\end{xitemize}

Then we have
$\forces{\poP}{\varphi(\ol{a})\mbox{ holds in a }\calP\mbox{-ground}}$.

The assertion 
\begin{xitemize}
\xitem[x-hierarchies-2] ``$\varphi(\ol{x})$ holds in a $\calP$-ground''
\end{xitemize} 
can be expressed in a $\Sigma_n$-formula $\varphi^{**}=\varphi^{**}(\ol{x})$ 
(see the remark after the proof of the present lemma). By 
\xitemof{x-hierarchies-1-0} and by the definition \xitemof{x-hierarchies-2} of
$\varphi^{**}$ we have $\forces{\poP\ast\utpoQ}{\varphi^{**}(\ol{a})}$ for all $\poP$-names 
$\utpoQ$ with $\forces{\poP}{\utpoQ\in\calP}$. Thus, by $\MP(\calP,A)_{\Sigma_n}$, it follows 
that $\uniV\models\varphi^{**}(\ol{a})$. By the definition of $\varphi^{**}$, it follows 
that there is a $\calP$-ground $\uniW_0$ \st\ $\uniW_0\models\varphi(\ol{a})$. 
\qedofLemma\qedskip

The fact that \xitemof{x-hierarchies-2} can be formulated in a
$\Sigma_n$-formula for $n\geq\max\ssetof{m,3}$, can be seen as follows.

First, let us notice the fact formulated in \Lemmaof{p-hierarchies-4-0} below.

In the following, for formulas $\varphi$ and $\psi=\psi(x,\ol{y})$ we
denote by $\varphi^{\psi(x;\ol{y})}$ the formula $\varphi$ restricted to $\psi(x,\ol{y})$
where $\psi(x,\ol{y})$ is thought to be the definition of the class
$\calA_{\ol{y}}=\setof{x}{\psi(x,\ol{y})}$ with parameters (or, more precisely,  place
holders for parameters) $\ol{y}$. The semi-colon in ``$\varphi^{\psi(x;\ol{y})}$'' should
remind this allocation of roles among the free variables of $\psi$.

Thus  
$\varphi^{\psi(x;\ol{y})}$ corresponds to the informal statement: $\calA_{\ol{y}}\models\varphi$.
\begin{Lemma}\Label{p-hierarchies-4-0}
  Suppose that $\psi=\psi(x,\ol{y})$ is a $\Sigma_m$-formula and $\varphi$ is
  $\Sigma_n$-formula (\/$\Pi_n$-formula, resp.). Then $\varphi^{\psi(x;\ol{y})}$ is a
  $\Sigma_{max\ssetof{m,n}}$-formula (a $\Pi_{max\ssetof{m,n}}$-formula,
  resp.).
  %% \ifextended\else\qed\fi
\end{Lemma}
{%% \ifextended\extendedcolor
\prf For quantifier free formula $\varphi$, the claim of the Lemma is true since
$\varphi^\psi=\varphi$.

Suppose that, for a $\Sigma_n$-formula $\varphi_0=\varphi_0(x_0,\ol{x})$ ($\Pi_n$-formula
$\varphi_1=\varphi_1(x_0,\ol{x})$), 
${\varphi_0}^{\psi(x;\ol{y})}$ is $\Sigma_k$ (${\varphi_1}^{\psi(x;\ol{y})}$ is $\Pi_k$
resp.).

Then
\begin{xitemize}
\item[] $(\forall x_0)(\psi(x_0,\ol{y})\rightarrow{\varphi_0}^{\psi(x;\ol{y})})$ is
  $\Pi_{\max\ssetof{m,k+1}}$, and 
\item[] $(\exists x_0)(\psi(x_0,\ol{y})\land {\varphi_1}^{\psi(x;\ol{y})})$ is
  $\Sigma_{\max\ssetof{m,k+1}}$.
\end{xitemize}

Using this fact, the claim of the Lemma can be proved now by induction on $n$. 
\qedofLemma\qedskip%% \fi
}

An examination of \cite{BHTU} and \cite{reitz} reveals a construction of a $\Pi_2$-formula\\
$\Phi(x,\ul{P},\ul{\delta},r,\ul{G})$
which says that 
\begin{xitemize}
\xitemx[] $\ul{P}$ is a \po, $\ul{\delta}$ is a regular cardinal in $\uniV$, 
  there is a uniquely determined inner model $\calM$ with $\ul{\delta}$-cover and
  $\ul{\delta}$-approximation properties  
  \st\ $r=(\fnsp{\ul{\delta}\GT}{2})^\calM$, $\ul{P}\in\calM$, $\calM\not=\uniV$, $\ul{G}$ is an
  $(\calM,\ul{P})$-generic set \st\ $\uniV=\calM[\ul{G}]$, and $x\in\calM$.
\end{xitemize}

Let $\psi=\psi(x)$ be a $\Sigma_m$-formula expressing ``$x\in\calP$''. Then, for a
$\Sigma_n$-formula $\varphi=\varphi(\ol{x})$ for $n\geq\max\ssetof{m,3}$, the formula
$\varphi^{**}(\ol{x})$ defined as
\begin{xitemize}
\xitemx[] $\exists\,\ul{P}\exists\,\ul{\delta}\,\exists\,r\,\exists\,\ul{G}\,
  (\Phi(\emptyset,\ul{P},\ul{\delta},r,\ul{G})\land \psi^{\Phi(x;\cdots)}(\ul{P})
  \land \varphi^{\Phi(x;\cdots)}(\ol{x}))$
\end{xitemize}
is $\Sigma_n$ by \Lemmaof{p-hierarchies-4-0}, and $\varphi^{**}(\ol{a})$ expresses 
``$\varphi(\ol{a})$ holds in a $\calP$-ground''. 

We shall also use the following variant of Maximality Principle. Let $\calP$, $A$, $\Gamma$ 
be as before.

\begin{xitemize}
\item[{\darkred $\MP^+(\calP,A)_\Gamma$ :}] For $\varphi\in\Gamma$, and $\ol{a}\in A$, if
  $\varphi(\ol{a})$ is a $\calP$-button, then $\ssetof{\bbone}$ is a push of the $\calP$-button
  $\varphi(\ol{a})$. 
\end{xitemize}

As before, we drop the subscript $\Gamma$ from $\MP^+(\calP,A)_\Gamma$ if $\Gamma=\Lin$.
\begin{Lemma}\Label{p-hierarchies-5}
  \wassertof{1} $\MP^+(\calP,A)_\Gamma$\ \ $\Rightarrow$\ \ $\MP(\calP,A)_\Gamma$. 
  \smallskip

  \wassert{2} {\rm (Hamkins \cite{hamkins})} For an iterable $\calP$, 
  $\MP^+(\calP,A)$\ \ $\Leftrightarrow$\ \ $\MP(\calP,A)$. More precisely, if $\calP$ is
  $\Sigma_m$-definable, then for any $n\geq\max\ssetof{m,1}$, we have  
  $\MP^+(\calP,A)_{\Pi_n}$\ \ $\Leftrightarrow$\ \ $\MP(\calP,A)_{\Pi_n}$.\smallskip\memo{??????}

  \wassert{3} For an iterable $\calP$, if\/ $\MP^+(\calP,A)_\Gamma$, then for any $\poP\in\calP$, we 
  have 
  \begin{xitemize}
  \item[] 
    $\forces{\poP}{\MP^+(\calP,A)_\Gamma}$. 
  \end{xitemize}
\end{Lemma}
\prf \assertof{1}: is clear by definition.\smallskip

\assertof{2}: By \assertof{1}, it is enough to show ``$\Leftarrow$'', Assume that
$\MP(\calP,A)$ holds, and suppose that $\varphi(\ol{a})$ is a $\calP$-button for 
an $\Lin$-formula $\varphi$ and $\ol{a}\in A$. 
Then $\varphi^*:=\forall Q\,(Q\in\calP\ \rightarrow\ \forces{Q}{\varphi(\ol{a})})$ is a
$\calP$-button. Hence, by $\MP(\calP,A)$, $\varphi^*$ holds in $\uniV$. But this means that
$\ssetof{\bbone}$ is a push for the button $\varphi$. \smallskip

\assertof{3}: Suppose that $\MP^+(\calP,A)_\Gamma$ holds (in $\uniV$).
For $\varphi\in\Gamma$, and $\ol{a}\in A$, let $\poP\in\calP$ be \st\ it forces that
$\varphi(\ol{a})$ is a $\calP$-button. By Maximal Principle (of forcing), there is 
a $\poP$-name $\utpoQ$ of a \po\ \st\
\begin{xitemize}
\item[] 
  $\forces{\poP}{\utpoQ\in\calP
    \ \land\ \forces{\utpoQ}{\forall R\,(R\in\calP\ \rightarrow\ 
      \forces{R}{\varphi(\ol{a})})}}$.
\end{xitemize}

Since $\calP$ is iterable, it follows that $\varphi(\ol{a})$ is a $\calP$-button over
$\uniV$. Thus, by $\MP^+(\calP,A)_\Gamma$, $\ssetof{\bbone}$ is a push of 
the $\calP$-button 
$\varphi(\ol{a})$ (in $\uniV$). 

Again since $\calP$ is iterable, it follows that
$\forces{\poP}{\ssetof{\bbone}\mbox{ is a push of the }\calP\mbox{-button }\varphi(\ol{a})}$. 
\qedofLemma\qedskip

The following should be folklore:\memo{\cite{recurrence}:p-Lg-RcA-0-0}

\begin{Lemma}\Label{p-hierarchies-5-0} \wassertof{1} 
  If $\alpha$ is a limit ordinal and $V_\alpha$ satisfies a sufficiently large finite 
  fragment of \ZFC, then for any $\poP\in V_\alpha$ and $(\uniV,\poP)$-generic $\genG$, we 
  have $V_\alpha[\genG]={V_\alpha}^{\uniV[\genG]}$.\smallskip 

  \wassert{2}  If $\alpha$ is a limit ordinal and $V_\alpha$ satisfies a sufficiently large finite 
  fragment of \ZFC, then for any direct limit $\poP$ of an iteration of length
  $\On^{V_\alpha}$ in $V_\alpha$, which is definable and preserving 
  cardinals in $V_\alpha$, if $\genG$ is $(\uniV,\poP)$-generic, then we 
  have
  \begin{xitemize}
    \xitemx[] 
    $V_\alpha[\genG]={V_\alpha}^{\uniV[\genG]}$. 
  \end{xitemize}
  \wassert{3} For each natural number $k$, there is a sufficiently large $k'>k$ \st\ for any
  $\alpha\in\On$ if $V_\alpha\prec_{\Sigma_{k'}}\uniV$ (i.e.\ $\alpha$ 
  is $\Sigma_{k'}$-correct), then for any \po\ $\poP\in V_\alpha$ 
  and $(\uniV,\poP)$-generic
  $\genG$, ${V_\alpha}^{\uniV[\genG]}\prec_{\Sigma_k}\uniV[\genG]$.\smallskip  

  \wassert{4} \memox{!!!!!!!} Suppose 
  that $\seqof{\poP_\alpha,\utpoQ_\alpha}{\alpha\in\On}$ is an Easton support class 
  iteration of increasingly directed closed \pos\ and $\clpoP$ is the class direct limit of the 
  iteration. If $k$ is a natural number and $\kappa$ is a regular cardinal which is
  $\Sigma_{k'}$-correct for a sufficiently large $k'>k$, then we have
  ${V_\kappa}^{\uniV[\genG_\kappa]}\prec_{\Sigma_{k}}\uniV[\clgenG]$ for 
  any $(\uniV,\clpoP)$-generic $\clgenG$ and $\genG_\kappa=\clgenG\cap\poP_\kappa$. 
\end{Lemma}
\prf \assertof{1}: See e.g.\ Lemma 3.2 in \cite{recurrence}. \quad
\assertof{2}: follows from \assertof{1}. \smallskip

\assertof{3}: see the proof of Lemma 4.8,\,\assertof{1} in the extended version of 
\cite{recurrence}.\smallskip

\assertof{4}: Let $\Phi=\Phi(x)$  be an $\Lin$-formula which defines $\clpoP$. Then by the 
choice of $\kappa$, we have $\poP_\kappa=\Phi^{{V_\kappa}^\uniV}$. The claim \assertof{4}
follows from this fact with an argument practically identical to that for \assertof{3}. 
\qedofLemma\qedskip 

By \Thmof{p-rec-GA-0-a}, we cannot replace $\Pi_2$ in the next theorem by $\Sigma_3$. 
\begin{Thm}\Label{p-hierarchies-6}
  Suppose 
  that $\calP$ is a $\Sigma_2$-definable iterable class of \pos\ containing 
  all $\sigma$-closed \pos%% \ \st\ the property ``$x\in\calP$'' is absolute for sufficiently 
  %% directed closed forcings
  , and that $\MP^+(\calP,\calH(\kappa_\refl))_{\Pi_2}$  holds. 
  Suppose further that there is a proper class $\calK$ of supercompact cardinals. %% \st\ 
  %% \begin{xitemize}
  %% \xitem[x-hierarchies-2-a] 
  %%   all $\kappa\in\calK$ are $\Sigma_k$-correct for a sufficiently large natural 
  %%   number $k$. 
  %%   %% $\calK$ is stationary in $\On$. 
  %%   %% ;\qquad and 
  %%   %% \xitem[x-hierarchies-2-a-0] 
  %%   %%   $\calK\cap\Lim(\calK)$ is cofinal in $\On$. 
  %% \end{xitemize}
  %% $\calK$ is stationary in $\On$ (i.e., it intersects with any (with parameters definable) 
  %% club class $\subseteq\On$). 

  If $\clpoP$ is the class 
  \po\ for  Laver preparation for $\calK$ (see the proof below for more details), then we have
  \begin{xitemize}
    \xitemx[] 
    $\forces{\clpoP}{{}
      \begin{array}[t]{@{}l}
        \GA\ +\ \MP(\calP,\calH(\kappa_\refl))_{\Pi_2}\\
        +\ \mbox{there are class many
          supercompact cardinals}}.
      \end{array}$ 
  \end{xitemize}
\end{Thm}
%% \footnotetext{Recall that a cardinal $\kappa$ is $\Sigma_k$-correct if 
%%   $V_\kappa\prec_{\Sigma_k}\uniV$ holds. The assumption of the theorem follows from a 
%%   tightly  
%%   super $C^{(k)}$ 
%%   $\calP$-Laver gen.\ hyperhuge cardinal for a sufficiently large $k$. However we 
%%   conjecture that the consistency 
%%   strength of the assumption is much lower than the existence of this generic large 
%%   cardinal.} 
\prf Let $\calP$, $\calK$ be as above.

Let $\mapping{f}{\On}{\uniV}$ be a universal Laver function for $\calK$. I.e., a class 
function $f$ \st\ 
\begin{xitemize}
\xitem[x-hierarchies-2-0] 
  for any $\kappa\in\calK$, we 
  have $\mapping{f\restr\kappa}{\kappa}{V_\kappa}$, and for any $x\in\uniV$ and any
  $\lambda\geq\max\ssetof{\kappa,\cardof{\trcl(x)}}$, there is a normal ultrafilter 
  $\calU_{\kappa,\lambda,x}$ over $\Pkl{}{}$ and associated elementary embedding
  $\Elembed{j_{\kappa,\lambda,x}}{\uniV}{M}{\kappa}$ 
  with $j_{\kappa,\lambda,x}(f)(\kappa)=x$. 
\end{xitemize}
We may also assume that 
\begin{xitemize}
\xitem[x-hierarchies-2-1] 
  $f(\alpha)=0$ for all
  $\alpha<\kappa_\refl$. 
\end{xitemize}
Note that $f\restr\kappa$, $\kappa\in\calK$ are uniformly definable across $V_\kappa$
(${}=\calH(\kappa)$) for all $\kappa\in\calK$.

A universal Laver function exists (see e.g.\ Apter \cite{apter-Laver-Indestr}, Lemma 1).

Let $\seqof{\poP_\alpha}{\alpha\in\On}$ be the Laver preparation along with $f$ making 
supercompactness of all $\kappa\in\calK$ indestructible by $\kappa$-directed closed 
forcing.{\ifextended\extendedcolor\footnote{\extendedcolor Note that $\kappa$-directed 
    closed means $\kappa$-directed closed.}\fi} I.e.,
$\seqof{\poP_\alpha}{\alpha\in\On}$ is defined as the iterative part of the  
Easton support{\ifextended\footnote{\extendedcolor I.e. direct limit at $\poP_\alpha$ for regular $\alpha$ 
    and inverse ($\approx$ full support) limit at singular $\alpha$.}\fi} iteration
$\seqof{\poP_\alpha, \utpoQ_\alpha}{\alpha\in\On}$ with a control sequence 
$\seqof{\lambda_\alpha}{\alpha\in\On}$ of cardinals defined recursively by
\begin{xitemize}
\xitem[x-hierarchies-3] If $\alpha\in\On$ is a limit and closed \wrt\
  $\seqof{\lambda_\beta}{\beta<\alpha}$, $f(\alpha)=\pairof{\utpoQ,\lambda}$ with
  $\forces{\poP_\alpha}{\utpoQ\mbox{ is }\alpha\mbox{-directed closed \po}}$, then
  $\utpoQ_\alpha=\utpoQ$ and $\lambda_\alpha=\lambda$;
\xitem[x-hierarchies-4] Otherwise $\lambda_\alpha=\sup\setof{\lambda_\beta}{\beta<\alpha}$ 
  and $\forces{\poP_\alpha}{\poQ_\alpha=\ssetof{\bbone}}$. 
\end{xitemize}

Let $\clpoP$ be the class forcing which is the direct limit of
$\seqof{\poP_\alpha}{\alpha\in\On}$. For each 
$\kappa\in\calK$, let $\clpoP_{\GT\kappa}$ be (the class $\calP_\kappa$-name of) 
the 
$\kappa$-directed closed tail part of the 
iteration. Thus $\clpoP\sim \poP_\kappa\ast\clpoP_{\GT\kappa}$ and
$\forces{\poP_\kappa}{\clpoP_{\GT\kappa}\mbox{ is a }\kappa\mbox{-directed closed class \po}}$. 

Let $\clgenG$ be $(\uniV, \clpoP)$-generic. For each $\kappa\in\calK$, let
$\genG_\kappa=\clgenG\cap\poP_\kappa$. 
Each $\kappa\in\calK$ is made indestructible under $\kappa$-directed closed forcing by $\poP_\kappa$ (see 
e.g. \cite{laver-indestr}). In particular, $\kappa$ remains supercompact in
$\uniV[\genG_\lambda]$ for 
all $\lambda\in\calK$. It follows that $\kappa$ remains supercompact also 
in $\uniV[\clgenG]$.

$\uniV[\clgenG]\models\GA$. This is because $\uniV[\clgenG]$ satisfies Continuum Coding Axiom 
(\CCA), and \GA\ follows from it (see \cite{reitz} Theorem 3.2). That $\uniV[\clgenG]$ satisfies 
\CCA\  
follows from the fact that in $\uniV[\clgenG]$ there are cofinally many indestructible 
supercompact cardinals. {\ifextended\extendedcolor[\![The first author learned the 
      following from G.\,Goldberg: \memo{\mbox{Scan\_2024-05-14--12.42 gappo - annotated} p.52} 
      Suppose that $a\subseteq\alpha$. Let $\kappa>\alpha$ be an indestructible 
      supercompact cardinal. Let $\poP$ be $\kappa$-directed closed \po\ \st\
      $\forces{\poP}{\forall \beta<\alpha\,( 2^{\aleph_{\kappa+\beta+1}}=(\aleph_{\kappa+\beta+1})^+
         \ \leftrightarrow\ \beta\in a)}$. Since
      $\forces{\poP}{\kappa\mbox{ is supercompact}}$ by assumption and a supercompact 
      cardinal is $\Sigma_2$-correct. It follows that
      \begin{xitemize}
      \item[] $\forces{\poP}{V_\kappa\models\exists\delta\,\forall \beta<\alpha\,
        (2^{\aleph_{\delta+\beta+1}}=(\aleph_{\delta+\beta+1})^+
         \ \leftrightarrow\ \beta\in a)}$
      \end{xitemize}
      Since $\forces{\poP}{V_\kappa=(V_\kappa)^\uniV}$ by $\kappa$-directed closedness of
      $\poP$, it follows that 
      \begin{xitemize}
      \item[] 
        $\uniV\modelof{V_\kappa\models\exists\delta\,\forall \beta<\alpha\,
        (2^{\aleph_{\delta+\beta+1}}=(\aleph_{\delta+\beta+1})^+
        \ \leftrightarrow\ \beta\in a)}$. 
      \end{xitemize}
      Hence the assertion holds in $\uniV$. 
]\!] \fi}

Thus, it is enough to show 
that $\uniV[\clgenG]\models\MP(\calP,\calH(\kappa_\refl))_{\Pi_2}$. 
Note that we have 
\begin{xitemize}
\xitem[x-hierarchies-5] 
  $\calH(\kappa_\refl)^\uniV=\calH(\kappa_\refl)^{\uniV[\genG_\kappa]}=\calH(\kappa_\refl)^{\uniV[\clgenG]}$ 
\end{xitemize}
by $\min(\calK)$-directed closedness of $\clpoP$, and 
\xitemof{x-hierarchies-2-1}. Also 
\begin{xitemize}
\xitem[x-hierarchies-6] ${V_\kappa}^{\uniV[\genG_\kappa]}={V_\kappa}^{\uniV[\clgenG]}$ 
\end{xitemize}
%% \quad and 
%% \xitem[x-hierarchies-7]
%%   $\calH(\kappa)^{\uniV[\genG_\kappa][\geng]}=\calH(\kappa)^{\uniV[\genG][\geng]}$    
%% \end{xitemize}
for all $\kappa\in\calK$
by $\kappa$-directed closedness of $\clpoP_{\GT\kappa}$ and 
\Lemmaof{p-hierarchies-5-0},\,\assertof{2}.  

Working in $\uniV[\clgenG]$, 
suppose that $\varphi=\varphi(\ol{x})$ is a $\Pi_2$-formula and
$\ol{a}\in\calH(\kappa_\refl)$ ($=\calH(\kappa_\refl)^\uniV$).
%%  --- henceforth we drop the superscript $\uniV$, $\uniV[\clgenG]$ etc.\ 
%% to ``$\calH(\kappa_\refl)$'' in virtue of \xitemof{x-hierarchies-5}. 

Further in $\uniV[\clgenG]$, suppose that $\poS\in\calP$ is \st\
\begin{xitemize}
\xitem[x-hierarchies-7-a-a] 
  $\uniV[\clgenG]\models\forces{\poS}{\forall T\in\calP\,(\forces{T}{\varphi(\ol{a})})}$. 
\end{xitemize}
We want to show that $\varphi(\ol{a})$ holds (in $\uniV[\clgenG]$). 

By replacing $\uniV$ by $\uniV[\genG_{\kappa_0}]$, and $\calK$ 
by $\calK\setminus\kappa_0+1$ for a large enough $\kappa_0\in\calK$ with
$\poS\in \uniV[\genG_{\kappa_0}]$, we   
may assume that $\poS\in {V_{\kappa_0}}^\uniV$ for a $\kappa_0<\min(\calK)$. Let $\geng$ be
$(\uniV[\clgenG],\poS)$-generic. Since 
$\poS\in\calP$ and since $\calP$ is iterable, we have 
\begin{xitemize}
\xitem[x-hierarchies-7-a] 
  $\uniV[\clgenG][\geng]\models\forall T\in\calP\,(\forces{T}{\varphi(\ol{a})})$.
\end{xitemize}

By \xitemof{x-hierarchies-6} (and \Lemmaof{p-hierarchies-5-0},\,\assertof{2}), we have
\begin{xitemize}
\xitem[x-hierarchies-6-0] 
  ${V_\kappa}^{\uniV[\genG_\kappa][\geng]}={V_\kappa}^{\uniV[\clgenG][\geng]}$ 
\end{xitemize}
for all $\kappa\in\calK$.

Note that each $\kappa\in\calK$ remains supercompact in all of $V[\genG_\kappa]$, 
  $\uniV[\genG_\kappa][\geng]$, $V[\clgenG]$, and $\uniV[\clgenG][\geng]$.  
Thus%% By \xitemof{x-hierarchies-7-a-0}
,

%% Let $\kappa\in\calK$ be a limit of $\calK\cap C^{(k)}$ for a large enough $k$
%% \st\ 
\begin{xitemize}
\xitem[x-hierarchies-7-0-a] 
  ${V_\kappa}^{\uniV[\clgenG]}\prec_{\Sigma_2}\uniV[\clgenG]$, 
\xitem[x-hierarchies-7-0] 
  ${V_\kappa}^{\uniV[\clgenG][\geng]}\prec_{\Sigma_2}\uniV[\clgenG][\geng]$,
\xitem[x-hierarchies-7-a-0-0]
  ${V_\kappa}^{\uniV[\genG_\kappa]}\prec_{\Sigma_2}\uniV[\genG_\kappa]$,\quad and
\xitem[x-hierarchies-7-a-1]
  ${V_\kappa}^{\uniV[\genG_\kappa][\geng]}\prec_{\Sigma_2}\uniV[\genG_\kappa][\geng]$ 
\end{xitemize}
for all $\kappa\in\calK$. 
By \xitemof{x-hierarchies-6}, and \xitemof{x-hierarchies-7-0-a}, we have 
\begin{xitemize}
\xitem[x-hierarchies-7-a-2] 
  ${V_\kappa}^{\uniV[\genG_\kappa]}\prec_{\Sigma_2}\uniV[\clgenG]$ 
\end{xitemize}
for all $\kappa\in\calK$. Similarly 
\begin{xitemize}
\xitem[x-hierarchies-7-0-0] 
  ${V_\kappa}^{\uniV[\genG_\kappa][\geng]}\prec_{\Sigma_2}\uniV[\clgenG][\geng]$
\end{xitemize}
holds for all
$\kappa\in\calK$ 
by \xitemof{x-hierarchies-6-0} and \xitemof{x-hierarchies-7-0}, 

``$\forall T\in\calP\,(\forces{T}{\varphi(\ol{a})})$'' is $\Pi_2$
(note that we need $\Sigma_2$-definability of $\calP$ for this). 
Hence
${V_\kappa}^{\uniV[\genG_\kappa][\geng]}\models\forall T\in\calP\,(\forces{T}{\varphi(\ol{a})})$
by \xitemof{x-hierarchies-7-a} and \xitemof{x-hierarchies-7-0-0}. By \xitemof{x-hierarchies-7-a-1}, 
it follows that
$\uniV[\genG_\kappa][\geng]\models\forall T\in\calP\,(\forces{T}{\varphi(\ol{a})})$. 
This implies that $\varphi(\ol{a})$ is a $\calP$-button in $\uniV[\genG_\kappa]$. 

Similarly, since $\uniV[\clgenG]\models\poS\in\calP$, and $\calP$ is $\Sigma_2$, we have
$\uniV[\genG_\kappa]\models\poS\in\calP$ by \xitemof{x-hierarchies-7-a-2} and
\xitemof{x-hierarchies-7-a-0-0}. 

Since we have $\uniV[\genG_\kappa]\models\MP^+(\calP,\calH(\kappa_\refl))_{\Pi_2}$ by 
\Lemmaof{p-hierarchies-5},\,\assertof{3}, it follows that
$V[\genG_\kappa]\models\varphi(\ol{a})$ for all $\kappa\in\calK$. 

Since $\varphi$ is $\Pi_2$ and $\calK$ is cofinal in $\On$, it follows that
$\uniV[\clgenG]\models\varphi(\ol{a})$ by \xitemof{x-hierarchies-7-a-2}.

This shows that
$\uniV[\clgenG]\models\MP(\calP,\calH(\kappa_\refl))_{\Pi_2}$ holds. 
%% Note that,  by \xitemof{x-hierarchies-2-a} and \Lemmaof{p-hierarchies-5-0},\,\assertof{4}, 
%% we have 
%% \begin{xitemize}
%% \xitem[x-hierarchies-8] ${V_\kappa}^{\uniV[\genG_\kappa][\geng]}\prec_{\Sigma_{k'}}
%% \uniV[\clgenG][\geng]$ for still large enough $k'$ for all $\kappa\in\calK$. 
%% \end{xitemize}
%% \memo{There is such 
%% $\kappa$ by \xitemof{x-hierarchies-2-a} and by Lévy-Montague Reflection Theorem. ?}
%% 
%% It follows that 
%% ${V_\kappa}^{\uniV[\clgenG][\geng]}\models\forall T\in\calP\,(\forces{T}{\varphi(\ol{a})})$  
%% By \xitemof{x-hierarchies-6}, it 
%% follows that ${V_\kappa}^{\uniV[\genG_\kappa][\geng]}\models\forall 
%% T\in\calP\,(\forces{T}{\varphi(\ol{a})})$.
%% 
%% Thus, by ****
%% \vspace{10zh}
%% 
%% It follows that $\uniV[\genG_\kappa][\geng]\models\forall 
%% T\in\calP\,(\forces{T}{\varphi(\ol{a})})$. By \Lemmaof{p-hierarchies-5},\,\assertof{3}, we 
%% have $V[\genG_\kappa]\models\MP^+(\calP,\calH(\kappa_\refl))_{\Pi_3}$. It follows that 
%% $V[\genG_\kappa]\models\varphi(\ol{a})$.
%% 
%% Since $\varphi$ is $\Pi_3$,  each $\kappa\in\calK$ is $\Sigma_2$-correct in $V[\clgenG]$, and
%% since $\calK$ is cofinal in $\On$, 
%% it follows that $\uniV[\clgenG]\models\varphi(\ol{a})$. 
\qedofThm\qedskip

If we start from a ground model with a proper class $\calK$ of $C^{(n)}$-supercompact cardinals 
(see Bagaria \cite{bagaria-Cn}) for sufficiently large $n$, we can improve the condition
``$\calP$ is $\Sigma_2$-definable'' in \Thmof{p-hierarchies-6} by ``$\calP$ is
$\Sigma_3$-definable'' (see the remark after the proof of \Propof{p-hierarchies-6-1}). 

%----------------------------------------------------------------------------------------asldk↓
``$\Pi_2$'' and  ``$\MP(\calP,\calH(\kappa_\refl))_{\Pi_2}$'' in \Thmof{p-hierarchies-6} can be
also replaced by ``$\Delta_3$'' and ``$\MP^*(\calP,\calH(\kappa_\refl))_{\Delta_3}$'' which is not
covered by $\MP(\calP,\calH(\kappa_\refl))_{\Pi_2}$ (\Propof{p-hierarchies-6-1}).  

In the following, we quickly review the definition and some needed facts about the variation
$\MP^*(\calP,A)_{\Gamma}$ of the Maximality Principle which will be further studied in
Gappo and Lietz \cite{future2}.\quad

For an iterable class $\calP$ of \pos, a set $A$ of parameters, and a set $\Gamma$ 
of $\Lin$-formulas, let 
{\darkred$(\Gamma)^*_\calP$} be the set of provably $\calP$-persistent formulas in $\Gamma$.
That is, the collection of all formulas $\varphi\in\Gamma$, $\varphi=\varphi(\ol{x})$ \st\ 
the $\Lin$-sentence 
\begin{xitemize}
\xitem[x-YAH-a] ${\darkred(\varphi)^*_\calP}
  :=\forall\ol{x}\,(\varphi(\ol{x})\rightarrow\forall\vsymb{\poP}\in\calP\,(\forces{\vsymb{\poP}}{\varphi(\ol{x})}))$. 
\end{xitemize}
is provable in \ZFC.

Note that if $\Gamma$ is closed \wrt\ equivalence (which is provable in \ZFC) and has a recursive 
representatives (modulo the equivalence), then the same holds for $(\Gamma)^*_\calP$ (as far 
as $\calP$ is a definable class but this is always assumed). 

Now {\darkred$\MP^*(\calP,A)_\Gamma$} is defined as the axiom scheme consisting of formulas of the form 
\begin{xitemize}
\xitemsub[x-YAH-0]{\varphi}
$(\forall\vsymb{\poP}\in\calP)\,(\forall\ol{x}\in A)\,(\forces{\vsymb{\poP}}{\varphi(\ol{x})}\ 
\rightarrow\ \varphi(\ol{x}))$
\end{xitemize}
for each $\varphi\in(\Gamma)^*_\calP$.

Similarly to the $\MP(\cdots)_\Gamma$ and $(\cdots)_\Gamma$-\RcAp\ hierarchies, we write
$\MP^*(\calP,A)$ for $\MP^*(\calP,A)_{\Lin}$. $\MP^*(\calP,\calH(\continuum))$ of the 
family $\calP$ consisting of ccc \pos\ corresponds to the principles considered in 
Stavi and Väänänen 
\cite{stavi-vaananen}. 

The main point of the definition of $\MP^*(\calP,A)_\Gamma$ is that the persistence is 
coded in the collection 
$(\Gamma)^*_\calP$ so that each formula \xitemsubof{x-YAH-0}{\varphi} in
$\MP^*(\calP,A)_\Gamma$ remains at about the same complexity of $\varphi$. 

\memo{\mbox{Scan\_2024-05-14--12.42 gappo - annotated} p.96}
\ifextended\else The following is easy to prove. 
\fi
\begin{Lemma}\Label{p-YAH-0} For an iterable class $\calP$ of \pos, arbitrary set $A$ of 
  parameters and set $\Gamma$ of formulas, \tfae: \smallskip\\
  \wassertof{a} $\MP(\calP,A)_{(\Gamma)^*_\calP}$,\quad
  \wassertof{b} $\MP^+(\calP,A)_{(\Gamma)^*_\calP}$,\quad 
  \wassertof{c} $(\calP,A)_{(\Gamma)^*_\calP}$-\RcAp\,,\\
  \wassertof{d} $\MP^*(\calP,A)_{\Gamma}$. 
  \ifextended\else\qed\fi
\end{Lemma}
{\ifextended\extendedcolor \prf We show \assertof{a} $\Leftrightarrow$ \assertof{d}. Other 
  equivalences can be proved similarly.\smallskip

  \assertof{a} $\Rightarrow$ \assertof{d}: Assume that $\MP(\calP,A)_{(\Gamma)^*_\calP}$ 
  holds. Suppose that $\varphi\in(\Gamma)^*_\calP$, $\ol{a}\in A$ 
  and $\forces{\poP}{\varphi(\ol{a})}$ holds for a $\poP\in\calP$. Then, since
  $\varphi\in(\Gamma)^*_\calP$, 
  $\varphi(\ol{a})$ is a $\calP$-button and $\poP$ is its push. By
  $\MP(\calP,A)_{(\Gamma)^*_\calP}$ it follows that $\uniV\models\varphi(\ol{a})$. This 
  shows that $\MP^*(\calP,A)_{\Gamma}$ holds. 
  \smallskip

  \assertof{d} $\Rightarrow$ \assertof{a}: Assume that $\MP^*(\calP,A)_{\Gamma}$ holds and 
  suppose that $\varphi\in(\Gamma)^*_\calP$, $\ol{a}\in A$, $\varphi(\ol{a})$ is 
  a $\calP$-button and $\poP\in\calP$ is its push. Then, 
  since $\ssetof{\bbone}\in\calP$, $\forces{\poP}{\varphi(\ol{a})}$.  
  By $\MP^*(\calP,A)_{\Gamma}$ it follows that $\uniV\models\varphi(\ol{a})$. This shows 
  that $\MP(\calP,A)_{(\Gamma)^*_\calP}$ holds. \qedofLemma\qedskip\fi}

The hierarchy of this type of restricted Maximality Principles also appears in Goodman 
\cite{goodman}. Our $\MP^*(\calP,A)_\Gamma$
is called ``$\Gamma$-$\MP_\calP(A)$'' in \cite{goodman}. 
Though the choice of symbols in \cite{goodman} is so that letter $\Gamma$ is used to denote the 
class of \pos\ and $\Phi$ to denote the class of formulas. 

The hierarchy of $\MP^*$ is actually a special 
case of the hierarchy of ``$\BFA(\calP,\Gamma)_{\kappa,\lambda}$'' in Asperó \cite{aspero}
(see \cite{future2}{\ifextended\extendedcolor: see also \Lemmaof{p-Yah-2}\fi}).\footnote{These principles are related but different from the Bounded
  Forcing Axioms $\BFA_{\LT\kappa}(\calP)$.}
%----------------------------------------------------------------------------------------asldk↑

The proof of \Lemmaof{p-hierarchies-4},\,\assertof{2}, shows also the implication:
\begin{xitemize}
\item[\assertof{2'}] $\MP^*(\calP,A)_{\Sigma_n}$\ \ $\Rightarrow$\ \ $(\calP,A)_{\Sigma_n}$-\RcAp, 
  for all $n\geq\max\ssetof{m,3}$ where $\calP$ is $\Sigma_m$. 
\end{xitemize}
Thus, by \Thmof{p-rec-GA-0-a}, $\MP^*(\calP,A)_{\Sigma_n}$ for $n\geq\max\ssetof{m,3}$ for 
$m$ as above implies $\neg\GA$. 
In particular, for $\Sigma_3$-definable $\calP$, $\MP^*(\calP,A)_{\Sigma_3}$ implies
$\neg\GA$. 
This shows that the condition $\Delta_3$ in the following proposition is (almost) optimal.

%% by $\MP(\calP,\calH(\kappa_\refl))_{\Gamma}$ for the class
%%     $\Gamma$ of formulas defined as in \: 
\begin{Prop}\Label{p-hierarchies-6-1}
  Suppose 
  that $\calP$ is a $\Sigma_3$-definable iterable class of \pos\ containing 
  all $\sigma$-closed \pos%% \ \st\ the property ``$x\in\calP$'' is absolute for sufficiently 
  %% directed closed forcings
  , and that $\MP^+(\calP,\calH(\kappa_\refl))_{\Delta_3}$  holds. 
  Suppose further that there is a proper class $\calK$ of $C^{(n)}$-supercompact cardinals
  for a sufficiently large $n$ 
  and\/ $\clpoP$ is the class \po\ defined as in the proof of 
  \Thmof{p-hierarchies-6} for this $\calK$.

  Then we have\vspace{-1.08ex}
  \begin{xitemize}
    \xitemx[] 
    $\forces{\clpoP}{{}
        \GA\ +\ \MP^*(\calP,\calH(\kappa_\refl))_{\Delta_3}}$. 
  \end{xitemize}
\end{Prop}
\prf Suppose that $\calK$ and $\clpoP$ are as above and
$\MP^+(\calP,\calH(\kappa_\refl))_{\Delta_3}$ holds. 

Let $\clgenG$ be 
a $(\uniV,\clpoP)$-generic filter.

As it has been already shown in the proof of \Thmof{p-hierarchies-6}, we have
$\uniV[\clgenG]\models\GA$. 

So we prove $\uniV[\clgenG]\models\MP^*(\calP,\calH(\kappa_\refl))_{\Delta_3}$. 
Working in $\uniV[\clgenG]$, 
suppose that $\varphi=\varphi(\ol{x})$ is a $(\Delta_3)^*_\calP$-formula and
$\ol{a}\in\calH(\kappa_\refl)$ ($=\calH(\kappa_\refl)^\uniV$).

Further in $\uniV[\clgenG]$, suppose that $\poS\in\calP$ is \st\
\begin{xitemize}
\xitem[x-hierarchies-8-0] 
  $\uniV[\clgenG]\models\forces{\poS}{\varphi(\ol{a})}$. 
\end{xitemize}
We want to show that $\varphi(\ol{a})$ holds (in $\uniV[\clgenG]$). 

Similarly to the proof of \Thmof{p-hierarchies-6}, we   
may assume that $\poS\in {V_{\kappa_0}}^\uniV$ for a $\kappa_0<\min(\calK)$. Let $\geng$ be
$(\uniV[\clgenG],\poS)$-generic. By  
the choice  \xitemof{x-hierarchies-8-0} of $\poS$, we have 
\begin{xitemize}
\xitem[x-hierarchies-8-1] 
  $\uniV[\clgenG][\geng]\models\varphi(\ol{a})$.
\end{xitemize}

By the choice of the ``sufficiently large'' $n$ (in terms of 
\Lemmaof{p-hierarchies-5-0},\,\assertof{4} and \assertof{3}), and by \xitemof{x-hierarchies-6},  
%% In Goldberg \cite{goldberg}, it is shown that, starting from a class $\calK$ of 
%% inaccessible $\Sigma_2$-correct cardinals, there is a class \po\ $\clpoP$ constructed
%% as the direct limit of an 
%% Easton support iteration of increasingly directed closed \pos\ which forces that 
%% elements of $\calK$ become indestructible inaccessible $\Sigma_2$-correct cardinals.
%% In \cite{goldberg} it 
%% is also shown that extendible cardinals in $\calK$ are preserved by $\clpoP$. Since 
%% extendible cardinals are $\Sigma_3$-correct, obtain 
\begin{xitemize}
\xitem[x-hierarchies-8-2] 
  ${V_\kappa}^{\uniV[\genG_\kappa]}\prec_{\Sigma_3}\uniV[\clgenG]$, 
\xitem[x-hierarchies-9] 
  ${V_\kappa}^{\uniV[\genG_\kappa][\geng]}\prec_{\Sigma_3}\uniV[\clgenG][\geng]$, 
\end{xitemize}
and, since $\kappa$ remains supercompact in $\uniV[\genG_\kappa]$ and $\uniV[\genG_\kappa][\geng]$, 
\begin{xitemize}
\xitem[x-hierarchies-9-0] 
    ${V_\kappa}^{\uniV[\genG_\kappa]}\prec_{\Sigma_2}\uniV[\genG_\kappa]$, 
\xitem[x-hierarchies-10]
  ${V_\kappa}^{\uniV[\genG_\kappa][\geng]}\prec_{\Sigma_2}\uniV[\genG_\kappa][\geng]$ 
\end{xitemize}
for all $\kappa\in\calK$.  

By \xitemof{x-hierarchies-8-2}, $\poS\in\calP$ (in $\uniV[\clgenG]$), and since
$\calP$ is $\Sigma_3$, ${V_\kappa}^{\uniV[\genG_\kappa]}\models\poS\in\calP$. Hence, by 
\xitemof{x-hierarchies-9-0}, $\uniV[\genG_\kappa]\models\poS\in\calP$.

Since $\varphi$ is $\Delta_3$, \xitemof{x-hierarchies-8-1} and \xitemof{x-hierarchies-9} 
implies ${V_\kappa}^{\uniV[\genG_\kappa][\geng]}\models\varphi(\ol{a})$. This and 
\xitemof{x-hierarchies-10} imply $\uniV[\genG_\kappa][\geng]\models\varphi(\ol{a})$.

By
$\varphi\in(\Delta_3)^*_\calP$, it follows that
$\uniV[\genG_\kappa][\geng]\models\forall\ul{\poQ}\in\calP\,(\forces{\ul{\poQ}}{\varphi(\ol{a})})$.
Thus $\uniV[\genG_\kappa]\models\exists\ul{\poS}\in\calP\,(
\forall\ul{\utpoQ}\,(\forces{\ul{\poS}}{\ul{\utpoQ}\in\calP\rightarrow\forces{\ul{\utpoQ}}{\varphi(\ol{a})}}))$. 

Since $\uniV[\genG_\kappa]\models\MP^+(\calP,\calH(\kappa_\refl))_{\Delta_3}$ by 
\Lemmaof{p-hierarchies-5},\,\assertof{3} (and since $\ssetof{\bbone}\in\calP$), it follows that
\begin{xitemize}
\xitem[x-hierarchies-11] 
  $\uniV[\genG_\kappa]\models\varphi(\ol{a})$. 
\end{xitemize}

Since $\varphi$ is $\Delta_3$, and hence $\Pi_3$ in particular,
${V_{\kappa}}^{\uniV[\genG_\kappa]}\models\varphi(\ol{a})$ by 
\xitemof{x-hierarchies-9-0}.
Thus, by 
\xitemof{x-hierarchies-11}, 
it follows that $\uniV[\clgenG]\models\varphi(\ol{a})$. 

This shows that
$\uniV[\clgenG]\models\MP^*(\calP,\calH(\kappa_\refl))_{\Delta_3}$ holds. 
\qedofProp\qedskip

The first half of the proof of \Propof{p-hierarchies-6-1} can be applied to the proof of 
\Thmof{p-hierarchies-6} to obtain: 
\begin{itemize}
\item[]\hspace{-2.4em}{\bf (A variant of \Thmof{p-hierarchies-6})}\quad\it 
   Suppose 
  that $\calP$ is a $\Sigma_3$-definable iterable class of \pos\ containing 
  all $\sigma$-closed \pos%% \ \st\ the property ``$x\in\calP$'' is absolute for sufficiently 
  %% directed closed forcings
  , and that $\MP^+(\calP,\calH(\kappa_\refl))_{\Pi_2}$  holds. 
  Suppose further that there is a proper class $\calK$ of $C^{(n)}$-supercompact cardinals
  for a sufficiently large $n$. %% \st\ 

  If $\clpoP$ is the class 
  \po\ for  Laver preparation for $\calK$, then we have
  \begin{xitemize}
    \xitemx[] 
    $\forces{\clpoP}{{}
      \begin{array}[t]{@{}l}
        \GA\ +\ \MP(\calP,\calH(\kappa_\refl))_{\Pi_2}\\
        +\ \mbox{there are class many
          supercompact cardinals}}.
      \end{array}$ 
  \end{xitemize}
  {\ifextended\extendedcolor (see the proof below for more details)\fi}
\end{itemize}

The following theorem is also obtained by combining the proofs of \Thmof{p-hierarchies-6} and
\Propof{p-hierarchies-6-1} taking \Lemmaof{p-YAH-0} into account. 
\begin{Thm}\Label{p-hierarchies-6-2}
  \wassertof{1} Suppose 
  that $\calP$ is a $\Sigma_2$-definable iterable class of \pos\ containing 
  all $\sigma$-closed \pos\ \st\ %% \ the property ``$x\in\calP$'' is absolute for sufficiently 
  %% directed closed forcings
  $\MP^+(\calP,\calH(\kappa_\refl))_{\Gamma}$  holds where 
  $\Gamma$ denotes here the set of all formulas representable as the conjunction of 
  a $\Sigma_2$-formula and a $(\Delta_3)^*_\calP$-formula.  
  Suppose further that there is a proper class $\calK$ of supercompact cardinals 
  and\/ $\clpoP$ is the class \po\ defined as in the proof of 
  \Thmof{p-hierarchies-6} for this $\calK$.

  Then we have\vspace{-1.08ex}%% -0.6*1.8=-1.08
  \begin{xitemize}
    \xitemx[] 
    $\forces{\clpoP}{{}
        \GA\ +\ \MP^*(\calP,\calH(\kappa_\refl))_{\Gamma}}$. 
  \end{xitemize}
  \wassertof{2} Suppose 
  that $\calP$ is a $\Sigma_3$-definable iterable class of \pos\ containing 
  all $\sigma$-closed \pos\ \st\ 
  %% \ the property ``$x\in\calP$'' is absolute for sufficiently 
  %% directed closed forcings
  $\MP^+(\calP,\calH(\kappa_\refl))_{\Gamma}$  holds where 
  $\Gamma$ is the set of formulas defined as in {\assertof{1}}. 
  Suppose further that there is a proper class $\calK$ of $C^{(n)}$-supercompact cardinals
  for 
  a sufficiently large $n$, 
  and\/ $\clpoP$ is the class \po\ defined as in the proof of 
  \Thmof{p-hierarchies-6} for this $\calK$.

  Then we have\vspace{-1.08ex}%% -0.6*1.8=-1.08
  \begin{xitemize}
    \xitemx[] 
    $\forces{\clpoP}{{}
        \GA\ +\ \MP^*(\calP,\calH(\kappa_\refl))_{\Gamma}}$. \qed
  \end{xitemize}
\end{Thm}

The following Corollary shows in particular that 
the implication in \Lemmaof{p-hierarchies-4}, \assertof{1} for $n=2$ cannot be 
reversed. %% , and 
%% ``$\max\ssetof{m,3}$'' in \Lemmaof{p-hierarchies-4}, \assertof{2} is optimal.

\begin{Cor}\Label{p-hierarchies-7} \wassertof{1}
  Suppose 
  that $\calP$ is a $\Sigma_2$-definable iterable class of \pos\ containing 
  all $\sigma$-closed \pos, and also a \po\ adding a real%%  \st\ the property ``$x\in\calP$'' 
  %% is absolute for sufficiently directed closed forcings
  . 
  Assume further that there is a proper class $\calK$ of supercompact cardinals, and
  $\MP^+(\calP,\calH(\kappa_\refl))_{\Gamma}$ holds where $\Gamma$ is defined just as in
  \Thmof{p-hierarchies-6-2}. Then, there is a class \po\/ $\clpoP$  
  \st\ we have\vspace*{-1.4ex}  
  \begin{xitemize}
  \xitemx[] $\forces{\clpoP}{\neg\,(\calP,\emptyset)_{\Pi_2}\mbox{-\RcA, }
    \neg\,(\calP,\emptyset)_{\Sigma_2}\mbox{-\RcA\ and\/ }
    \MP(\calP,\calH(\kappa_\refl))_{\Gamma}}$.
  \end{xitemize}

  \wassert{2}
  Suppose 
  that $\calP$ is a $\Sigma_3$-definable iterable class of \pos\ containing 
  all $\sigma$-closed \pos, and a \po\ adding a real%%  \st\ the property ``$x\in\calP$'' is 
  %% absolute for sufficiently directed closed forcings
  . 
  Assume further that there is a proper class $\calK$ of $C^{(n)}$-supercompact
  cardinals for sufficiently large $n$, and 
  $\MP^+(\calP,\calH(\kappa_\refl))_{\Gamma}$ holds for the set of formulas $\Gamma$ as defined
  in \Thmof{p-hierarchies-6-2}. Then, there is a class \po\/ $\clpoP$ 
  \st\ we have
  \begin{xitemize}
  \xitemx[] $\forces{\clpoP}{\neg\,(\calP,\emptyset)_{\Pi_2}\mbox{-\RcA, }
    \neg\,(\calP,\emptyset)_{\Sigma_2}\mbox{-\RcA\ and\/ }
    \MP^*(\calP,\calH(\kappa_\refl))_{\Gamma}}$.
  \end{xitemize}
\end{Cor}
\prf Note that
$\MP^+(\calP,\calH(\kappa_\refl))_{\Sigma_1}$ for $\calP$ as here implies $\neg\CH$ (see 
Fuchino and Usuba \cite{recurrence}, Theorem 3.3). 

\assertof{1}: By \Propof{p-rec-GA-1-1-0},\,\assertof{1} and 
\Propof{p-hierarchies-6-2},\,\assertof{1}.\smallskip

\assertof{2}: By \Propof{p-rec-GA-1-1-0},\,\assertof{1} and 
\Thmof{p-hierarchies-6-2},\,\assertof{2}.  
\qedofCor\qedskip

Typical instances of $\calP$ in \Corof{p-hierarchies-7} are when $\calP$ is the class of 
all proper \pos, the  class of all semi-proper \pos, or the class of all 
stationary preserving \pos.  

\begin{Problem}
  Do some theorems hold which would imply certain non-implications similar to those in 
  \Corof{p-hierarchies-7} for $\calP=$ $\sigma$-closed \pos, or $\calP=$ ccc \pos? 
\end{Problem}

\section{Generic absoluteness under Recurrence Axioms}
\Label{genabs-rec}
The conclusion of the following \Thmof{p-genabs-rec-0} generalizes that of Viale's \Thmof{p-theintro-a} 
(Theorem 1.4 in 
\cite{viale-revisited}). Note that %% in contrast to Viale's Theorem, 
the  
assumption in our \Thmof{p-genabs-rec-0}, the Recurrence Axiom $(\calP,\calH(\kappa))$-\RcAp\ 
  for an uncountable cardinal $\kappa$ and an iterable class $\calP$ of \pos, namely, 
is of much lower consistency strength than the 
assumptions in Viale's Theorem for some instances of $\calP$.
\memo{Schindler: Semi-proper forcing, remarkable cardinals, and Bounded Martin’s Maximum}
Actually the assumption of \Thmof{p-genabs-rec-0} is even compatible with $V=L$ 
for many ``natural'' classes $\calP$ of \pos\ including the cases ``$\calP=$ all ccc \pos'' 
or ``$\calP=$ all proper \pos'' (see Theorems 5.6, 
5.10 in \cite{hamkins}). 
%% Recurrence Axioms have relatively low consistency strength. See Theorems 5.6, 5.10 
%% in \cite{hamkins} in case of ``$\calP=$ all ccc \pos'' or ``$\calP=$ all proper \pos'' 
%% or Theorem 1.6 in \cite{ikegami-trang} 
For ``$\calP=$ all stationary preserving 
\pos'',  Theorem 1.6 in Ikegami and Trang \cite{ikegami-trang} proves that 
the existence of proper class many strongly compact cardinals plus a reflecting cardinal 
is an upper bound of the 
consistency strength of the Maximality Principle for the $\calP$.

Known lower bound of this Recurrence Axiom is also large. By Ikegami-Trang 
\Thmof{p-theintro-1}, $(\mbox{stationary preserving},\calH(\continuum))$-\RcA\ is 
equivalent with \BMM. Schindler \cite{schindler2} shows that \BMM\ implies that there is an 
inner model with a strong cardinal.

%% Even though this upper bound seems to be stronger than the assumption of 
%% Viale's \Thmof{p-theintro-a}, 
%% it is still a mild assumption compared with the consistency strength of the 
%% assumption 
%% of 
In contrast, the Maximality Principle for 
$\calP={}$ semi-proper \pos, the consistency strength is much lower than this 
by Asperó \cite{aspero}{\ifprivate\privatecolor\ (see also \subsectionof{Yah} below)\fi}. Note that, in general, 
semi-proper and stationary preserving are not identical notions. 

The existence of the tightly 
super-$C^{(\infty)}$ $\calP$-Laver-gen.\ hyperhuge cardinal $\kappa$ is the known 
Laver-generic large cardinal axiom which implies the full Recurrence Axiom for $\calP$ and
$\calH(\kappa)$ (see Fuchino and Usuba \cite{recurrence}).  

There is 
practically no (consistent) generic large cardinal axiom formalizable in a single formula which 
also implies
$(\calP,\emptyset)$-\RcA\ for any sufficiently general class $\calP$ of \pos\ 
(\cite{future}{\ifprivate\privatecolor, see also \sectionof{class-many} below\fi}). 

\memox{define $\kappa_\refl$.}
In \cite{recurrence}, it is proved that for an iterable class $\calP$ of \pos, 
the existence of the tightly 
super-$C^{(\infty)}$ $\calP$-Laver-gen.\ hyperhuge cardinal $\kappa$ 
with $\kappa=\kappa_\refl$ implies
$(\calP,\calH(\kappa_\refl))$-\RcAp\ where $\kappa_\refl$ is defined as
$\kappa_\refl:=\max\ssetof{\continuum,\aleph_2}$. 
Note that this does not 
contradict what we mentioned in the last paragraph since  the tightly
super-$C^{(\infty)}$ $\calP$-Laver-gen.\ hyperhugeness of $\kappa_\refl$ is only expressed 
by an axiom schema. Note that by \cite{recurrence} we know the exact consistency strength 
of this principle (as that of $\kappa_\refl$ being super-$C^{(\infty)}$ hyperhuge in the 
bedrock). 

In \sectionof{genabs-Laver} we show that the generalization of the conclusion of Viale's 
\Thmof{p-theintro-a} 
(like that of the following \Thmof{p-genabs-rec-0}) already 
follows from tight $\calP$-Laver-gen.\ hugeness. This assumption is still much stronger 
than that of Viale's \Thmof{p-theintro-a} but the upper bound of the consistency strength of 
this Laver-genericity is far below the consistency strength of a tight 
super-$C^{(\infty)}$ $\calP$-Laver-gen.\ hyperhuge cardinal. 
\ifextended\smallskip\fi

\begin{Thm}\Label{p-genabs-rec-0}
  Suppose that $\calP$ is an iterable $\Sigma_n$-definable class of \pos\ for $n\geq 2$ and 
  $(\calP,\calH(\kappa))_{\Sigma_n\cup\Gamma}$-\RcAp\ 
  holds for an uncountable cardinal $\kappa$ where $\Gamma$ is the set of all formulas which are 
  conjunction of a $\Sigma_2$-formula and a $\Pi_2$-formula. Then, 
  for any $\poP\in\calP$ \st\ 
  $\forces{\poP}{\BFA_{\LT\kappa}(\calP)}$, 
  \begin{xitemize}
  \item[] 
    $\calH(\mu^+)^\uniV\prec_{\Sigma_2}\calH(\mu^+)^{\uniV[\genG]}$\quad 
    holds 
    for all $\mu<\kappa$ and for $(\uniV,\poP)$-generic $\genG$. 
  \end{xitemize}
  %% Thus $\calH(\kappa)^\uniV\prec_{\Sigma_2}\calH(\kappa)^{\uniV[\genG]}$ holds for all
  %% $(\uniV,\poP)$-generic filter $\genG$. 
  Thus, we have
  $\calH(\kappa)^\uniV\prec_{\Sigma_2}\calH((\kappa^{(+)})^{\uniV[\genG]})^{\uniV[\genG]}$. 
\end{Thm}\memox{For $n\geq2$ if $\calP$ is $\Sigma_n$-definable then isn't 
  $(\calP,\calH(\kappa))_{\Sigma_n}$-\RcAp\ enough to prove this?}
\prf Suppose that $\poP\in\calP$ is \st\
$\forces{\poP}{\BFA_{\LT\kappa}(\calP)}$ and $\genG$ is a 
$(\uniV,\poP)$-generic filter.  
Let 
$\varphi=\varphi(x)$ be a $\Sigma_2$-formula in $\Lin$, and
$\varphi(x)=\exists y\,\psi(x,y)$ for a $\Pi_1$-formula $\psi$ in $\Lin$. Let $\mu<\kappa$ 
and $a\in\calH(\mu^+)$ ($\subseteq\calH(\kappa)$). We have to show that
$\calH(\mu^+)^\uniV\models\varphi(a)$ $\Leftrightarrow$
$\calH((\mu^+)^{\uniV[\genG]})^{\uniV[\genG]}\models\varphi(a)$. 

Suppose first that $\calH(\mu^+)^\uniV\models\varphi(a)$. Let $b\in\calH(\mu^+)^\uniV$ be 
\st\ $\calH((\mu^+)^{\uniV})^{\uniV}\models\psi(a,b)$. Since we have
$\uniV\models\BFA_{\LT\kappa}(\calP)$ by Ikegami-Trang \Thmof{p-theintro-1}, it follows 
that $\calH((\mu^+)^{\uniV[\genG]})^{\uniV[\genG]}\models\psi(a,b)$ by Bagaria's 
Absoluteness \Thmof{p-theintro-0}, and thus
$\calH((\mu^+)^{\uniV[\genG]})^{\uniV[\genG]}\models\varphi(a)$. 

Note that we did not use the assumption ``$\forces{\poP}{\BFA_{\LT\kappa}(\calP)}$'' for 
this direction. 

Suppose now $\calH((\mu^+)^{\uniV[\genG]})^{\uniV[\genG]}\models\varphi(a)$. By
$(\calP,\calH(\kappa))_{\Sigma_n\cup\Gamma}$-\RcAp, there is a $\calP$-ground $\uniW$ of $\uniV$ \st\
\begin{xitemize}
\xitem[x-genabs-rec-a-0] 
  $\uniW\modelof{\BFA_{\LT\mu^+}(\calP)\ \land\ \calH(\mu^+)\models\varphi(a)}$. 
\end{xitemize}
Note that the formula in \xitemof{x-genabs-rec-a-0} is $\Sigma_n$ if $n\geq 3$ and $\Gamma$ 
if $n=2$. \memo{add some more details the extended version \ see\\
  \mbox{generic-absoluteness-revisited-xx - annotated 30.pdf} p.29\\
  \mbox{Scan\_2024-05-14--12.42 gappo - annotated} pp.58-59}

Let $b\in\calH((\mu^+)^\uniW)^\uniW$ be \st\ $\uniW\modelof{\calH(\mu^+)\models\psi(a,b)}$.
By Bagaria's Absoluteness \Thmof{p-theintro-0}, and since $\uniV$ is a $\calP$-generic 
extension of $\uniW$, it follows that $\uniV\modelof{\calH(\mu^+)\models\psi(a,b)}$ and 
hence $\calH(\mu^+)^\uniV\models\varphi(a)$. 

For the last statement of the present theorem, let $\varphi$ be a $\Sigma_2$-formula, and
$a\in\calH(\kappa)$. If $\calH(\kappa)\models\varphi(a)$, then, by 
\xitemof{x-theintro-0}, there is $\mu<\kappa$ 
\st\ $\calH(\mu^+)\models\varphi(a)$. By the first part of the theorem, it follows that
$\calH((\mu^+)^{\uniV[\genG]})^{\uniV[\genG]}\models\varphi(a)$. Thus
$\calH((\kappa^{(+)})^{\uniV[\genG]})^{\uniV[\genG]}\models\varphi(a)$ by \xitemof{x-theintro-0}. 

If $\calH((\kappa^{(+)})^{\uniV[\genG]})^{\uniV[\genG]}\models\varphi(a)$, then there is $\mu<\kappa$ 
\st\ $\calH((\mu^+)^{\uniV[\genG]})^{\uniV[\genG]}\models\varphi(a)$ (this is also shown using 
\xitemof{x-theintro-0}). Hence $\calH((\mu^+)^{\uniV})\models\varphi(a)$ by the first 
part of the theorem. \qedofThm\qedskip

\memox{$\INS$ version of Bagaria's Theorem $\INS$ extension of \Thmabove.}
Note that by \Lemmaof{p-theintro-4-0}, the conclusion \xitemof{x-genabs-rec-a-0} of 
\Thmof{p-genabs-rec-0} can be yet strengthened to   
\begin{xitemize}
\xitem[x-genabs-rec-a-1] 
  $\pairof{\calH(\mu^+),\in, I_\NS}^\uniV
  \prec_{\Sigma_2}\pairof{\calH(\mu^+), \in, I_\NS}^{\uniV[\genG]}$\quad 
  holds for all $\mu<\kappa$ and for $(\uniV,\poP)$-generic $\genG$. 
\end{xitemize}

\section{Generic absoluteness under Laver-genericity}
\Label{genabs-Laver} Laver-genericity also implies a conclusion similar to that of Viale's 
\Thmof{p-theintro-a}.  
Although this fact does not have an advantage in terms of consistency strength in 
comparison with \Thmof{p-genabs-rec-0}, the 
Laver-generic large cardinal we need in \Thmof{p-genabs-Laver-1} below is much weaker than 
the tight super $C^{(\infty)}$ $\calP$-Laver-generic hyperhugeness, 
the generic large cardinal property which is known to imply the corresponding Recurrence 
Axiom used in \Thmof{p-genabs-rec-0} (see \cite{future} and \cite{recurrence}). \memo{more  
  precise localization of results in \cite{future} and \cite{recurrence} needed!}

In Viale \cite{viale-revisited}, the absoluteness statement of his \Thmof{p-theintro-a} is 
also discussed in connection with the Resurrection Axiom (see \Thmof{p-genabs-Laver-a-1}). 

Adopting the generalized setting introduced in Fuchino \cite{future}, we define the  
Resurrection Axiom as follows: 
for an iterable class $\calP$ of \pos\ and a definition $\mu^\bullet$ of an uncountable cardinal (e.g.\ as 
  $\aleph_1$, $\aleph_2$, $\continuum$, $(\continuum)^+$, $\kappa_\refl$ 
etc.),\footnote{More precisely, when we say ``$\mu^\bullet$ is a definition of an 
  uncountable cardinal'' we mean that \ZF\ or \ZFC\ proves the statement ``$\mu^\bullet$ 
  uniquely exists and $\mu^\bullet$ is an uncountable cardinal''. }
  {\It the Resurrection Axiom} for $\calP$ and $\mu^\bullet$ is the statement: 
\begin{xitemize}
\item[{\darkred$\RA(\calP, \mu^\bullet)$:}\quad] For any $\poP\in\calP$, there is 
  a $\poP$-name $\utpoQ$  
  of a \po\ \st\ $\forces{\poP}{\utpoQ\in\calP}$ and 
  $\calH(\mu^\bullet)^\uniV\prec\calH(\mu^\bullet)^{\uniV[\genH]}$ for any
  $(\uniV,\poP\ast\utpoQ)$-generic $\genH$.\footnote{Here we mean with 
    ``$\calH(\mu^\bullet)^\uniV\prec\calH(\mu^\bullet)^{\uniV[\genH]}$'' the elementarity $\calH((\mu^\bullet)^\uniV)^\uniV\prec\calH((\mu^\bullet)^{\uniV[\genH]})^{\uniV[\genH]}$.}
\end{xitemize}

\begin{Lemma}{\rm (Hamkins, and Johnstone \cite{hamkins-johnstone})}\Label{p-genabs-Laver-a-0} 
  Assume that $\calP$ 
  is an iterable class of \pos, $\mu^\bullet$ is a definition of an uncountable cardinal, 
  and $\RA(\calP,\mu^\bullet)$ holds. Then \wassertof{1} $\BFA_{\LT\mu^\bullet}(\calP)$ 
  holds.\smallskip

  \wassert{2} If all elements of\/ $\calP$ preserve stationarity of subsets of $\omega_1$, 
  $2^{\aleph_0}=2^{\aleph_1}$, and $\mu^\bullet={}$``\/$\continuum$'', then 
  $\BFA^{+\LT\mu^\bullet}_{\LT\mu^\bullet}(\calP)$ holds.
\end{Lemma}
\prf \assertof{1}: It is easy to check that \assertof{c} of Bagaria's \Thmof{p-theintro-0} 
holds.\footnote{Actually we do not need this condition in 
    this part of the proof and hence we can obtain the Bounded Forcing Axiom under a weaker 
    notion of Resurrection Axiom in which the second step $\utpoQ$ may be anything.}

 {\ifextended\extendedcolor Suppose that $\overline{a}\in\calH(\mu^\bullet)$ and 
  $\varphi$ is a $\Sigma_1$-formula in $\Lin$. For $\poP\in\calP$, let $\genG$ be 
  a $(V,\poP)$-generic filter.

  If $\calH(\mu^\bullet)^\uniV\models\varphi(\overline{a})$, then we have 
  $\uniV\models\varphi(\overline{a})$ and hence $\uniV[\genG]\models\varphi(\overline{a})$. Thus 
  we have $\calH(\mu^\bullet)^{\uniV[\genG]}\models\varphi(\overline{a})$ by 
  \xitemof{x-theintro-0}.

  Suppose now that $\calH(\mu^\bullet)^{\uniV[\genG]}\models\varphi(\overline{a})$. Let 
  $\utpoQ$ be a $\poP$-name of a \po\ 
  \st\ $\forces{\poP}{\utpoQ\in\calP}$
  \st\ for a $(\uniV[\genG], \utpoQ^\genG)$-generic $\genH$,
  \begin{xitemize}
  \xitemA[x-genabs-Laver-a-a] 
    $\calH(\mu^\bullet)^\uniV\prec\calH(\mu^\bullet)^{\uniV[\genG\ast\genH]}$.
  \end{xitemize}

  Since $\uniV[\genG]\models\varphi(\ol{a})$, we have 
  $\calH(\mu^\bullet)^{\uniV[\genG\ast\genH]}\models\varphi(\overline{a})$
  by the same argument as in the first part of this proof. Thus, by 
  \xitemof{x-genabs-Laver-a-a}, it follows that
  $\calH(\mu^\bullet)^\uniV\models\varphi(\overline{a})$. 
\fi}
\smallskip

\assertof{2}: Similarly to \assertof{1} using \assertof{c} of 
\Thmof{p-theintro-5}.\footnote{Note that for this proof, the weaker variant of $\RA(\calP)$ 
  as in the proof of \assertof{1} is apparently not sufficient. }
\qedofLemma

\begin{Thm}{\rm (A generalization of Theorem 5.1 in 
    Viale \cite{viale-revisited})}\Label{p-genabs-Laver-a-1} 
  Suppose that $\calP$ is an iterable class of \pos, $\mu^\bullet$ is a definition of an 
  uncountable cardinal and $\RA(\calP,\mu^\bullet)$ holds. Then we have
  \begin{xitemize}
  \item[] $\calH(\mu^\bullet)^\uniV\prec_{\Sigma_2}\calH(\mu^\bullet)^{\uniV[\genG]}$
  \end{xitemize}
  for any $\poP\in\calP$ \st\ $\forces{\poP}{\BFA_{\LT\mu^\bullet}(\calP)}$, and
  $(\uniV,\calP)$-generic $\genG$.
\end{Thm}
\prf An argument similar to that of \Lemmaof{p-genabs-Laver-a-0} will 
do.\memo{\mbox{Scan\_2024-05-14--12.42\ gappo\ -\ annotated,} pp.25-26}
\qedofThm\qedskip

In \cite{future}, the boldface version of the following is proved: 

\begin{Thm}{\rm (Fuchino \cite{future})}\Label{p-genabs-Laver-a-2}
  For an iterable class $\calP$ of \pos, and a definition $\mu^\bullet$ of a cardinal, 
  if $\mu^\bullet$ is tightly $\calP$-Laver gen.\ 
  superhuge, then $\RA(\calP,\mu^\bullet)$ holds.\qed
\end{Thm}\memo{\cite{future} Theorem 7.1}
\begin{Cor}\Label{p-genabs-Laver-a-3}
  For an iterable class $\calP$ of \pos, and a definition $\mu^\bullet$ of a cardinal, if
  $\mu^\bullet$ is tightly $\calP$-Laver gen.\  
  superhuge, then we have
  \begin{xitemize}
  \item[] $\calH(\mu^\bullet)^\uniV\prec_{\Sigma_2}\calH(\mu^\bullet)^{\uniV[\genG]}$
  \end{xitemize}
  for any $\poP\in\calP$ \st\ $\forces{\poP}{\BFA_{\LT\mu^\bullet}(\calP)}$ and
  $(\uniV,\poP)$-generic $\genG$.
\end{Cor}
\prf By \Thmof{p-genabs-Laver-a-1} and \Thmof{p-genabs-Laver-a-2}. \qedofCor\qedskip

Note that for many cases (with natural setting of $\calP$ and the notion of large 
cardinal), if $\kappa$ is $\calP$-Laver-gen.\ large cardinal, 
then $\kappa=\kappa_\refl$ (see Fuchino, Ottenbreit Maschio Rodrigues, and Sakai \cite{sfetal-II}).

In the following, we show in a direct proof, that \Corof{p-genabs-Laver-a-3} can be yet 
slightly improved (see \Thmof{p-genabs-Laver-1} below).

\memox{Define the 
  notion ``iterable''}

It is known that 
Laver-gen.\ large cardinal axiom proves strong forms of  
double-plus forcing axioms (see Theorem 5.7 in \cite{sfetal-II}). In the 
\Propof{p-genabs-Laver-a} below we  
only recap a part of this result needed for the following argument.

\memox{iterable class of pos\\Laver-genericity (supercompact, huge, tightly)}

For a class $\calP$ of \pos, and cardinals $\kappa$ and $\lambda$,
\begin{xitemize}
\item[{\rm\darkred($\FA^{+\LT\lambda}_{\LT\kappa}(\calP)$): }] For any $\poP\in\calP$, any family
  $\calD$ 
  of dense subsets of $\calP$ with $\cardof{\calD}<\kappa$, and any family $\calS$ 
  of $\poP$-names of stationary subsets of $\omega_1$ with $\cardof{\calS}<\lambda$, there 
  is a $\calD$-generic filter $\genG$ on $\poP$ \st\ $\utilde{S}[\genG]$ is a stationary 
  subset of $\omega_1$ for all $\utilde{S}\in\calS$. 
\end{xitemize}
Note that \MMpp\ is just
$\FA^{+\LT\aleph_2}_{\LT\aleph_2}(\mbox{stationary preserving \pos})$.  
$\FA_{\LT\kappa}(\calP)$ is the principle we obtain by dropping the mention about $\calS$ 
from the definition of $\FA^{+\LT\lambda}_{\LT\kappa}(\calP)$.

\begin{Prop}\Label{p-genabs-Laver-a} \wassertof{1} Suppose that $\kappa$ is
  $\calP$-Laver-gen.\ 
  supercompact. Then\\
  $\FA_{\LT\kappa}(\calP)$ holds.\smallskip

  \wassert{2} If all elements of the class $\calP$ of \pos\ are stationary preserving and
  $\kappa$ is $\calP$-Laver-gen.\ 
  supercompact, then $\FA^{+\LT\kappa}_{\LT\kappa}(\calP)$ holds.
\end{Prop}
\prf We prove \assertof{2}. \assertof{1} can be shown by a subset of this proof. 

Assume that $\kappa$ is a $\calP$-Laver-gen.\ supercompact cardinal, and 
let $\calP$, $\calD$, $\calS$ be as in the definition of
$\FA^{+\LT\kappa}_{\LT\kappa}(\calP)$. Let $D_\alpha$, $\alpha<\mu$ 
and $\utilde{S}_\alpha$, $\alpha<\mu'$ be enumerations of $\calD$ and $\calS$ respectively. 

Let $\lambda=\cardof{\poP}$. \Wolog, we may assume that $\lambda>\kappa$ and the underlying 
set of $\poP$ is $\lambda$. Let $\utpoQ$ be a $\poP$-name with
$\forces{\poP}{\utpoQ\in\calP}$ and \st\ for any $(\uniV,\poP\ast\utpoQ)$-generic $\genH$, 
there are $j$, $M\subseteq\uniV[\genG]$ \st\ $\Elembed{j}{\uniV}{M}{\kappa}$,
$j(\kappa)>\lambda$, $j\imageof{\lambda}$, $\poP$, $\poP\ast\utpoQ$, $\genH\in M$.

Let $\genG$ be the $\poP$ part of $\genH$. Then, since $j(\mu)=\mu$, $j(\mu')=\mu'$,
$j(\calD)=\setof{j(D_\alpha)}{\alpha<\mu}$,  and
$j(\calS)=\setof{j(\utilde{S}_\alpha)}{\alpha<\mu}$, we have
\begin{xitemize}
\xitem[x-genabs-rec-0] $M\modelof{\,
  \begin{array}[t]{@{}l}
    j\imageof{\genG}\mbox{ generates a filter on }j(\poP)\mbox{ which is}\\
    j(\calD)\mbox{-generic, and realizes elements of }j(\calS)\mbox{ to be stationary}}.
  \end{array}$
\end{xitemize}
Note that we need here the condition that $\calP$ is stationary preserving since otherwise 
the stationary set 
$\utilde{S}[\genG]$ in $\uniV[\genG]$ might be no more stationary in $\uniV[\genH]$. 

\xitemof{x-genabs-rec-0} implies that 
\begin{xitemize}
\item[] $M\modelof{\,
  \begin{array}[t]{@{}l}
    \mbox{there is a }j(\calD)\mbox{-generic filter on }j(\poP)\\
    \mbox{which realizes all elements of }j(\calS)\mbox{ to be stationary}}.
  \end{array}$
\end{xitemize}

By elementarity, it follows that 
\begin{xitemize}
\item[] $\uniV\modelof{\,
  \begin{array}[t]{@{}l}
    \mbox{there is a }\calD\mbox{-generic filter on }\poP\\
    \mbox{which realizes all elements of }\calS\mbox{ to be stationary}}.
  \end{array}$
\end{xitemize}
\qedofProp
\qedskip

\begin{Lemma}\Label{p-genabs-Laver-0}
  Suppose that $\kappa$ is $\calP$-Laver-gen.\ supercompact for an iterable $\calP$.
  Then we have
  $\calH(\kappa)^\uniV\prec_{\Sigma_1}\calH((\kappa^{(+)})^{\uniV[\genG]})^{\uniV[\genG]}$ 
  for any  
  $\poP\in\calP$ and $(\uniV,\poP)$-generic $\genG$. 
\end{Lemma}
\prf By \Propof{p-genabs-Laver-a},\,\assertof{1} and Bagaria's Absoluteness \Thmof{p-theintro-0}. 
\qedofLemma\qedskip

The following theorem improves \Corof{p-genabs-Laver-a-3}.

\begin{Thm}\Label{p-genabs-Laver-1} For an iterable class $\calP$ of \pos, suppose that 
  $\BFA_{\LT\kappa}(\calP)$ holds, and $\kappa$ is 
  tightly $\calP$-Laver-gen.\  
  huge.\footnotemark\ Then, 
  for any $\poP\in\calP$ \st\  $\forces{\poP}{\BFA_{\LT\kappa}(\calP)}$, 
  \begin{xitemize}
  \item[] 
    $\calH(\mu^+)^\uniV\prec_{\Sigma_2}\calH(\mu^+)^{\uniV[\genG]}$\quad
    holds 
    for all $\mu<\kappa$ and for $(\uniV,\poP)$-generic $\genG$. 
  \end{xitemize}
  Thus, we have
  $\calH(\kappa)^\uniV\prec_{\Sigma_2}\calH((\kappa^{(+)})^{\uniV[\genG]})^{\uniV[\genG]}$. 
\end{Thm}
\footnotetext{\Label{fn-genabs-Laver-0}Note that by \Propof{p-genabs-Laver-a},
  $\BFA_{\LT\kappa}(\calP)$ follows from the assumption that $\kappa$ 
  is $\calP$-Laver-generic supercompact. Thus the conclusion of the theorem follows from 
  the combination of the assumption  $\kappa$ 
  being $\calP$-Laver-generic supercompact and tightly $\calP$-Laver-gen.\ huge. Note also 
  that this combination follows from the tight $\calP$-Laver-gen.\ superhugeness of $\kappa$. }
\prf Suppose that $\forces{\poP}{\calH(\mu^+)\models\varphi(\ol{a})}$ for $\poP\in\calP$ 
with $\forces{\poP}{\BFA_{\LT\kappa}(\calP)}$, $\mu<\kappa$, $\Sigma_2$-formula $\varphi$ and
for $\ol{a}\in\calH(\mu^+)$. Let $\genG$ be a $(\uniV,\poP)$-generic filter. Then we have
\begin{xitemize}
\xitem[x-genabs-Laver-a-0] 
  $\uniV[\genG]\modelof{\BFA_{\LT\kappa}(\calP)\land\calH(\mu^+)\models\varphi(\ol{a})}$. 
\end{xitemize}
Let $\varphi=\exists y\psi(\ol{x},y)$ where $\psi$ is a $\Pi_1$-formula in $\Lin$. Let
$b\in\calH((\mu^+)^{\uniV[\genG]})^{\uniV[\genG]}$ be \st\
$\calH((\mu^+)^{\uniV[\genG]})^{\uniV[\genG]}\models\psi(\ol{a},b)$. 

Let $\utpoQ$ be a $\poP$-name with $\forces{\poP}{\utpoQ\in\calP}$ \st, for
$(\uniV,\poP\ast\utpoQ)$-generic $\genH$ with 
\begin{xitemize}
\xitem[x-genabs-Laver-a-1] 
  $\genG\subseteq\genH$ (under the 
  identification $\poP\circleq\poP\ast\utpoQ$), 
\end{xitemize}
there are $j$, $M\subseteq\uniV[\genH]$ \st\
$\Elembed{j}{\uniV}{M}{\kappa}$, 
\begin{xitemize}
\xitem[x-genabs-Laver-0] $\cardof{\poP\ast\utpoQ}\leq j(\kappa)$\qquad (by tightness), 
\xitem[x-genabs-Laver-1] $\poP$, $\poP\ast\utpoQ$, $\genH\in M$ and
\xitem[x-genabs-Laver-2] $j\imageof{j(\kappa)}\in M$. 
\end{xitemize}
\memo{[25.05.19(Mo09:40(JST))] \xitemof{x-genabs-Laver-2} can be replaced by ${V_{j(\lambda)}}^{\uniV[\genH]}\in M$. 
  Thus the assumption of the theorem can be replaced by $\calP$-LgLCA for extendible.}

By \xitemof{x-genabs-Laver-a-0}, \xitemof{x-genabs-Laver-a-1} and Bagaria's Absoluteness 
\Thmof{p-theintro-0} (applied to $V[\genG]$), we have $\uniV[\genH]\modelof{\psi(\ol{a},b)}$  and hence
$\uniV[\genH]\modelof{\calH(\mu^+)\models\psi(\ol{a},b)}$.

By \xitemof{x-genabs-Laver-0}, and \xitemof{x-genabs-Laver-2}, 
there is a $\poP$-name of $b$ in $M$. By \xitemof{x-genabs-Laver-1}, it follows that
$b\in M$. By similar argument, we have $\calH((\mu^+)^{\uniV[\genH]})^{\uniV[\genH]}\subseteq M$ and 
hence $\calH((\mu^+)^{\uniV[\genH]})^{\uniV[\genH]}=\calH((\mu^+)^M)^M\in M$. Thus we have
$M\modelof{\calH(\mu^+)\models\psi(\ol{a},b)}$.

By elementarity, it follows that
$\uniV\modelof{(\exists\ul{b}\in\calH(\mu^+))\ \calH(\mu^+)\models\psi(\ol{a},\ul{b})}$, and  
hence $\uniV\modelof{\calH(\mu^+)\models\varphi(\ol{a})}$ as desired. 
\smallskip

Suppose now that $\poP$, $\mu$, $\varphi$,
$\ol{a}$ are as above and assume that 
$\uniV\modelof{\calH(\mu^+)\models\varphi(\ol{a})}$ holds. For a $\Pi_1$-formula $\psi$ as 
above, let $b\in\calH(\mu^+)^\uniV$ be \st\ $V\modelof{\calH(\mu^+)\models\psi(\ol{a},b)}$.

Since $\uniV\models\BFA_{\LT\kappa}(\calP)$ by assumption, 
it follows that $\uniV[\genG]\models\psi(\ol{a},b)$ by Bagaria's Absoluteness 
\Thmof{p-theintro-0}, and hence $\uniV[\genG]\modelof{\calH(\mu^+)\models\varphi(\ol{a})}$. 

The last assertion of the theorem follows by the same argument as that given at the end of the 
proof of \Thmof{p-genabs-rec-0}. 
\qedofThm

\memo{Corollary Laver-gen.\ super huge ... $\calH(\kappa)\models$ ... }

\section{Some more remarks and open questions}\Label{misc}
In this final section we collect some observations we could not put in the appropriate 
places in previous sections, and mention some open questions. 

Since most of the claims in this section are easy consequences of known results, 
some of them may be already known. 

\subsection{Restricted Recurrence Axioms under Laver-genericity}\Label{RRAL}
As we already mentioned at the end of \sectionof{rec-GA}, the existence of a tightly $\calP$-Laver gen.\ ultrahuge 
cardinal $\kappa$ implies $(\calP,\calH(\kappa))_{\Sigma_2}$-\RcAp\ (Theorem 21 in \cite{janos}). 
This result can be slightly improved so that its conclusion stands in line with the assumptions of 
\Thmof{p-genabs-rec-0} for the case of $n=2$ %% (see \Corof{p-Lg-RcA-1})
. 

\begin{Thm}{\rm (A slightly improved version of Theorem 21 in Fuchino \cite{janos})}
  \Label{p-Lg-RcA-0} Suppose that $\kappa$ is tightly $\calP$-Laver-generically ultrahuge  
  for an 
  iterable class $\calP$ of \pos%% \ and there are stationarily many 
  %% inaccessible cardinals 
  . Then
  $(\calP,\calH(\kappa))_{\Gamma}$-\RcAp\ holds where $\Gamma$ is the set of all formulas 
  which are conjunctions of a $\Sigma_2$-formula and a $\Pi_2$-formula.
\end{Thm}
\prf A slight modification of the proof given in \cite{janos} will do. Nevertheless, we 
present the proof here because of the subtlety of the modification of the proof around 
\xitemof{x-Lg-RcA-2-0} below. 
%%{\ifextended\extendedcolor

Assume that $\kappa$ is tightly $\calP$-Laver generically ultrahuge for an 
iterable class $\calP$ of \pos. 

Suppose that $\varphi=\varphi(\overline{x})$ is $\Sigma_2$ formula (in $\Lin$), 
$\psi=\psi(\overline{x})$ is $\Pi_2$ formula (in $\Lin$), 
$\overline{a}\in\calH(\kappa)$, and  
$\poP\in\calP$ is \st\ 
\begin{xitemize}
\xitem[x-Lg-RcA-a] 
  $\uniV\models\forces{\poP}{\varphi(\overline{a})\land \psi(\overline{a})}$. 
\end{xitemize}

Let $\lambda>\kappa$ be %% an inaccessible 
\st\ $\poP\in\uniV_\lambda$ and 
\begin{xitemize}
\xitem[x-Lg-RcA-0] $V_\lambda\prec^{}_{\Sigma_n}\uniV\mbox{ for a sufficiently large }n$. 
\end{xitemize}
In particular, we may assume that we have chosen the $n$ above so that a sufficiently large 
fragment of \ZFC\ holds in 
$V_\lambda$ (``sufficiently large'' means here, in particular, in terms of 
\Lemmaof{p-hierarchies-5-0},\,\assertof{1} and that the argument at the end of this proof is possible).  

Let $\utpoQ$ be a $\poP$-name \st\ $\forces{\poP}{\utpoQ\in\calP}$, and for
$(\uniV,\poP\ast\utpoQ)$-generic $\genH$, there are $j$, $M\subseteq\uniV[\genH]$ with 
\begin{xitemize}
  \xitem[x-Lg-RcA-1] 
  $\Elembed{j}{\uniV}{M}{\kappa}$, 
\end{xitemize}
\begin{xitemize}
  \xitem[x-Lg-RcA-1-0] 
  $j(\kappa)>\lambda$, 
\end{xitemize}
\begin{xitemize}
  \xitem[x-Lg-RcA-1-1] 
  $\poP\ast\utpoQ,\ \poP,\ \genH,\ {V_{j(\lambda)}}^{\uniV[\genH]}\in M,\mbox{ and }$
\end{xitemize}
\begin{xitemize}
  \xitem[x-Lg-RcA-1-2]
  $\cardof{\poP\ast\utpoQ}\leq j(\kappa)$. 
\end{xitemize}
By \xitemof{x-Lg-RcA-1-2}, we may assume that the underlying set 
of $\poP\ast\utpoQ$ is $j(\kappa)$ and $\poP\ast\utpoQ\in {V_{j(\lambda)}}^\uniV$. 

Let $\genG:=\genH\cap\poP$. 
Note that $\genG\in M$ by \xitemof{x-Lg-RcA-1-1}. 
We have 
\begin{xitemize}
\iffalse
0.8*1.8=1.4400000000000002
\fi
\xitem[x-Lg-RcA-2] 
  ${V_{j(\lambda)}}^M\ubecause{=}{}{by \xitemof{x-Lg-RcA-1-1}}
  {V_{j(\lambda)}}^{\uniV[\genH]}
  \obecause{=}{1.44ex}{\hspace{20em}
    \vbox{\hbox{Since ${V_{j(\lambda)}}^M$ ($={V_{j(\lambda)}^{\uniV[\genH]}}$) satisfies a 
    sufficiently large fragment of 
    \ZFC\vspace{-0.18ex}}\hbox{by elementarity of $j$, and by 
        \Lemmaof{p-hierarchies-5-0},\,\assertof{1}}}} 
  {V_{j(\lambda)}}^\uniV[\genH].$ 
\end{xitemize}
Thus, by \xitemof{x-Lg-RcA-1-1} and by the definability of grounds, we have 
${V_{j(\lambda)}}^\uniV\in M$ and ${V_{j(\lambda)}}^\uniV[\genG]\in M$. 
We may assume that $V_{j(\lambda)}^{\uniV}$ as a ground of $V_{j(\lambda)}^M$ satisfies a 
large enough fragment of \ZFC. 

\begin{Claim}
  \Label{cl-Lg-RcA-0}
  ${V_{j(\lambda)}}^\uniV[\genG]\models\varphi(\overline{a})\land\psi(\overline{a})$.  
\end{Claim}
\noindent
\prfofClaim
%% By elementarity \xitemof{x-Lg-RcA-1} and \xitemof{x-Lg-RcA-1-1}, $j(\lambda)$ is 
%% inaccessible in $\uniV[\genH]$ and hence also in $\uniV[\genG]$ and $\uniV$.
%%
By \Lemmaof{p-hierarchies-5-0},\,\assertof{1}, ${V_\lambda}^{\uniV}[\genG]={V_\lambda}^{\uniV[\genG]}$, and
${V_{j(\lambda)}}^{\uniV}[\genG]={V_{j(\lambda)}}^{\uniV[\genG]}$. %% by 
%% \xitemof{x-Lg-RcA-2}.  
By \xitemof{x-Lg-RcA-0}, both ${V_\lambda}^\uniV[\genG]$ and $V_{j(\lambda)}^\uniV[\genG]$ 
satisfy still large enough fragment of \ZFC. Thus, by \Lemmaof{p-Lg-RcA-0-0} below, 
it follows that 
\begin{xitemize}
  \xitem[x-Lg-RcA-2-0] 
  ${V_\lambda}^\uniV[\genG]\prec_{\Sigma_1}{V_{j(\lambda)}}^\uniV[\genG]\prec_{\Sigma_1}V[\genG]$. 
\end{xitemize}
By \xitemof{x-Lg-RcA-a} and \xitemof{x-Lg-RcA-0}, we have
${V_\lambda}^\uniV[\genG]\models\varphi(\overline{a})$ and $\uniV[\genG]\models\psi(\ol{a})$. By 
\xitemof{x-Lg-RcA-2-0} and since 
$\varphi$ is $\Sigma_2$, and $\psi$ is $\Pi_2$, it follows that
${V_{j(\lambda)}}^\uniV[\genG]\models\varphi(\overline{a})\land\psi(\overline{a})$.  
\qedofClaim\qedskip

Thus we have 
\begin{xitemize}
  \xitem[x-Lg-RcA-3] 
  $M\modelof{\mbox{there is a }
  \calP\mbox{-ground }N\mbox{ of }V_{j(\lambda)}\mbox{ with }N\models\varphi(\overline{a})\land\psi(\overline{a})}$.
\end{xitemize}

By the elementarity \xitemof{x-Lg-RcA-1}, it follows that 
\begin{xitemize}
  \xitem[x-Lg-RcA-4] 
  $\uniV\modelof{\mbox{there is a }
  \calP\mbox{-ground }N\mbox{ of }V_{\lambda}\mbox{ with }N\models\varphi(\overline{a})\land\psi(\overline{a})}$.
\end{xitemize}

Now by \xitemof{x-Lg-RcA-0}, it follows that there is a $\calP$-ground $\uniW$ of $\uniV$ \st\\
$\uniW\models\varphi(\overline{a})\land\psi(\overline{a})$. %% \fi}
  \qedofThm\qedskip

%%{\ifextended\extendedcolor
We used the following variation of \xitemof{x-theintro-0}
in the proof of \Thmof{p-Lg-RcA-0} to obtain 
\xitemof{x-Lg-RcA-2-0}: 
\begin{Lemma}\Label{p-Lg-RcA-0-0}
  Suppose that $\delta$, $\delta'\in\On$, $\delta<\delta'$ and both $V_\delta$ and
  $V_{\delta'}$ satisfy a sufficiently large fragment of \ZFC. Then we have
  $V_\delta\prec_{\Sigma_1}V_{\delta'}\prec_{\Sigma_1}\uniV$.\ifextended\else\qed\fi
\end{Lemma}
\prf Suppose that $\ol{a}\in V_\delta$ and $\psi(\ol{x},\ol{y})$ is a bounded formula in
$\Lin$.

If $V_\delta\models\exists\ol{y}\psi(\ol{a},\ol{y})$, then there are $\ol{b}\in V_\delta$ 
\st\ $V_\delta\models\psi(\ol{a},\ol{b})$. It follows that
$V_{\delta'}\models\psi(\ol{a},\ol{b})$ and hence
$V_{\delta'}\models\exists\ol{y}\psi(\ol{a},\ol{y})$.

Suppose now that $V_{\delta'}\models\exists\ol{y}\psi(\ol{a},\ol{y})$. 
Since $V_{\delta'}$ satisfies a sufficiently large fragment of \ZFC, there is
$M\in V_{\delta'}$ \st\ 
\begin{xitemize}
\item[] 
  $V_{\delta'}\modelof{\,{}
  \begin{array}[t]{@{}l}
    \cardof{\trcl(\ol{a})}=\cardof{M},\,
    M\mbox{ is transitive},\,\ol{a}\in M,\,\mbox{there are }\ol{b}\in M\\
   \mbox{\st\ }M\models\psi(\ol{a},\ol{b})}.
  \end{array}$ 
\end{xitemize}

But then such $M$ as above must be an element of $V_\delta$ and thus 
\begin{xitemize}
\item[] 
  $V_{\delta}\modelof{\,{}
  \begin{array}[t]{@{}l}
    \cardof{\trcl(\ol{a})}=\cardof{M},\,
    M\mbox{ is transitive},\,\ol{a}\in M,\,\mbox{there are }\ol{b}\in M\\
   \mbox{\st\ }M\models\psi(\ol{a},\ol{b})}.
  \end{array}$ 
\end{xitemize}

It follows that $V_\delta\models\exists\ol{y}\,\psi(\ol{a},\ol{y})$. 

The argument above shows that $V_\delta\prec_{\Sigma_1}V_{\delta'}$. 
$V_{\delta'}\prec_{\Sigma_1}\uniV$ can be shown with practically the same argument. 
\qedofLemma%%\fi}

%% \begin{Cor}\Label{p-Lg-RcA-1}
%%   Suppose that $\calP$ is a 
%% \end{Cor}

\subsection{Separation of some other axioms and assertions}\Label{Sep}
In \sectionof{hierarchies}, we separated some instances of $(\calP,\calH(\kappa))_\Gamma$-\RcA\ 
and $\MP(\calP,\calH(\kappa))_\Gamma$ by compatibility with the Ground Axiom (\GA). The same idea 
can be also used to separate some other principles and axioms. 

\begin{Thm}\Label{p-Lg-RcA-1-0}
 ``\/$\MM^{++}$ $+$ there are class many supercompact cardinals'' (or even class many 
  extendible cardinals) is consistent with \GA.
\end{Thm}
\prf  Sean Cox \cite{cox2} proved that $\MM^{++}$ is preserved by $\omega_2$-directed 
closed forcing (Theorem 4.7 in \cite{cox2}). Starting from a model with cofinally many 
supercompact cardinals, use the first supercompact to force $\MM^{++}$. Then the class 
forcing just like that in the proof of \Thmof{p-hierarchies-6} (or 
like the one in \cite{goldberg}) will produce a desired model. \qedofThm\qedskip

\begin{Cor}\Label{p-Lg-RcA-2} $\MM^{++}$ or even $\MM^{++}$ $+$ ``there are class many 
  super compact cardinals'' does not imply that the continuum is a tightly $\calP$-Laver 
  gen.\ ultrahuge cardinal for any of the large enough subclass $\calP$ of the class of all 
  semiproper \pos.  
\end{Cor}
\prf Let $\calP=$ semiproper \pos. 
Note that, if $\kappa$ is
$\calP$-Laver generically supercompact, then $\kappa=\continuum$ follows (see e.g.\ Theorem 
5 and Lemma 6 in Fuchino \cite{janos}).

If $\kappa$ is the tightly $\calP$-Laver generically ultrahuge continuum, then 
\Thmof{p-Lg-RcA-0} together with each one of the Propositions \ref{p-rec-GA-1}, \ref{p-rec-GA-1-0}, 
\ref{p-rec-GA-1-1-0} 
implies that \GA\ does not hold. Thus the model of $\MM^{++}$ $+$ ``there are class many 
super compact cardinals'' $+$ \GA\ of \Thmof{p-Lg-RcA-1-0} witnesses the desired 
non-implication. \qedofCor\qedskip

Note that \Corof{p-Lg-RcA-2} with ``tightly $\calP$-Laver gen.\ ultrahuge'' replaced by ``tightly
$\calP$-Laver gen.\ hyperhuge'' is trivial. This is because consistency strength of the 
existence of the tightly $\calP$-Laver gen.\ hyperhuge cardinal is known to be that of the 
existence of a (genuinely) hyperhuge cardinal (see the remark right after \Propof{p-rec-GA-1-1-1}).

Our Theorems \ref{p-genabs-rec-0} and \ref{p-genabs-Laver-1} generalize Viale's 
\Thmof{p-theintro-a} in terms of possible instances of the class $\calP$ not covered by
\Thmof{p-theintro-a} and also in terms of the cardinal $\kappa$ in the conclusion of the 
theorems which  
can be strictly bigger than $\aleph_2$ (which can really happen if e.g.\ 
 $\calP$ is the class of ccc \pos). 

On the other hand, the premise of Viale's \Thmof{p-theintro-a} is consistent with \GA\ 
by \Thmof{p-Lg-RcA-1-0}
while this is not the case with \Thmof{p-genabs-rec-0} for many natural instances 
of $\calP$ by Propositions \ref{p-rec-GA-1}, \ref{p-rec-GA-1-0}, \ref{p-rec-GA-1-1-0} and 
unclear in case of $\calP=$ stationary preserving \pos\ with \Thmof{p-genabs-Laver-1}. 

Viale's \Thmof{p-theintro-a} implies, in particular: 
\begin{Cor}{\rm(to Viale's \Thmof{p-theintro-a} and \Thmof{p-Lg-RcA-1-0})}\Label{p-misc-0}
  The assertion 
  \begin{xitemize}
  \item[] $\calH(\aleph_2)^\uniV\prec_{\Sigma_2}\calH(\aleph_2)^{\uniV[\genG]}$\quad 
    for any 
    stationary preserving \po\ $\poP$ with\\ $\forces{\poP}{\BMM}$, and $(\uniV,\poP)$-generic $\genG$. 
  \end{xitemize}
  is consistent with \GA. \qed
\end{Cor}

Concerning  \Thmof{p-genabs-Laver-1}, it is open at the moment if the 
existence of a tightly $\calP$-Laver-gen.\ huge cardinal is inconsistent with \GA. However 
some of its strengthenings do contradict \GA\ for many instances of the class $\calP$ 
of \pos\ as we saw in \Propof{p-rec-GA-1-2} and \Thmof{p-rec-GA-2}. 

The positive answer to the following question would give a clear separation of 
Laver-genericity from the corresponding forcing axiom with double plus: 

\begin{Problem}
  Does the (tightly) $\calP$-Laver gen.\ supercompact cardinal axiom
  (i.e., the existential statement of such a cardinal, e.g.\ for $\calP$ as in 
  \Corof{p-hierarchies-7}) imply the negation of \GA? 
\end{Problem}

Though \Corof{p-Lg-RcA-2} makes the positive answer to the following problem rather unpromising, 
Theorem 2.53 of Woodin \cite{woodin-book} and its variants (e.g.\ Theorem 4.5 in 
\cite{cox2}) seem to suggest a positive answer:

\begin{Problem}
  Is there any reasonable assumption under which $\MM^{++}$ and (tightly) $\calP$-Laver 
  gen.\ supercompact cardinal axiom are equivalent?
\end{Problem}

{\ifextended\extendedcolor\mbox{}\bigskip

  The next subsection has been removed from the version for publication of the current paper as 
  we have decided to move the material to the subsequent paper which is currently 
  in preparation. 
\subsection{Yet another hierarchy of restricted Recurrence Axioms}\Label{Yah}

The hierarchy $\MP(\cdot,\cdot)_{\cdot}$ of Maximality Principles ($\Leftrightarrow$ 
Recurrence Axioms) introduced in \sectionof{hierarchies} is 
suitable for the analysis of consistency strength and strictness of the hierarchy of 
restricted form of these principles.

{}

\begin{Lemma}\Label{p-Yah-0} Suppose that $\calP$ is an iterable class of \pos. \smallskip
  
  \wassert{1} $(\calP,A)_\Gamma$-\RcAp\ \ $\Rightarrow$\ \ 
  $\MP^*(\calP,A)_\Gamma$.\smallskip

  \wassert{2} If $A\subseteq A'$, and $\Gamma\subseteq\Gamma'$ then 
  $(\calP,A')_{\Gamma'}$-\RcAp\ \ $\Rightarrow$\ \ $(\calP,A)_\Gamma$-\RcAp, and 
  $\MA^*(\calP,A')_{\Gamma'}$\ \ $\Rightarrow$\ \ $\MA^*(\calP,A)_{\Gamma}$.\smallskip

  \wassert{3} If $\calP\subseteq\calP'$, $A\subseteq A'$, and $\Gamma\subseteq\Gamma'$ then 
  $(\calP',A')_{\Gamma'}$-\RcA\ \ $\Rightarrow$\ \ $(\calP,A)_\Gamma$-\RcA.

  \wassert{4} $\MP^*(\calP,A)_{\Pi_1}$ holds (in \ZFC).\smallskip

  \wassert{5} $(\calP,A)_{\Sigma_1}$-\RcA\ $\Leftrightarrow$\ $(\calP,A)_{\Sigma_1}$-\RcAp\
  $\Leftrightarrow$\ $\MP^*(\calP,A)_{\Sigma_1}$. 
\end{Lemma}
\prf \assertof{1}: Suppose that $(\calP,A)_\Gamma$-\RcAp\ holds. Let $\varphi$ be a 
provably $\calP$-persistent formula in $\Gamma$ and $\ol{a}\in A$. If $\varphi(\ol{a})$ is 
forced by $\poP\in\calP$, then $\poP$ is a push of the $\calP$-button
$\forall\ul{\poP}(\forces{\ul{\poP}}{\varphi(\ol{a})})$ by $(\varphi)^\ast_\poP$ which is 
provable by assumption. 
By $(\calP,A)_\Gamma$-\RcAp, it follows that $\varphi(\ol{a})$ holds (in $\uniV$). This 
shows that \xitemsubof{p-Yah-0}{\varphi} holds. 
\smallskip

\assertof{2}, \assertof{3}: by definitions. \smallskip

\assertof{4}: $(\calP,A)_{\Pi_1}$-\RcAp\ holds by \Lemmaof{p-hierarchies-3}. Thus 
\assertof{1} implies $\MP^*(\calP,A)_{\Pi_1}$.\smallskip

\assertof{5}: The first equivalence is a part of \Thmof{p-theintro-1}. The second 
equivalence is proved from this and argument similar to the proof of \Thmof{p-theintro-1}. 
\qedofLemma\qedskip

The monotonicity \Lemmaof{p-Yah-0},\,\assertof{3} is used in Fuchino \cite{janos} in the 
argument to single out the maximal instance of Recurrence Axiom
$(\mbox{stationary preserving},\calH(\kappa_\refl))$-\RcA\ under $\continuum=\aleph_2$ and the 
other maximal instance $(\mbox{all \pos}, \calH(\continuum))$-\RcA\ under \CH.

\Lemmaof{p-Yah-0},\,\assertof{4} and \assertof{5} shows that the list of equivalent 
assertions in \Corof{p-hierarchies-2} and \Lemmaof{p-hierarchies-3} can also include
$\MP^*(\calP,\calH(\kappa))_{\Sigma_1}$ and $\MP^*(\calP,A)_{\Pi_1}$ respectively. 

$\MP^*(\cdots)_\Gamma$ is almost identical with $(\cdots)_\Gamma$-\RcAp.

\begin{Prop}\Label{p-Yah-1} Suppose that $\calP$ is an iterable class of \pos\ and $A$ is 
  any set. \smallskip

  \wassert{1} If $\calP$ is $\Sigma_m$-definable then for any natural number $n$ with
  $\max\ssetof{4,m}\leq n$, we 
  have $(\calP,A)_{\Sigma_n}$-\RcAp\ $\Leftrightarrow$\ $\MP^*(\calP,A)_{\Sigma_n}$. 
  \smallskip

  \wassert{2} $(\calP,A)$-\RcAp\ ($\Leftrightarrow$\ $\MP(\calP,A)$)\ $\Leftrightarrow$\
  $\MP^*(\calP,A)$.\smallskip

  \wassert{3} If $\calP$ is $\Sigma_4$-definable, and one of the  
  conditions in \Propof{p-rec-GA-1-1-0} holds, then $\MP^*(\calP,\emptyset)_{\Sigma_4}$ 
  implies $\neg\GA$ (cf.\ \Thmof{p-hierarchies-6-1}).
\end{Prop}
\prf \assertof{1}: Suppose that $\calP$ is $\Sigma_n$-definable iterable class of \pos, and
$n\geq\max\ssetof{3,m}$. If $(\calP,A)_{\Sigma_n}$-\RcAp\ holds, then we have
$\MP^*(\calP,A)_{\Sigma_n}$ by \Lemmaof{p-Yah-0},\,\assertof{1}.

Assume now that $\MP^*(\calP,A)_{\Sigma_n}$ holds, and suppose 
that $\Sigma_n$ formula $\varphi=\varphi(\ol{x})$ and $\ol{a}\in A$ are \st\
$\forces{\poP}{\varphi(\ol{a})}$ for a $\poP\in\calP$.

With the proof of definability of grounds in mind (see e.g.\ \cite{reitz}), let
$\varphi^*(\ol{x})$ be the formula claiming \memo{see math-notes-20\\Bukovsk\'y's and 
  Laver's theorems revisited}
\begin{xitemize}
\item[] $\exists X\,(\,{}
  \begin{array}[t]{@{}l}
    \mbox{``\,}X\mbox{ is the parameter which codes a }\calP\mbox{-ground''}\,\\[\jot]\land\,
    \mbox{``}\ol{x}\in\mbox{ the }\calP\mbox{-ground coded by }X\,\\[\jot]\land\,
    \mbox{``the ground coded by }X\mbox{ satisfies }\varphi(\ol{x})\mbox{''}). 
  \end{array}$
\end{xitemize}
By the choice of $n$, $m$, $\calP$, $\varphi^*$ can be expressed as $\Sigma_n$-formula and,  
it is provably $\calP$-persistent since $\calP$ is iterable. We also have $\forces{\poP}{\varphi^*(\ol{a})}$.
By $\MP^*(\calP,A)_{\Sigma_n}$, it follows that $\uniV\models\varphi^*(\ol{a})$. By 
definition of $\varphi^*$, this means that there is a $\calP$-ground $\uniW$ of $\uniV$ 
\st\ $\uniW\models\varphi(\ol{a})$.

This shows that $(\calP,A)_{\Sigma_n}$-\RcAp\ holds. 
\smallskip

\assertof{2}: follows from \assertof{1} (the first equivalence in parentheses is by 
\Propof{p-intro-1}, \assertof{1}). \smallskip

\assertof{3}: By \assertof{1} and \Propof{p-rec-GA-1-1-0}. 
\qedofProp\qedskip

%% The hierarchy $\MP^*(\calP,\calH(\kappa))_{\Sigma_n}$, $n\in\omega$ is a special cases of 
%% the hierarchy of forcing axioms introduced in Asperó \cite{aspero}:

For a class of \pos\ $\calP$, a set $\Gamma$ of formulas, and infinite cardinals $\kappa$, 
$\lambda$ with $\kappa\leq\lambda$, the principle $\BFA(\calP,\Gamma)_{\kappa,\lambda}$ is 
defined by

\begin{xitemize}
\item[{\darkred$\BFA(\calP,\Gamma)_{\kappa,\lambda}$}\,: ] For any $\ol{a}\in\calH(\lambda)$ and
  $\varphi=\varphi(\ol{x})\in\Gamma$ with $\forces{\poP}{\varphi(\ol{a})}$ for some
  $\poP\in\calP$, there are stationarily many $X\in[\calH(\lambda)]^{\LT\kappa}$ \st\ $X$ 
  (as an $\in$-model) satisfies the Axiom of 
  Extensionality, $\ol{a}\in X$, and $\uniV\models\varphi(m_X(\ol{a}))$ where $m_X$ denotes 
  the Mostowski collapse of $X$. 
\end{xitemize}

Note that $X$ as above should satisfy the Axiom of Extensionality so that the Mostowski 
collapse of $X$ is defined (and injective). 
$\BFA(\calP,\Gamma)_{\kappa,\lambda}$ was introduced in Asperó \cite{aspero}. 
Note that this principle is not formulated as a generalization of 
the original definition of $\BFA_{\LT\kappa}(\calP)$ (see  
\sectionof{theintro}) but rather a statement which stands in analogy with the 
property in Bagaria's Absoluteness \Thmof{p-theintro-0} characterizing $\BFA_{\LT\kappa}(\calP)$. 

\begin{Lemma}\Label{p-Yah-2} {\rm (Goodman \cite{goodman})} Suppose that $\calP$ is a class 
  of \pos\ and $\Gamma$  
  a set of formulas $\subseteq\Lin$. \wassertof{1} $\BFA(\calP,\Gamma)_{\kappa,\kappa}$ 
  holds if and only if 
  \begin{xitemize}
  \xitem[x-Yah-1] For any $\ol{a}\in\calH(\kappa)$ and $\varphi=\varphi(\ol{x})\in\Gamma$ 
    with $\forces{\poP}{\varphi(\ol{a})}$, $\uniV\models\varphi(\ol{a})$ holds.
  \end{xitemize}
  
  \wassert{2} $\BFA(\calP,(\Gamma)^*_\calP)_{\kappa,\kappa}$\ $\Leftrightarrow$\ 
  $\MP^*(\calP,\calH(\kappa))_\Gamma$ holds. 
\end{Lemma}
\prf \assertof{1}: Suppose that $\BFA(\calP,\Gamma)_{\kappa,\kappa}$ holds. For
$\ol{a}\in\calH(\kappa)$ and $\varphi=\varphi(\ol{x})\in\Gamma$, suppose that
$\forces{\poP}{\varphi(\ol{a})}$. Then there are stationarily many
$X\in[\calH(\kappa)]^{\LT\kappa}=\calH(\kappa)$ \st\ $\ol{a}\in X$ and
$\uniV\models\varphi(m_x(\ol{a}))$. In particular, there is such $X$ that
$\trcl^+(\ol{a})\subseteq X$. Then $m_X(\ol{a})=\ol{a}$ and $\uniV\models\varphi(\ol{a})$.
Thus \xitemof{x-Yah-1} holds. 

Conversely, if \xitemof{x-Yah-1} holds, then
$\calX:=\setof{X\in[\calH(\kappa)]^{\LT\kappa}}{\trcl^+(\ol{a})\subseteq X}$ is stationary 
(actually club) in $[\calH(\kappa)]^{\LT\kappa}$. For each $X\in\calX$, we have
$m_X(\ol{a})=\ol{a}$. This shows that $\BFA(\calP,\Gamma)_{\kappa,\kappa}$ holds. \smallskip

\assertof{2}: follows from \assertof{1}.
\qedofLemma\qedskip

The following argument leading to \Propof{p-Yah-6} and \Propof{p-Yah-9}\memox{\darkred!!!} is mostly a combination of ideas already 
utilized somewhere in Asperó\  
\cite{aspero}, {Goodman} \cite{goodman}, and/or Hamkins \cite{hamkins}. We include the 
details of the proofs here for convenience of the reader.

\begin{Lemma}\Label{p-Yah-3}
  Suppose that $\MP^*(\calP,\calH(\kappa_\refl))_{\Sigma_2}$ holds and $\continuum$ is regular. If, 
  either \wassertof{a} there is a \po\ in $\calP$ collapsing $\kappa_\refl$\,, or 
  \wassertof{b} %% $\neg\CH$ holds and 
  there is a \po\ in $\calP$ adding ${\kappa_\refl}^+$ reals 
  without collapsing ${\kappa_\refl}^+$, then ${\kappa_\refl}^\uniV$ is inaccessible 
  in $\uniL$. 
\end{Lemma}
\prf Assume that $\MP^*(\calP,\calH(\kappa_\refl))_{\Sigma_2}$ holds.

%% First, suppose that \assertof{a} holds.
%% %% Note that $\kappa_\refl$ is a regular cardinal in $L$. 
%% If 
%% $\mu<\kappa_\refl$ is a cardinal in $\uniL$, 
Consider the sentence $\varphi(x)$ saying 
\begin{xitemize}
\item[] 
  $\exists\vsymb{\mu}'\,(\uniL\modelof{\vsymb{\mu}'\mbox{ is a cardinal}}\,\land\,
  x<\vsymb{\mu}'<\kappa_\refl)$. 
\end{xitemize}
$\varphi(x)$ is $\Sigma_2$ and it is provably $\calP$-persistent. 
\memo{``$\ul{\kappa}=\kappa_\refl$'' is $\Sigma_2$ (\mbox{Scan\_2024-05-14--12.42 gappo - 
    annotated} p.95)} 

In both of the cases, there is 
$\poP\in\calP$ \st\ $\forces{\poP}{\cardof{{\kappa_\refl}^\uniV}<\kappa_\refl}$. 
Suppose that $\mu<\kappa_\refl$ is a cardinal in $\uniL$. 

Since 
$\kappa_\refl$ is a regular cardinal in $\uniL$, it follows that $\forces{\poP}{\varphi(\mu)}$. 

By $\MP^*(\calP,\calH(\kappa_\refl))_{\Sigma_2}$, it follows 
that $\uniV\models\varphi(\mu)$. I.e., $(\mu^+)^\uniL<\kappa_\refl$. 
\qedofLemma

\begin{Lemma}\Label{p-Yah-4}
  Suppose that $\continuum$ is regular and $\calP$ is a class of \pos\ \st\ either 
  \wassertof{a'} %% $\CH$ holds and, 
  for any cardinal $\lambda\geq\kappa_\refl$, there is a \po\/ $\poP\in\calP$ which collapses
  $\lambda$ to cardinality $\aleph_1$\,, or  
  \wassertof{b'} %% $\neg\CH$ holds and, 
  for cofinally many cardinals $\lambda>\continuum$, there is
  $\poP\in\calP$ adding at least $\lambda$ many reals without collapsing $\lambda$. 
  Then, for any $n\geq 2$, $\MP^*(\calP,\calH(\kappa_\refl))_{\Sigma_n}$ implies
  $\uniL_{{\kappa_\refl}^\uniV}\prec_{\Sigma_n}\uniL$. 
\end{Lemma}
\prf For $n=2$, $\MP^*(\calP,\calH(\kappa_\refl))_{\Sigma_n}$ implies that $\kappa_\refl$ 
is an inaccessible cardinal in $\uniL$ by \Lemmaof{p-Yah-3}. Thus
$\uniL\models L_{{\kappa_\refl}^\uniV}=\calH({\kappa_\refl}^\uniV)$. By 
\xitemof{x-theintro-0}, it follows 
that 
\begin{xitemize}
\xitem[x-Yah-1-0] 
  $\uniL_{{\kappa_\refl}^\uniV}\prec_{\Sigma_1}\uniL$.
\end{xitemize}

Thus it is enough to show that $\uniL_{{\kappa_\refl}^\uniV}\prec_{\Sigma_{n}}\uniL$ holds  
for $n\geq 2$ assuming 
\begin{xitemize}
\xitem[x-Yah-2] 
  $\MP^*(\calP,\calH(\kappa_\refl))_{\Sigma_n}$, and
\xitem[x-Yah-3] 
  $\uniL_{{\kappa_\refl}^\uniV}\prec_{\Sigma_{n-1}}\uniL$.
\end{xitemize}
Note that for $n=2$, \xitemof{x-Yah-3} is just \xitemof{x-Yah-1-0}.

To prove this claim, assume \xitemof{x-Yah-2} and \xitemof{x-Yah-3}, and 
let $\delta:={\kappa_\refl}^\uniV$, 
$\psi=\psi(\ol{x},\ol{y})$ a $\Pi_{n-1}$-formula, and $\ol{a}\in\uniL_\delta$.
Since $\delta$ is a limit ordinal, there is $\delta_0<\delta$ \st\
$\ol{a}\in L_{\delta_0}$.

If $\uniL_\delta\models\exists\ol{y}\psi(\ol{a},\ol{y})$, then
$\uniL\models\exists\ol{y}\psi(\ol{a},\ol{y})$ by \xitemof{x-Yah-3}. 

Assume now that $\uniL\models\exists\ol{y}\psi(\ol{a},\ol{y})$. Let $\eta=\eta(u)$  be the 
$\Sigma_n$-formula
\begin{xitemize}
\item[] $\exists v\,(\mbox{``}v=\uniL_{\kappa_\refl}\mbox{''}\,\land\,
  (\forall x\in u)\,(\mbox{``}\uniL\models\exists\ol{y}\psi(\ol{x},\ol{y})\mbox{''} 
  \rightarrow(\exists\ol{y}\in v)\,(\mbox{``}\uniL\models\psi(\ol{x},\ol{y})\mbox{''})))$.
\end{xitemize}

This formula is clearly provably $\calP$-persistent. By the assumption 
\assertof{a'} or \assertof{b'}, and by an argument similar to the proof of  
\Lemmaof{p-Yah-3}, it follows (in both of the cases \assertof{a'} and \assertof{b'}) that there is
$\poP\in\calP$ \st\ $\forces{\poP}{\eta(\uniL_{\delta_0})}$. Thus by \xitemof{x-Yah-2},
$\uniV\models\eta(\uniL_{\delta_0})$.

In particular, there is $\ol{b}\in\uniL_{\delta}$ \st\ $\uniL\models\psi(\ol{a},\ol{b})$.
By \xitemof{x-Yah-3}, it follows that $\uniL_\delta\models\psi(\ol{a},\ol{b})$ and thus
$\uniL_\delta\models\varphi(\ol{a})$. 
\qedofLemma\qedskip

The following proposition can be seen as a subset of {\sc Theorem} 2.6 in 
\cite{aspero}{\ifextended\extendedcolor\ and 
the proof given here 
is also more or less identical with the one in \cite{aspero}\,\fi}:
\begin{Prop}\Label{p-Yah-5}
  For a natural number $n$, assume that $\lambda$ is an infinite cardinal and 
$\kappa\geq\lambda$ is an inaccessible $\Sigma_n$-correct cardinal (i.e. a 
  cardinal $\kappa$ with the property $V_\kappa (=\calH(\kappa))\prec_{\Sigma_n}\uniV$).

  Suppose that $\calP$ is a $\Sigma_n$-definable iterable class of \pos\ \st
  \begin{xitemize}
  \xitem[x-Yah-4] all\/ $\poP\in\calP$ preserve cardinals $\leq\lambda$, and if
    $\lambda<\kappa$, then the equation\medskip\\
    \qquad 
    $\lambda=\max\setof{\mu\in\Card}{\mbox{all\/ }\poP\in\calP
      \mbox{\/ preserves cardinals }\leq\mu}$\medskip\\
    is 
    provable (in \ZFC);\/\footnotemark
  \xitem[x-Yah-5] $\calP$ admits iteration 
    $\seqof{\poP_\alpha,\utpoQ_\beta}{\alpha\leq\kappa,\beta<\kappa}$ of length $\kappa$ 
    with some appropriate kind of iteration
    (e.g.\ either FS-, CS-, or Easton-support iteration) \st\ 
    \begin{itemize}
    \xxitem[x-Yah-5][a] $\poP_\kappa$ is the direct limit of
      $\seqof{\poP_\alpha}{\alpha<\kappa}$. 
    \xxitem[x-Yah-5][b]       
      $\poP_\kappa\in\calP$;
    \xxitem[x-Yah-5][c]
      $\forces{\poP_\alpha}{\poP_\kappa/\utgenG_\alpha\in\calP}$ for all $\alpha<\kappa$;
    \xxitem[x-Yah-5][d]
      $\poP_\alpha\in V_\kappa$ for all $\alpha<\kappa$, and\/ $\poP_\kappa$ satisfies the
      $\kappa$-cc.
    \end{itemize}
  \end{xitemize}

  %% Let
  %% \begin{xitemize}
  %% \item[] 
  %%   $\lambda^{(+)}=\left\{\,{}
  %%   \begin{array}{@{}ll}
  %%     \lambda^+, &\mbox{if }\lambda<\kappa;\\[\jot]
  %%     \kappa\ (=\lambda), &\mbox{otherwise.}
  %%   \end{array}
  %%   \right.$ 
  %% \end{xitemize}

  Then there is an iteration
  $\seqof{\poP_\alpha,\utpoQ_\beta}{\alpha\leq\kappa,\beta<\kappa}$ satisfying 
  \xitemof{x-Yah-5} above \st\
  $\forces{\poP_\kappa}{\MP^*(\calP,\calH(\kappa))_{\Sigma_n}}$. \ifextended\else\qed\fi 
\end{Prop}
\footnotetext{In this case, we assume that  $\lambda$ is definable e.g.\ as
  $\aleph_1$ so that we can formulate the condition without introducing a new constant symbol.}
{\ifextended\extendedcolor
\prf Let $\mapping{f}{\kappa}{\kappa\times\omega\times\kappa}$ be a surjection \st\ each
$\pairof{\alpha_0,n\alpha_1}\in\kappa\times\omega\times\kappa$ appears $\kappa$ times as a 
value of $f$. Let $\seqof{\varphi_m}{m\in\omega}$ be an enumeration 
of $(\Sigma_n)^*_\calP$ (as a set in \ZFC) corresponding to the recursive enumeration of
$(\Sigma_n)^*_\calP$ in meta-mathematics. 

Let $\seqof{\poP_\alpha,\utpoQ_\beta}{\alpha\leq\kappa,\beta<\kappa}$ be an iteration in
$\calP$ satisfying \xitemof{x-Yah-5} defined along with the sequences 
$\seqof{\utab^\alpha_\xi}{\xi<\kappa}$ in $V_\kappa$ for each $\alpha<\kappa$ \st\
\begin{xitemize}
\xitem[x-Yah-6]
  $\forces{\poP_\alpha}{\setof{\utab^\alpha_\xi}{\xi<\kappa}^\bullet=\calH(\kappa)}$.\footnotemark  
\end{xitemize}
\footnotetext{With the superscript bullet "$\cdots^\bullet$" in connection with forcing with a \po\
  $\poP$, we denote the canonical $\calP$-name corresponding to the object ``$\cdots$'' 
  describes where we assume that $\poP$-names are introduced as in Kunen \cite{kunen-2011}.}

The successor step of the iteration is set by the following:
\begin{xitemize}
\xitem[x-Yah-7] If $f(\alpha)=\pairof{\alpha_0,m,\alpha_1}$, $\alpha_0\leq\alpha$, and\\[\jot]
  \qquad$\forces{\poP_\alpha}{(\exists\vsymb{\poQ}\in\calP)\,(\vsymb{\poQ}\in\calH(\kappa)\,\land\,
  \forces{\vsymb{\poQ}}{\varphi_m(\utab^{\alpha_0}_{\alpha_1})})}$,\footnotemark\\[\jot]
  then $\utpoQ_\alpha\in V_\kappa$ is such a $\poP_\alpha$-name as $\vsymb{\poQ}$ as above;\\[\jot]
  otherwise $\utpoQ_\alpha=\ssetof{\bbone}^\bullet$.
\end{xitemize}
\footnotetext{Strictly speaking we mean with $\utab^{\alpha_0}_{\alpha_1}$ the ($\vsymb{\poQ}$-check 
  names) of the $\poP_\alpha$-names corresponding to the $\poP_{\alpha_0}$-names. }

We show that this iteration with $\poP_\kappa$ is as desired. Let $\genG_\kappa$ be a
$(\uniV,\poP_\kappa)$-generic filter. 

Suppose that $\utab$ is a 
tuple of $\poP_\kappa$-names of elements of $\calH(\kappa)^{\uniV[\genG_\kappa]}$
($=\calH(\kappa)^{\uniV[\genG_\kappa]}$), and $\utpoQ$ be $\poP_\kappa$-name of \po\ in 
$\calP$ and $\forces{\poP_\kappa\ast\utpoQ}{\varphi(\utab)}$ for a (concretely given)
provably $\calP$-persistent $\Sigma_n$-formula $\varphi$. By definition of the sequence 
$\seqof{\varphi_m}{m\in\omega}$ there is a number $m^*$ \st\ $\varphi=\varphi_{m^*}$. 
By \xxitemof{x-Yah-5}{d}, there is $\gamma<\kappa$ \st\ $\utab$ corresponds 
to $\poP_\gamma$-names which we shall also denote with $\utab$. Thus, there is 
$\beta^*<\kappa$ \st\ $\utab=\utab^\gamma_{\beta^*}$. 

Note that for all $\alpha\in\kappa\setminus\gamma$ we have 
\begin{xitemize}
\xitem[x-Yah-8]
  $\forces{\poP_\alpha}{\forces{(\poP_\kappa/\utgenG_\alpha)\ast\utpoQ}{\varphi(\utab)}}$
  and 
  $\forces{\poP_\alpha}{\poP_\kappa/\utgenG_\alpha\ast\utpoQ\in\calP}$
\end{xitemize}
by \xxitemof{x-Yah-5}{c} and by iterability of $\calP$. 

Let $\alpha^*\in\kappa\setminus\gamma$ be \st\ 
\begin{xitemize}
\xitem[x-Yah-9] 
  $f(\alpha^*)=\pairof{\gamma,n^*,\beta^*}$. 
\end{xitemize}
By \xitemof{x-Yah-8} we have   $\forces{\poP_{\alpha^*}}{
  \forces{(\poP_\kappa/\utgenG_{\alpha^*})\ast\utpoQ}{\varphi_{m^*}(\utab^\gamma_{\beta^*})}}$, and
$\forces{\poP_{\alpha^*}}{\poP_\kappa/\utgenG_{\alpha^*}\ast\utpoQ\in\calP}$. 

Thus we have
\begin{xitemize}
\xitem[x-Yah-10]
  $\uniV\models\exists\vsymb{\utpoQ}\,(\forces{\poP_{\alpha^*}}{\vsymb{\utpoQ}\in\calP
  \ \land\ \forces{\vsymb{\utpoQ}}{\varphi_{m^*}(\utab^\gamma_{\beta^*})}})$.
\end{xitemize}
Since the property in \xitemof{x-Yah-10} is formalizable in $\Sigma_n$, and $\kappa$ is
$\Sigma_n$-correct, it follows that 
\begin{xitemize}
\xitemx[x-Yah-11]
  $V_\kappa\models
  \exists\vsymb{\utpoQ}\,(\forces{\poP_{\alpha^*}}{\vsymb{\utpoQ}\in\calP
    \ \land\ \forces{\vsymb{\utpoQ}}{\varphi_{m^*}(\utab^\gamma_{\beta^*})}})$.
\end{xitemize}
By \xitemof{x-Yah-7} it follows that
$\forces{\poP_{\alpha^*}}{\forces{\utpoQ_{\alpha^*}}{\varphi(\utab)}}$. Since $\varphi$ is 
provably $\calP$-persistent, it follows by \xxitemof{x-Yah-5}{c} that
$\forces{\poP_{\alpha^*}}{\forces{\poP_\kappa/\utgenG_{\alpha^*}}{\varphi(\utab)}}$. This 
is equivalent to $\forces{\poP_\kappa}{\varphi(\utab)}$. 
\qedofProp\fi}

\begin{Lemma}{\rm ({\sc Fact 2.2} in Asperó \cite{aspero})}\Label{p-Yah-5-0}
  If $\kappa$ is inaccessible and $\Sigma_{n+1}$-correct for some $n\geq 1$ then there are 
  unboundedly many inaccessible $\Sigma_n$-correct cardinals below $\kappa$. \qed
\end{Lemma}

\begin{Prop}\Label{p-Yah-6} Assume that $n\geq 2$. 
  Suppose $\calP$ is a $\Sigma_n$-definable iterable class of \pos\ \st\ $\calP$
  satisfies the condition \assertof{a'} or \assertof{b'} of \Lemmaof{p-Yah-4}, and it is also 
  provable that $\calP$ satisfies \xitemof{x-Yah-4} and \xitemof{x-Yah-5} of \Propof{p-Yah-5} for 
  any infinite cardinal $\lambda$ and inaccessible $\kappa\geq\lambda$. 
 
  Then, assuming the consistency of the theory 
  \begin{xitemize}
  \xitem[x-Yah-10-0] 
    \ZFC\ $+$ ``\/$\continuum$ is regular'' $+$ $\MP^*(\calP,\calH(\kappa_\refl))_{\Sigma_n}$, 
  \end{xitemize}
  this theory does not imply
  $\MP^*(\calP,\calH(\kappa_\refl))_{\Sigma_{n+1}}$. 
\end{Prop}
\prf Suppose, toward a contradiction, that \ZFC\ proves that\\
$\MP^*(\calP,\calH(\kappa_\refl))_{\Sigma_n}$ implies 
$\MP^*(\calP,\calH(\kappa_\refl))_{\Sigma_{n+1}}$.

Working in the theory \ZFC\ $+$
$\MP^*(\calP,\calH(\kappa))_{\Sigma_n}$, we also have
$\MP^*(\calP,\calH(\kappa))_{\Sigma_{n+1}}$ by the assumption. By \Lemmaof{p-Yah-4} applied to 
$n+1$, and then by \Lemmaof{p-Yah-5-0}, combined 
with Downward Löwenheim-Skolem Theorem and Mostowski Collapsing Lemma, we obtain a 
countable transitive model $M$ of \ZFC\ $+$ $\exists\kappa\,(\kappa\mbox{ is inaccessible}\,\land\,
V_\kappa\prec_{\Sigma_n}\uniV)$. 

By \Propof{p-Yah-5}, there is a generic extension $M[\genG]$ which is a model of \ZFC\,$+$\,
$\MP^*(\calP,\calH(\kappa))_{\Sigma_n}$.\footnote{The proof of \Propof{p-Yah-5}  
  is written as a proof of $\uniV[\genG_\kappa]\models\MP^*(\cdots)_{\Sigma_n}$
  for the axiom scheme $\MP^*(\calP,\calH(\kappa_\refl))_{\Sigma_n}$ in the sense of 
  meta-mathematics, but for the set model $M$  
  this proof can be almost verbosely adopted to prove $M[\genG_\kappa]\models\MP^*(\cdots)_{\Sigma_n}$ 
  for the axiom scheme $\MP^*(\calP,\calH(\kappa_\refl))_{\Sigma_n}$ in the sense of \ZFC.} 
By \assertof{a'} and \assertof{b'} we have 
$M[\genG]\models\kappa=\kappa_\refl$.

Thus we obtained a proof of 
$consis(\ZFC\,+\,\MP^*(\calP,\calH(\kappa_\refl))_{\Sigma_n})$ in
\ZFC\ $+$ $\MP^*(\calP,\calH(\kappa_\refl))_{\Sigma_n}$.
This is a contradiction by   
The Second Incompleteness Theorem. \qedofProp

\begin{Cor}\Label{p-Yah-6-0} Under the same assumption on $n$ and $\calP$ as in 
  \Propof{p-Yah-6}, if
  \begin{xitemize}
  \xitemd[x-Yah-10-0]{'}
    \ZFC\ $+$ ``\/$\continuum$ is regular'' $+$ $(\calP,\calH(\kappa_\refl))_{\Sigma_n}$-\RcAp, 
  \end{xitemize}
  is consistent then this theory does not prove $(\calP,\calH(\kappa_\refl))_{\Sigma_{n+1}}$-\RcAp, 
\end{Cor}
\prf By \Propof{p-Yah-1} and \Propof{p-Yah-6}.\qedofCor\qedskip

In case of $n\geq 3$ \Propof{p-Yah-6} can be further improved (see \Propof{p-Yah-9} below). 

\begin{Lemma}\Label{p-Yah-7} 
  Suppose $\calP$ is a class of \pos\ satisfying \assertof{a'} or 
  \assertof{b'} of \rlap{\Lemmaof{p-Yah-4}.}\smallskip

  \wassert{1} 
  If $n\geq 3$, $\MP^*(\calP,\emptyset)_{\Pi_n}$ 
  implies that there are unboundedly many $\Sigma_{n-1}$-correct inaccessible cardinals in
  $\uniL$.\smallskip 
  
  \wassert{2} Suppose that $\MP^*(\calP,\emptyset)_{\Pi_2}$ holds. If there is at least one  
  inaccessible cardinal, then there are cofinally many inaccessible cardinals.
\end{Lemma}
\prf \assertof{1}: Assume, toward a contradiction that 
\begin{xitemize}
\xitem[x-Yah-11] $\MP^*(\calP,\emptyset)_{\Pi_n}$ holds 
\end{xitemize}
but the class $\calB$ of all $\Sigma_{n-1}$-correct cardinals in $\uniL$ is bounded.

By \assertof{a'} or \assertof{b'} there is a \po\ $\poP\in\calP$ \st\ 
\begin{xitemize}
\item[] $\forces{\poP}{\kappa_\refl>\calB}$. 
\end{xitemize}
The forced statement is $\Pi_n$ and it is provably $\calP$-persistent. Thus, by \xitemof{x-Yah-11},
$\uniV\modelof{\kappa_\refl>\calB}$. This is a contradiction to \Lemmaof{p-Yah-4} with $n$ 
replaced by $n-1$. \smallskip

\assertof{2}: Suppose that $\MP^*(\calP,\emptyset)_{\Pi_2}$ holds but there are only
boundedly many but at least one inaccessible cardinals. Then, by the condition on $\calP$, 
there is $\poP\in\calP$ \st\ 
\begin{xitemize}
\xitem[x-Yah-11-0] 
 $\forces{\poP}{\mbox{there are no inaccessible cardinals}}$.
\end{xitemize}
Since the statement in \xitemof{x-Yah-11-0} is $\Pi_2$ and provably $\calP$-persistent. It 
follows that\\
$\uniV\modelof{\mbox{there are no inaccessible cardinals}}$. This is a contradiction. 
\qedofLemma

\begin{Lemma}\Label{p-Yah-8}
  For $n\geq 2$, if the theory 
  \begin{xitemize}
  \item[] 
    $\ZFC$ $+$
    ``\/$\exists\vsymb{\kappa}\,(\vsymb{\kappa}\mbox{ is }\Sigma_n\mbox{-correct and inaccessible\,})$''
  \end{xitemize}
  is consistent, then this theory does not prove 
  \begin{xitemize}
  \item[] 
    $\exists\vsymb{\kappa}\exists\vsymb{\lambda}\,(\,{}
    \begin{array}[t]{@{}l}
      \vsymb{\kappa}<\vsymb{\lambda}\ \land\ 
      \vsymb{\kappa}\mbox{ is }\Sigma_n\mbox{-correct and inaccessible}\\
      \land\ \vsymb{\lambda}\mbox{ is }\Sigma_{n-1}\mbox{-correct and inaccessible\,}
      ).
    \end{array}$ 
  \end{xitemize}
\end{Lemma}
\prf Suppose otherwise. Working in the theory \ZFC\ $+$ \memo{Goodman thesis: Lemma 2.1.2}
``$\exists\,\vsymb{\kappa}\,(\vsymb{\kappa}\mbox{ is }\Sigma_n\xmbox{-correct and inaccessible\,})$'', 
assume that $\kappa$ is $\Sigma_n$-correct and inaccessible and there is 
a $\Sigma_{n-1}$-correct inaccessible cardinal $\lambda>\kappa$.

\begin{Claim}\Label{Cl-Yah-0}
  $\uniV_\lambda\modelof{\kappa\mbox{ is }\Sigma_n\mbox{-correct and inaccessible}}$. 
\end{Claim}
\prfofClaim Suppose that $\psi(\ol{x},\ol{y})$ is $\Pi_{n-1}$ and $\ol{a}\in V_\kappa$.
If $V_\kappa\models\exists\ol{y}\psi(\ol{a},\ol{y})$, then there are $\ol{b}\in V_\kappa$ 
\st\ $V_\kappa\models\psi(\ol{a},\ol{b})$. Since $\psi$ is $\Pi_{n-1}$, it follows that
$V_\lambda\models\psi(\ol{a},\ol{b})$. Thus
$V_\lambda\models\exists\ol{y}\psi(\ol{a},\ol{y})$.

If $V_\lambda\models\exists\ol{y}\psi(\ol{a},\ol{y})$, then there are $\ol{y}\in V_\lambda$ 
\st\ $V_\lambda\models\psi(\ol{a},\ol{b})$. Ti follows 
that $\uniV\models\psi(\ol{a},\ol{b})$. 
Hence $\uniV\models\exists\ol{y}\psi(\ol{a},\ol{y})$. Now since $\kappa$ 
is $\Sigma_n$-correct, it follows that $V_\kappa\exists\ol{y}\psi(\ol{a},\ol{y})$. 
\qedofClaim\qedskip

Thus we proved $consis(\ZFC\,+\,
\exists\vsymb{\kappa}\,(\vsymb{\kappa}\mbox{ is }\Sigma_n\mbox{-correct and inaccessible\,}))$.
But this is a contradiction by The Second Incompleteness Theorem. 
\qedofProp\qedskip

\Propof{p-Yah-8} can be still improved as follows: 
\begin{Prop}\Label{p-Yah-9}
  Suppose $\calP$ is a $\Sigma_n$-definable iterable class of \pos\ \st\ $\calP$
  satisfies the condition \assertof{a'} or \assertof{b'} of \Lemmaof{p-Yah-4}, and it is also 
  provable that $\calP$ satisfies \xitemof{x-Yah-4} and \xitemof{x-Yah-5} of \Propof{p-Yah-5} for 
  any infinite cardinal $\lambda$ and inaccessible $\kappa\geq\lambda$. 
 
  \wassert{1} For $n\geq 3$, assuming the consistency of the theory
  \begin{xitemize}
  \xitem[x-Yah-12] 
    \ZFC\ $+$ ``\/$\continuum$ is regular'' $+$ $\MP^*(\calP,\calH(\kappa_\refl))_{\Sigma_n}$, 
  \end{xitemize}
  this theory does not prove  $\MP^*(\calP,\calH(\kappa_\refl))_{\Pi_{n}}$.\smallskip

  \wassert{1'} For $n\geq 3$, assuming the consistency of the theory
  \begin{xitemize}
  \xitem[x-Yah-12-a] 
    \ZFC\ $+$ ``\/$\continuum$ is regular'' $+$ $(\calP,\calH(\kappa_\refl))_{\Sigma_n}$-\RcAp, 
  \end{xitemize}
  this theory does not prove  $(\calP,\calH(\kappa_\refl))_{\Pi_{n}}$-\RcAp.\smallskip

  \wassert{2} Assume the consistency of the theory 
  \begin{xitemize}
  \xitem[x-Yah-12-0] \ZFC\ $+$ ``\/there is a supercompact cardinal and an 
    inaccessible above it.'' 
  \end{xitemize}
  Then $\MP^*(\calP,\calH(\kappa_\refl))_{\Sigma_2}$ does not 
  imply $\MP^*(\calP,\emptyset)_{\Pi_2}$. 
\end{Prop}
\prf \assertof{1}: 
Assume otherwise and suppose that the theory \xitemof{x-Yah-12} proves 
$\MP^*(\calP,\calH(\kappa_\refl))_{\Pi_{n}}$ for some $n\geq3$. 

By \Lemmaof{p-Yah-4}, in the theory \xitemof{x-Yah-12}, we have
$\uniL_{{\kappa_\refl}^V}\prec\uniL$. Thus from the assumption of the consistency of the 
theory \xitemof{x-Yah-12}, it follows the consistency of
\begin{xitemize}
\xitemx[] $\ZFC$ $+$ ``$\kappa_\refl$ is inaccessible in $\uniL$'' $+$
  $\uniL_{{\kappa_\refl}^V}\prec\uniL$. 
\end{xitemize}

By \Lemmaof{p-Yah-8}, we obtain the consistency of the following theory:
\begin{xitemize}
\xitem[x-Yah-13] \ZFC\ $+$ $\uniV=\uniL$ $+$\\[\jot]
  $\exists\vsymb{\kappa}\,({}
  \begin{array}[t]{@{}l}
    \mbox{``}\vsymb{\kappa}\mbox{ is inaccessible''}\,\land\,\uniL_{\vsymb{\kappa}}\prec_{\Sigma_n}\uniL\\[\jot]
    \land\,\ubecause{\mbox{``\/there is no }\Sigma_{n-1}\mbox{-correct inaccessible cardinals 
      above }\vsymb{\kappa}\mbox{''}}{}{\ixitemq[x-Yah-14]}).
  \end{array}$
\end{xitemize}
Working in this theory \xitemof{x-Yah-13}, we find a \po\ 
 $\poP_\kappa\in\calP$ \st\\
$\forces{\poP_\kappa}{\MP^*(\calP,\calH(\kappa_\refl))_{\Sigma_n}}$ by \Propof{p-Yah-5}. By 
\xitemof{x-Yah-14} and \Lemmaof{p-Yah-7},\,\assertof{1}
\\$\forces{\poP_\kappa}{\neg\MP^*(\calP,\calH(\kappa_\refl))_{\Pi_n}}$. This is a 
contradiction to the assumption we set at the beginning of the proof.\smallskip

\assertof{1'}: By \assertof{1} and \Propof{p-Yah-1},\,\assertof{1}.\smallskip

\assertof{2}: We work in the theory \ZFC\ $+$ ``\/there is a supercompact $\kappa$ and a 
single inaccessible $\mu$ above it''. The consistency of this theory follows form 
\xitemof{x-Yah-12-0}. In the following we shall use some notions and results from Goodman 
\cite{goodman} (see also the paragraph right before \Lemmaof{p-Yah-0} above).

The supercompact $\kappa$ is also supercompact for $C^{(1)}$ by Lemma 2.2.6 in 
\cite{goodman}. By Theorem 3.1.6 in \cite{goodman}, there is a \po\ $\poP\in\calP$ of size 
$\kappa$ \st\ $\uniV[\genG]\models\Sigma_2\mbox{-CFA}_{\LT\kappa}(\calP)$ for a
$(\uniV,\poP)$-generic $\genG$ 
(for the 
definition of this principle, see Definition 3.1.2 in \cite{goodman} --- this Theorem 3.1.6 
is proved similarly to our \Propof{p-Yah-5}). 
By modifying the construction of $\poP$ slightly if necessary, we also obtain
$\uniV[\genG]\models\kappa=\kappa_\refl$. 

By Theorem 3.1.4 in 
\cite{goodman}, this implies
$\uniV[\genG]\modelof{\MP^*(\calP,\calH(\kappa_\refl))_{\Sigma_2}}$. 

On the other hand, since $\mu$ is the unique inaccessible cardinal in $\uniV[\genG]$ above
$\kappa$, \Lemmaof{p-Yah-7},\,\assertof{2} implies
$\uniV[\genG]\modelof{\neg\MP^*(\calP,\emptyset)_{\Pi_2}}$. 
\qedofThm\qedskip

Some possible non-implications still remain. For example:

\begin{Problem}
  Is $\MP^*(\calP,\calH(\kappa_\refl))_{\Pi_2}$ $+$
  $\neg\MP^*(\calP,\calH(\kappa_\refl))_{\Sigma_2}$ consistent?
\end{Problem}
Cf.\ \Propof{p-Yah-1},\,\assertof{3}. \fi}

{\ifprivate\privatecolor\mbox{}\bigskip

%-------------------------------------------------------------------------------------------
  The following sections are sketches of the materials which will appear in a sequel to
  this paper:

\section{Universe with class many large cardinals}
\Label{class-many}

\memo{Woodin book: Theorem 2.53\\
 upward reflection under tightly generic hyperhuge cardinal}

%-------------------------------------------------------------------------------------------
\fi}\bigskip %%\ifprivate

\noindent
{\bf Acknowledgments}\quad
The authors would like to thank Gabriel Goldberg, Daisuke Ikegami, 
and Hiroshi Sakai for helpful remarks and comments. The authors also would like to thank Andreas 
Lietz for pointing out a flaw in the proof of an earlier version of 
\Thmof{p-hierarchies-6}. A part of the statement of the original version of this theorem is now 
recovered in \Thmof{p-hierarchies-6-1} and in the subsequent remark. They would like to 
thank as well the anonymous referee for pointing out numerous slips.  
\bigskip

\noindent
{\bf Competing Interests}\quad The research of the first author was supported by Kakenhi Grant-in-Aid for 
  Scientific Research (C) 20K03717. 

  For the second author, this research was funded in whole, or in part, by the Austrian Science
  Fund (FWF) [Grant numbers I6087]. For the purpose of open access, the
  authors have applied a CC BY public copyright license to any Author
  Accepted Manuscript version arising from this submission.

  The third author is an International Research Fellow of the Japan Society for the Promotion 
  of Science. 

%% \\The author would like to thank Kaethe Minden for making 
%%   him aware of \cite{minden}. He also would like to thank Gunter Fuchs, Takehiko Gappo, 
%%   Paul Larson, Hiroshi Sakai, and Kostas Tsaprounis for many helpful remarks and comments.

%-------------------------------------------------------------------------------------------


\begin{thebibliography}{99}
\addcontentsline{toc}{section}{References}
\Label{ref}%
\newcommand{\bysame}[1]{\underline{\phantom{#1}}}%
\newcommand{\bysamex}[1]{#1}%
\iftesting
\newcommand{\Bibitem}[1]{\bibitem{#1}\marginnote{{\color{cyan}%
\renewcommand{\baselinestretch}{0.4}\tiny \rlap{#1}}}}
\else
\newcommand{\Bibitem}[1]{\bibitem{#1}}
\fi
%%
{\ifprivate\privatecolor
\Bibitem{abraham-shelah} Uri Abraham, and Saharon Shelah, Forcing closed unbounded sets, 
  Journal of Symbolic Logic,  Vol.48(3), (1983), 643--657.\quad

\fi}

%% \Bibitem{alexandroff-uryson}  Paul Alexandroff and Paul Urysohn, Zur Theorie der 
%%   topologischen R\"aume, Mathematische Annalen 92 (1924), 258--266. 


\Bibitem{apter-Laver-Indestr} Arthur W.\ Apter, Laver indestructibility and the class of compact 
  cardinals, The Journal of Symbolic Logic, Vol.63, No.1, (1998), 149--157.\quad

\Bibitem{aspero} David Asperó, A Maximal bounded forcing axiom, The Journal of Symbolic Logic, 
  Vol.67, No.1 (2002), 130--142.\quad

{\ifextended\extendedcolor
\Bibitem{aspero-schindler} David Asperó, and Ralf Schindler, Martin's Maximum
  implies Woodin's axiom $(*)$, Annals of Mathematics, Vol.193 (2021), 793--835.\quad

  \fi}


\Bibitem{bagaria-bounded} Joan Bagaria, Bounded forcing axioms as principles of generic 
  absoluteness, Archive for Mathematical Logic, Vol.39, (2000), 393--401.\quad


\Bibitem{bagaria-Cn} \bysame{Joan Bagaria}, $C^{(n)}$-cardinals, Archive for Mathematical 
  Logic,  Vol.51, (2012), 213--240.\quad


\Bibitem{BHTU} Joan Bagaria, Joel David Hamkins,  Konstantinos Tsaprounis, Toshimichi Usuba,
  Superstrong and other large cardinals are never Laver indestructible, 
  Archive for Mathematical Logic, Vol.55 (2016), 19--35. 

\Bibitem{5a} Neil Barton, Andrés Eduardo Caicedo, Gunter Fuchs, Joel David Hamkins, Jonas 
  Reitz, and Ralf Schindler, Inner-Model Reflection Principles, Studia Logica, 108 (2020), 
  573--595.\quad

%% \Bibitem{BMR} J.\ Baumgartner, J.\ Malitz, and W.\ Reinhardt, Embedding Trees in the Rationals, 
%%  Proceedings of the National Academy of Sciences, Vol.67, No.4, (1970), 1748-1753. 
%% %% pp. 1748-1753, December 1970

%% \Bibitem{baumgartner-taylor-wagon} 
%%   James E.\,Baumgartner, Alan D.\,Taylor and Stanley Wagon,
%%   On Splitting Stationary Subsets of Large Cardinals, 
%%   The Journal of Symbolic Logic,  Vol.42, No.2 (1977), 203--214. 

{\ifprivate\privatecolor
\Bibitem{boney} Will Boney, Model Theoretic Characterizations of Large Cardinals,
  Israel Journal of Mathematics, 236, (2020), 133--181.\quad

\fi}
  
%% \Bibitem{cohenPE}  Paul {\color{red} E.} Cohen, A Large Power Set Axiom, The Journal of Symbolic 
%% Logic, Vol.40, No.1, (1975), 48--54.

{\ifprivate\privatecolor
\Bibitem{cox} Sean Cox, The diagonal reflection principle, 
Proceedings of the American Mathematical Society, Vol.140, No.8 (2012), 
2893--2902.\quad

\fi}

\Bibitem{cox2} Sean D.\,Cox, Forcing axioms, approachability, and stationary set 
  reflection, The Journal of Symbolic Logic Volume 86, Number 2, June 2021, 499--530.\quad  


{\ifprivate\privatecolor
\Bibitem{dtw}  Alan Dow, Franklin D.\ Tall and W.A.R.\ Weiss, New proofs of the consistency 
of  
  the normal Moore space conjecture I, II, Topology and its Applications, 37 (1990) 33--51, 
  115-129.\quad

\fi}

\Bibitem{eda-kada-yuasa} Katsuya Eda, Masaru Kada, and Yoshifumi Yuasa, The tightness about 
  sequential fans and combinatorial properties, Journal of the Mathematical Society of Japan, 
  Vol.49, No.1, (1997), 181--187.\quad


%%{\ifprivate\privatecolor
\Bibitem{friedman-sy} Sy-David Friedman, Internal Consistency and the Inner Model 
  Hypothesis, The Bulletin of Symbolic Logic, Vol.\,12, No.\,4 (Dec., 2006), 591--600.\quad

%%\fi}

\Bibitem{potential} Sakaé Fuchino, On potential embedding and versions of Martin's axiom, 
Notre Dame Journal of Logic, Vol.33, No.4, (1992), 481--492.\quad


{\ifprivate\privatecolor
\Bibitem{fuchino-left} \bysame{Sakaé Fuchino}, Left-separated topological spaces under 
  Fodor-type Reflection Principle, RIMS Kôkyûroku, No.1619,
(2008), 32--42.\quad

\fi}

{\ifprivate\privatecolor
\Bibitem{balogh} \bysame{Sakaé Fuchino}, Fodor-type Reflection Principle and Balogh's 
  reflection theorems, RIMS Kôkyûroku, No.1686, (2010), 41--58.\quad

\fi}

{\ifprivate\privatecolor
\Bibitem{sf-note} \bysame{Sakaé Fuchino}, Rado's Conjecture implies
the Fodor-type Reflection Principle, Note, \\
\href{https://fuchino.ddo.jp/notes/RCimpliesFRP2.pdf}{\scalebox{0.89}[1]{\small\tt 
  https://fuchino.ddo.jp/notes/RCimpliesFRP2.pdf}}\quad

\fi}

{\ifprivate\privatecolor
\Bibitem{RIMS15} \bysame{Sakaé Fuchino}, On reflection numbers under large 
  continuum, RIMS Kôkyû\-roku, No.1988, 1--16, (2016).\quad

\fi}

{\ifprivate\privatecolor
\Bibitem{pre-hilbert} \bysame{Sakaé Fuchino}, Pre-Hilbert spaces without orthonormal 
  bases, preprint.\quad

\fi}
  
{\ifprivate\privatecolor
\Bibitem{fuchino-reverse} \bysame{Sakaé Fuchino}, A reflection principle as a 
  reverse-mathematical  
  fixed point over the base theory ZFC, Annals of the Japan Association for the Philosophy 
  of Science,  
  Vol.25, (2017), 67--77.\quad

\fi}

{\ifprivate\privatecolor
\Bibitem{nagoya} \bysame{Sakaé Fuchino}, Images of the white board of a seminar talk given 
  at Nagoya University on May 31, 2019,\\
  \href{https://fuchino.ddo.jp/talks/talk-nagoya-2019-05-31.pdf}{%
    \scalebox{0.89}[1]{\small\tt 
      https://fuchino.ddo.jp/talks/talk-nagoya-2019-05-31.pdf}}\quad

\fi}

{\ifprivate\privatecolor
\Bibitem{slides} \bysame{Sakaé Fuchino}, Resurrection and Maximality under a/the tightly 
  Laver-generically ultrahuge cardinal, slides of a series zoom talks given at May 29, 
  and June 5, 2023 at Kobe Set Theory Seminar.\\
  \href{https://fuchino.ddo.jp/slides/kobe2023-05-29-pf.pdf}{\scalebox{0.89}[1]{\small\tt 
    https://fuchino.ddo.jp/slides/kobe2023-05-29-pf.pdf}}\quad 

\fi}

{\ifprivate\privatecolor
\Bibitem{slides-x} \bysame{Sakaé Fuchino}, Resurrection and Maximality under a/the tightly 
  Laver-generically ultrahuge cardinal --- additional slides, additional slides of the zoom 
  talks given at June 5, 2023 at Kobe Set Theory Seminar.\\ 
  \href{https://fuchino.ddo.jp/slides/kobe2023-06-05a-pf.pdf}{\scalebox{0.89}[1]{\small\tt 
    https://fuchino.ddo.jp/slides/kobe2023-06-05a-pf.pdf}}\quad 

\fi}

\Bibitem{future} \bysame{Sakaé Fuchino}, 
Maximality Principles and Resurrection Axioms  
  in light of a Laver-generic large cardinal, 
  Pre-preprint.\\   \href{https://fuchino.ddo.jp/papers/RIMS2022-RA-MP-x.pdf}{%
    \scalebox{0.89}[1]{\small\tt 
    https://fuchino.ddo.jp/papers/RIMS2022-RA-MP-x.pdf}}\quad 

\Bibitem{janos} \bysame{Sakaé Fuchino}, Reflection and Recurrence, Model Theory, Computer 
  Science, and Graph Polynomials 
  Festschrift in Honor of Johann A. Makowsky (2025).\\
\href{https://fuchino.ddo.jp/papers/reflection\_and\_recurrence-Janos-Festschrift-x.pdf}{%
\scalebox{0.89}[1]{\small\tt https://fuchino.ddo.jp/papers/reflection\_and\_recurrence-Janos-Festschrift-x.pdf}}\quad


\Bibitem{future2} Takehiko Gappo, and Andreas Lietz, 
  Separating Maximality Principles, preprint.

{\ifprivate\privatecolor
\Bibitem{fjetal} Sakaé Fuchino, István Juhász, Lajos Soukup,
Zoltán Szentmikl\'ossy, and Toshimichi Usuba, Fodor-type Reflection Principle and 
reflection of metrizability and
meta-Lindelöfness, Topology and its Applications, Vol.157, 8 (2010), 1415--1429.\quad

\fi}

{\ifprivate\privatecolor
\Bibitem{sfetal-I} Sakaé 
  Fuchino, André Ottenbreit Maschio Rodrigues, and Hiroshi Sakai, 
  Strong Löwenheim-Skolem theorems for stationary logics, I, 
  Archive for Mathematical Logic,  Volume 60, issue 1-2, (2021), 17--47.\\
  \href{https://fuchino.ddo.jp/papers/SDLS-x.pdf}{\scalebox{0.89}[1]{%%
      \small\tt https://fuchino.ddo.jp/papers/SDLS-x.pdf}}\quad

\fi}  

\Bibitem{sfetal-II} \ifprivate\bysame{Sakaé 
  Fuchino, Maschio Rodrigues and Hiroshi, Sakai}\else\bysamex{Sakaé 
  Fuchino, André Ottenbreit Maschio Rodrigues, and Hiroshi Sakai}\fi, Strong downward 
  Löwenheim-Skolem theorems for stationary logics, II --- reflection down  
  to the continuum,
  Archive for Mathematical Logic,  Volume 60, issue 3-4, (2021), 495--523.\\
  \href{https://fuchino.ddo.jp/papers/SDLS-II-x.pdf}{\scalebox{0.89}[1]{\small\tt 
      https://fuchino.ddo.jp/papers/SDLS-II-x.pdf}}\quad

{\ifextended\extendedcolor
\Bibitem{sfetal-III} \bysame{Sakaé 
  Fuchino, Maschio Rodrigues, and Hiroshi Sakai}, Strong downward Löwenheim-Skolem 
  theorems for stationary logics, III --- mixed support iteration, 
  to appear in the Proceedings of the Asian Logic Conference 2019.\\
  \href{https://fuchino.ddo.jp/papers/SDLS-III-xx.pdf}{\scalebox{0.89}[1]{\small\tt 
      https://fuchino.ddo.jp/papers/SDLS-III-xx.pdf}}\quad 

\fi}

{\ifextended\extendedcolor
\Bibitem{FuOt} Sakaé Fuchino, André Ottenbreit Maschio Rodrigues, Reflection 
  principles, generic large cardinals, and the Continuum Problem, 
  Proceedings of the Symposium on Advances in Mathematical Logic 2018, Springer (2021), 1--26.\\
  \href{https://fuchino.ddo.jp/papers/refl\_principles\_gen_large\_cardinals\_continuum\_problem-x.pdf}{\scalebox{0.79}[1]{\small\tt 
https://fuchino.ddo.jp/papers/refl\_principles\_gen\_large\_cardinals\_continuum\_problem-x.pdf}}\quad

\fi}

{\ifprivate\privatecolor
\Bibitem{rinot} Sakaé Fuchino and Assaf Rinot, Openly generated Boolean algebras and 
  the Fodor-type Reflection Principle, Fundamenta Mathematicae 212, (2011), 261-283.\quad

\fi}

\Bibitem{fuchino-sakai} Fuchino, Sakaé, and Sakai, Hiroshi, The first-order definability 
  of generic large cardinals, submitted.  Extended version of the paper:\\  
  \href{https://fuchino.ddo.jp/papers/definability-of-glc-x.pdf}{\scalebox{0.89}[1]{\small\tt 
      https://fuchino.ddo.jp/papers/definability-of-glc-x.pdf}}\quad 


{\ifextended\extendedcolor
\Bibitem{fuchino-sakai-2} Sakaé Fuchino, and Hiroshi Sakai, Generically
  super\-com\-pact cardinals by forcing with chain 
  conditions, RIMS Kôkûroku,\\No.2213, (2022), 94--111. \\ 
  \href{https://repository.kulib.kyoto-u.ac.jp/dspace/bitstream/2433/275453/1/2213-08.pdf}{%%
    \scalebox{0.85}[1]{\small\tt
      https://repository.kulib.kyoto-u.ac.jp/dspace/bitstream/2433/275453/1/2213-08.pdf}}\quad

\fi}
  %% https://fuchino.ddo.jp/papers/RIMS2021-ccc-gen-supercompact-x.pdf

{\ifprivate\privatecolor
\Bibitem{more} Sakaé Fuchino, Hiroshi Sakai, Lajos Soukup and Toshimichi Usuba, 
  More about Fodor-type Reflection Principle, \\
\href{https://fuchino.ddo.jp/papers/moreFRP-x.pdf}{\scalebox{0.89}[1]{\small\tt 
  https://fuchino.ddo.jp/papers/moreFRP-x.pdf}}\quad

\fi}


{\ifprivate\privatecolor
\Bibitem{fuchino-torres-etal} Sakaé Fuchino, Hiroshi Sakai, Victor Torres Perez, and 
  Toshimichi Usuba, Rado's 
  Conjecture and the Fodor-type Reflection Principle, in preparation.\quad

\fi}

\Bibitem{recurrence} Sakaé Fuchino, and Toshimichi Usuba, On recurrence axioms,
  Annals of Pure and Applied Logic, Vol.176, (10), (2025).\\
  \href{https://fuchino.ddo.jp/papers/recurrence-axioms-x.pdf}{\scalebox{0.89}[1]{\small\tt 
  https://fuchino.ddo.jp/papers/recurrence-axioms-x.pdf}}\quad


{\ifprivate\privatecolor
\Bibitem{fuchs} Gunter Fuchs Closed Maximality Principles: Implications, Separations and 
  Combinations, The Journal of Symbolic Logic, Vol.\,73, No.\,1, (2008), 276--308.\quad

\fi}

%%{\ifextended\extendedcolor
\Bibitem{goldberg} Gabriel Goldberg, Usuba's theorem is optimal, Journal of Mathematical 
  Logic, (2024).

%%\fi}

\Bibitem{goldstern-shelah} Martin Goldstern, and Saharon Shelah, The Bounded Proper Forcing 
Axiom, The Journal of Symbolic Logic, Vol.60, No.1 (1995), 58--73.\quad


%%{\ifextended\extendedcolor
\Bibitem{goodman} Benjamin P.\ Goodman, Sigma$_{\,}$\_$_{\,}$n-correct Forcing Axioms, CUNY Academic 
  Works (2024).\\
  \href{https://academicworks.cuny.edu/gc\_etds/5808/}{\scalebox{0.89}[1]{\small\tt 
    https://academicworks.cuny.edu/gc$_{\,}$\_$_{\,}$etds/5808/}}\quad

%%\fi}

{\ifextended\extendedcolor
\Bibitem{fhr} Gunter Fuchs, Joel David Hamkins, and Jonas Reitz, Set-theoretic geology, 
  Annals of Pure and Applied Logic 166(4), (2015), 464--501. 

\fi}

\Bibitem{hamkins}  Joel David Hamkins, A simple maximality principle, The Journal of 
Symbolic Logic Vol.68, No.7, (2003), 527--550.\quad 


\Bibitem{hamkins-johnstone}  Joel David Hamkins, and Thomas A.\,Johnstone, Resurrection 
  axioms and uplifting cardinals, Archive for Mathematical Logic, Vol.53, Iss.3-4, (2014), 
  463--485.\quad


{\ifprivate\privatecolor
\Bibitem{hamkins-johnstone2} Joel David Hamkins, and Thomas A.\,Johnstone, Strongly 
  uplifting cardinals and the boldface resurrection axioms, Archive for Mathematical Logic 
  Vol.\,56, (2017), 1115--1133.\quad

\fi}

%% \Bibitem{halbeisen} Lorenz J.\,Halbeisen, Combinatorial Set Theory, With a Gentle 
%%   Introduction to Forcing, Springer (2012).\quad


%% \Bibitem{hausdorff}  Felix Hausdorff, Grundz\"uge der Mengenlehre, Veit\& Comp., Leipzig (1914). 

%% \Bibitem{herrlich}  Horst Herrlich, Axiom of Choice, Lecture Notes in Mathematics, 1876, 
%%   Springer (2006).\quad

\Bibitem{ikegami-trang} Daisuke Ikegami and Nam Trang, On a class of maximality principles, 
  Archive for Mathematical Logic, Vol.\,57, (2018), 713--725.\quad

%% \Bibitem{jech}  Thomas Jech, The Axiom of Choice, Dover Publications, Inc., (1973/2008). 

\ifextended
\Bibitem{millennium-book} Thomas Jech, Set Theory, The Third Millennium
	Edition, Springer, (2001/2003).\quad

\else
\Bibitem{millennium-book} Thomas Jech, Set Theory, The Third Millennium
	Edition, Springer, (2003).\quad

\fi
\Bibitem{higher-inf} Akihiro Kanamori, The Higher Infinite, 
  Springer--Verlag (1994/2003).\quad

%% \Bibitem{kelley} John L.\ Kelley, Tychonoff's theorem implies AC, Fundamenta Mathematicae 37, 
%%   (1950), 75--76.

{\ifprivate\privatecolor
\Bibitem{koenig} Bernhard König, Generic compactness reformulated, Archive for  
  Mathematical Logic 43, (2004), 311--326.\quad

\fi}

\Bibitem{kunen-2011} Kenneth Kunen, Set Theory, College Publications (2011).\quad

%% \Bibitem{los-ryll-nardzewski} Jerzy \L o\'s and Czes\l aw Ryll-Nardzewski, Effectiveness of 
%%   the representation theory for Boolean algebras, Fundamenta Mathematicae 41, (1955), 
%%   49--56. 

{\ifprivate\privatecolor
\Bibitem{larson-lumsdaine-yin} Paul B.\,Larson, Peter Lumsdaine and Yimu Yin, An 
  introduction to $\poP_{\rm max}$ forcing, in: James Cummings, and Ernest Schimmerling 
  (eds.): Appalachian Set Theory 2006--2012, Cambridge University Press (2013), 5--23.\quad

\fi}

\Bibitem{laver-indestr} Richard Laver, Making the supercompactness of $\kappa$ indestructible 
  under $\kappa$-directed closed forcing, Israel Journal of Mathematics, Vol.29, No.4 (1978), 
  385--388.\quad

\Bibitem{levy} Azriel Lévy, A Hierarchy of Formulas in Set Theory. Memoirs of the American 
  Mathematical Society 57 (1965).\quad

\Bibitem{reitz}  Jonas Reitz, The Ground Axiom, The Journal of Symbolic Logic, Vol.\,72, 
  No.\,4 (2007), 1299--1317.\quad

{\ifprivate\privatecolor
\Bibitem{minden} Kaethe Minden, Combining resurrection and maximality, The Journal of 
  Symbolic Logic, Vol.\,86, No.\,1, %March
  (2021), 397--414.\quad

\fi}

%% \Bibitem{schindler} Ralf Schindler, Semi-proper forcing, remarkable cardinals, and Bounded Martin's 
%%   Maximum, Mathematical Logic Quarterly, Vol.50, (6), (2004), 527--532.\quad 

\Bibitem{schindler2} Ralf Schindler, Bounded Martin's Maximum and Strong Cardinals, in:
  Joan Bagaria, Stevo Todorcevic (Eds.), Set Theory
Centre de Recerca Matemàtica Barcelona, 2003-2004, Birkhäuser, (2006), 401--406. 



%% \Bibitem{proper-improper} Saharon Shelah, Proper and Improper Forcing, Second edition, 
%%   Springer--Verlag (1998). 

\Bibitem{stavi-vaananen} Jonathan Stavi, and Jouko Väänänen, Reflection principles for the 
  continuum, Logic and algebra, 59--84, Contemporary Mathematics, 302, American 
  Mathematical Society, Providence, RI, (2002). 
  
%% \Bibitem{tarski} Alfred Tarski, Sur quelques theorems qui equivalent a l'axiome du choix, 
%%   Fundamenta Mathematicae, 5, (1924), 147–154.  

%% \Bibitem{stevo1983} Stevo Todor\v{c}evi\'c, On a conjecture of R.Rado, 
%%   Journal London Mathematical Society
%%   Vol.s2-27, (1) (1983), 1--8. 

{\ifprivate\privatecolor
\Bibitem{stevo1993} Stevo Todor\v{c}evi\'c, Conjectures of Rado and Chang and cardinal 
  arithmetic, In: Sauer, N.W., Woodrow, R.E., Sands, B. (eds.) Finite and Infinite 
  Combinatorics in Sets and Logic. NATO ASI Series, vol 411, (1993), 385–398. 
  
\fi}


%% \Bibitem{tychonoff} A.\ Tychonoff, \"Uber die topologische Erweitrerung von R\"aumen, 
%%   Mathematische Annalen 102 (1930), 544--561.

{\ifextended\extendedcolor
\Bibitem{tsaprounis1} Konstantinos Tsaprounis, On resurrection axioms, 
	The Journal of Symbolic Logic, Vol.80, No.2, %%June
  (2015), 587--608.\quad

\fi}

{\ifprivate\privatecolor
\Bibitem{tsaprounis2} \bysame{Konstantinos Tsaprounis}, Ultrahuge cardinals, Mathematical 
  Logic Quarterly, Vol.62, No.1-2, (2016), 1--2.\quad

\fi}

{\ifextended\extendedcolor
\Bibitem{usuba} Toshimichi Usuba, The downward directed grounds hypothesis 
and very large cardinals, Journal of Mathematical Logic, Vol.\,17(2) (2017), 1--24.\quad

\fi}

{\ifextended\extendedcolor
\Bibitem{usuba2} \bysame{Toshimichi Usuba}, Extendible cardinals and the mantle, Archive for 
  Mathematical Logic, Vol.58, (2019), 71-75.\quad

\fi}

%% \Bibitem{viale} Matteo Viale, The Proper Forcing Axiom and the Singular Cardinal
%% Hypothesis, The Journal of Symbolic Logic, Vol.71,\,(2), (2006), 473--479. \quad

\Bibitem{venturi-viale} Giorgio Venturi, and  Matteo Viale, Second order arithmetic as the 
  model companion of set theory, Archive for Mathematical Logic Vol.62, (2023), 
  29--53.\quad 


\Bibitem{viale-revisited} Matteo Viale, Martin's maximum revisited, Archive for Mathematical 
  Logic, Vol.55, (2016), 295--317.\quad


%% \Bibitem{wiki-heine-borel}  Wikipedia, Heine-Borel Theorem, \\
%%   \href{https://en.wikipedia.org/wiki/Heine-Borel_theorem}{https://en.wikipedia.org/wiki/Heine-Borel\_theorem}. 

\Bibitem{woodin-book} Hugh Woodin, The Axiom of Determinacy, Forcing Axioms, and the 
  Nonstationary Ideal, De Gruyter, 2nd rev. ed. Edition (2010). 


%% \Bibitem{zermelo1904} Ernst Zermelo, Beweis, da\ss\ jede Menge wohlgeordnet werden kann 
%%   (Aus einem an Herrn Hilbert gerichteten Briefe), Mathematische Annalen, Vol.\,59, (1904), 514--516.





\ifprivate{\extendedcolor
\bigskip\bigskip\mbox{}

\mbox{}\hspace{-1.5em}{\large\bf Additional references cited in the extended version}

{\ifprivate\privatecolor
\Bibitem{math-20} 渕野 昌，数学ノート (2020 -), \\
  \href{https://fuchino.ddo.jp/notes/math-notes-20.pdf}{\scalebox{0.89}[1]{\small\tt 
    https://fuchino.ddo.jp/notes/math-notes-20.pdf}}\quad

\fi}

%% \Bibitem{MP-RA} Sakaé Fuchino, Maximality Principles and Resurrection Axioms in 
%%   light of Laver-generic large cardinals, in preparation. 

%% \Bibitem{FS}   S.\ Fuchino, H.\ Sakai, The first-order definability of generic large cardinals, pre-preprint.\\
%% \href{https://fuchino.ddo.jp/papers/definability-of-glc-x.pdf}{https://fuchino.ddo.jp/papers/definability-of-glc-x.pdf}%




%% {\ifprivate\privatecolor
%% \Bibitem{rietz} Jonas Reitz, The ground axiom, Journal of Symbolic Logic Vol.72, No.4, (2007), 
%%   1299--1317. 
%% 
%% \fi}

%% \Bibitem{sh794} Saharon Shelah, Reflection implies the SCH,  Fundamenta Mathematicae 198, 
%%   (2008), 95--110. 



%% \Bibitem{takeuti}  Gaishi Takeuti, Hypotheses on power set, Proceedings of Symposia in Pure Mathematics,
%% Vol.13, Part I, American Mathematical Society, Providence, R.I., (1971), 439--446.
%% \Bibitem{takeuti-article} \bysame{Gaishi Takeuti} (竹内外史)，ランダム実数と連続体仮説，数学セミナー，\\
%% 1999年 5月号，(1999), 34--37.
}\fi




\end{thebibliography}
\end{document}
%%% Local Variables: %%%
%%% mode:LaTeX %%%
%%% my-local-language: "non-Japanese" %%%
%%% tex-main-file: "generic-absoluteness-revisited-xx.tex" %%%
%%% End: %%%